\documentclass[11pt]{article}

\usepackage{color}
\usepackage{latexsym}
\usepackage{amssymb}
\usepackage{amsmath, amsfonts,amssymb,theorem,euscript,array,enumerate,amsfonts,mathrsfs}
\usepackage{graphicx}

\newtheorem{Theorem}{Theorem}[part]
\newtheorem{Definition}{Definition}[part]
\newtheorem{Proposition}{Proposition}[part]

\newtheorem{Lemma}{Lemma}[part]

\newtheorem{Remark}{Remark}[part]
\newtheorem{Example}{Example}[part]

\def\esssup_#1{\underset{#1}{\mathrm{ess\,sup\, }}}
\def\essinf_#1{\underset{#1}{\mathrm{ess\,inf\, }}}

\def \trans{^{\scriptscriptstyle{\intercal}}}

\def\<{\langle }
\def\>{\rangle }

\def \Frac{\displaystyle\frac}

\def \trans{^{\scriptscriptstyle{\intercal }}}

\def \N{\mathbb{N}}
\def \R{\mathbb{R}}

\def \E{\mathbb{E}}
\def \F{\mathbb{F}}

\def \P{\mathbb{P}}

\def \Ac{{\cal A}}
\def \Bc{{\cal B}}

\def \Ec{{\cal E}}
\def \Fc{{\cal F}}

\def \Lc{{\cal L}}
\def \Pc{{\cal P}}
\def \Mc{{\cal M}}

\def \Vc{{\cal V}}

\def \Vc{{\cal V}}

\def \eps{\varepsilon}

\def \ep{\hbox{ }\hfill$\Box$}

\def\Dt#1{\Frac{\partial #1}{\partial t}}
\def\DT#1{\Frac{\partial #1}{\partial T}}

\def\reff#1{{\rm(\ref{#1})}}

\def\beqs{\begin{eqnarray*}}
\def\enqs{\end{eqnarray*}}
\def\beq{\begin{eqnarray}}
\def\enq{\end{eqnarray}}

\allowdisplaybreaks

\begin{document}

\title{Long time asymptotics for fully nonlinear Bellman equations:\\ a Backward SDE approach}

\author{Andrea COSSO\thanks{Laboratoire de Probabilit\'es et Mod\`eles Al\'eatoires, CNRS, UMR 7599, Universit\'e Paris Diderot,
		\sf  cosso at math.univ-paris-diderot.fr}~~~
		Marco FUHRMAN\thanks{Dipartimento di Matematica, Politecnico di Milano,
               \sf  marco.fuhrman at polimi.it}~~~
               Huy{\^e}n PHAM\thanks{Laboratoire de Probabilit\'es et Mod\`eles Al\'eatoires, CNRS, UMR 7599, Universit{\'e} Paris Diderot, and
               CREST-ENSAE,  \sf pham at math.univ-paris-diderot.fr}
             }

\date{\today}

\maketitle

\begin{abstract}
We study the large time behavior of solutions to fully nonlinear parabolic equations of Hamilton-Jacobi-Bellman type arising typically in stochastic control theory with control both  on drift and diffusion coefficients. 
We prove that, as time horizon goes to infi\-nity, the long run average solution is characterized by a nonlinear ergodic equation. Our results hold under dissipativity conditions, and without any nondegeneracy assumption on the diffusion term. Our approach uses mainly probabilistic arguments relying on new backward SDE representation for  nonlinear parabolic, elliptic and ergodic equations. 
\end{abstract}

\vspace{5mm}

\noindent {\bf MSC Classification}: 60H30, 60J60, 49J20.

\vspace{5mm}

\noindent {\bf Keywords}: Hamilton-Jacobi-Bellman type equation, stochastic control, large time be\-havior, backward SDE, ergodic type Bellman equation.

\newpage

\section{Introduction}

Let us consider the fully nonlinear parabolic equation of the form
\beq \label{GHJBP} 
\DT{v}  -  \sup_{a \in A}\big[ \Lc^a v + f\big(x,a,\frac{v}{T+1}\big) \big] & =& 0, \;\;\; \mbox{ on } \;  (0,\infty)\times\R^d,
\enq
with $v(0,.)$ $=$ $h$ on $\R^d$, where  $\Lc^a$ is the second-order differential operator:
\beqs
\Lc^a v &=& b(x,a).D_x v + \frac{1}{2}{\rm tr}(\sigma\sigma\trans(x,a)D_x^2 v).
\enqs
Here $A$ is some Borel subset of $\R^q$, $b,\sigma$ are continuous functions on $\R^d\times\R^q$, and  $f$ $=$ $f(x,a,y)$ 
is a measurable function on $\R^d\times\R^q\times\R$ satisfying some  conditions to be specified later on.   When the generator $f$ $=$ $f(x,a)$ does not depend on $y$,  equation \reff{GHJBP} is the dynamic programming equation, also called Hamilton-Jacobi-Bellman (HJB) equation, associated to the stochastic control problem:
\beqs
v(T,x) & := & \sup_{\alpha\in\Ac} \E_x\Big[ \int_0^T f(X_t^{\alpha},\alpha_t) dt + h(X_T^\alpha) \Big], 
\enqs 
where $X^{\alpha}$ is the controlled diffusion process 
 \beq \label{controlX}
 dX_t^{\alpha} &=&   b(X_t^{\alpha},\alpha_t) dt +   \sigma(X_t^{\alpha},\alpha_t) dW_t,
 \enq
driven by a $d$-dimensional Brownian motion $W$ on a probability space $(\Omega,\Fc,\P)$ equipped with the natural filtration of $W$, and  given a control  $\alpha$ $\in$ $\Ac$, i.e., 
an $A$-valued adapted process.  In the general case $f$ $=$ $f(x,a,y)$, we shall see that under suitable conditions, there exists a unique viscosity solution $v$ $=$ $v(T,x)$ 
to the generalized parabolic HJB equation \reff{GHJBP}, and our aim is to investigate the large time behavior of $v(T,.)$ as $T$ goes to infinity.  It turns out that this asymptotic problem is related to the 
generalized {\it ergodic}  HJB equation: 
\beq \label{GEHJB}
\lambda - \sup_{a \in A} \big[ \Lc^a \phi + f(x,a,\lambda) \big] & =& 0, \;\;\; \mbox{ on } \; \R^d. 
\enq

Asymptotics for stochastic control and related HJB equation have been studied in va\-rious settings by many authors since the works  \cite{benfre92}  and \cite{arilio98}.  In the PDE literature, we refer for instance to \cite{barsou01} in a periodic setting, \cite{souzha06} under Dirichlet conditions, or \cite{fujishlor06} in the whole space. In these cited  papers, the HJB equation is semi-linear, i.e., the nonlinearity   
appears only in the first order derivative. Recently, by combining PDE and stochastic analysis arguments, the papers  \cite{ichihara12}, \cite{ichshe13} and \cite{robxin13} studied large time behavior of semi-linear HJB equations with quadratic nonlinearity in gradients. We would like also to point out the recent paper \cite{humadric14}, which studied large time behavior of solutions to semi-linear HJB equations by a probabilistic approach relying on ergodic BSDE introduced in \cite{fuhrman_hu_tess09}. Interestingly, the authors are able to prove in their context a rate of convergence for the solution to the parabolic equation towards the ergodic equation under weak dissipativity conditions. Long time asymptotics of solutions to HJB equations has been also considered in the context of risk-sensitive stochastic control and utility maximization problem, see e.g. \cite{fleshe99} and  \cite{nag03}.

Our motivation is to develop a systematic study applicable to a large class of fully nonlinear HJB type equation, and to give natural conditions on the dynamics of the control system ensuring ergodicity.  
The principal novelty of this paper is to consider control on both drift and diffusion coefficients $b(x,a)$, $\sigma(x,a)$, possibly degenerate, and satisfying  dissipativity conditions,  instead of periodicity condition.  
In this case, we do not have in general smooth solution to the HJB equation.  Another  original  feature of our framework is the dependence of   $f(x,a,y)$ on $y$, which occurs for example in stochastic control with recursive utility functions.   
Our first main result is to prove the existence of a viscosity solution pair $(\lambda,\phi)$ $\in$ $\R\times C(\R^d)$, with $\phi$ Lipschitz, to the ergodic fully nonlinear equation \reff{GEHJB}.   We adopt the following approach.  We consider the sequence of fully nonlinear elliptic HJB equation for $\beta$ $>$ $0$: 
\beq \label{GHJBE} 
\beta v^\beta  -  \sup_{a \in A}\big[ \Lc^a v^\beta + f\big(x,a,\beta v^\beta\big) \big] & =& 0, \;\;\; \mbox{ on } \;  \R^d,  
\enq
and obtain the existence and uniqueness of a solution $v^\beta$ to \reff{GHJBE} by combining analytical and probabilistic methods. More precisely, following the randomization approach of \cite{khapha12} for representing parabolic HJB equations,  we introduce a class of Backward Stochastic Diffe\-rential Equations (BSDEs) with nonpositive jumps over an infinite horizon, and supported by a  forward regime switching process $(X,I)$ where 
\beq \label{regimeX}
dX_t &=& b(X_t,I_t) dt + \sigma(X_t,I_t) dW_t,
\enq
and $I$ is a pure jump process valued in $A$.  The minimal solution $Y^\beta$ 
to this class of elliptic BSDEs is shown to exist and to provide the unique  (viscosity) solution $v^\beta$ to \reff{GHJBE}, and  the key point is to derive 
uniform Lipschitz estimate for the sequence  $(v^\beta)_\beta$. This is achieved by ergodicity properties on the forward process $X$, and suitable estimation on the minimal solution $Y^\beta$. 
Then, by standard analytical approximation procedures (when $\beta$ goes to zero) as in \cite{ichihara12} or \cite{fuhrman_hu_tess09},  we obtain the existence of a pair $(\lambda,\phi)$ solution to \reff{GEHJB}.  Moreover, the function $\phi$ admits a probabilistic representation in terms of a new class of BSDEs, namely ergodic BSDEs with nonpositive jumps. Ergodic BSDEs  have been introduced in   \cite{fuhrman_hu_tess09} and then in \cite{ric09},  and related to optimal control on the drift of diffusions.  We extend this connection to the case of controlled diffusion coefficient by imposing a nonpositive jump constraint on the ergodic BSDE.

Next, our main theorem  is to prove that for any solution $(\lambda,\phi)$ to \reff{GEHJB}, we have the convergence  of the solution to the parabolic generalized HJB equation \reff{GHJBP}: 
\beq \label{VTlambda}
\frac{v(T,.)}{T} & \longrightarrow & \lambda, \;\;\;\; \mbox{ in } \; C(\R^d), \;\; \mbox{ as } \; T \rightarrow \infty. 
\enq
Here, convergence  ``in $C(\R^d)$" stands for locally uniform convergence in $\R^d$.  This shows as a  byproduct that $\lambda$ in \reff{GEHJB} is unique.  
The main difficulty with respect to the semi-linear HJB case is that we do not have in general a smooth solution and an optimal control for the finite horizon and ergodic stochastic control, and the classical arguments as 
in \cite{ichihara12} or \cite{ichshe13} do not work anymore. Moreover, when $f$ $=$ $f(x,a,y)$ depends also on $y$, we do not even have a stochastic control representation of the function $v$. Our arguments for proving  \reff{VTlambda} rely on the BSDE representation of solution to  \reff{GHJBP} and \reff{GEHJB}, corresponding comparison theorems, and dual representation of such BSDEs in terms of equivalent probability measures introduced in \eqref{Pnu}. 
Furthermore, we can strengthen the  convergence result \reff{VTlambda} by a verification theorem:  under the  condition that the  ergodic equation \reff{GEHJB} admits a classical component solution  $\phi$, we have
\beq \label{VTPHI}
v(T,.)  - \big( \lambda T + \phi) & \longrightarrow & c, \;\;\;\; \mbox{ in } \; C(\R^d), \;\; \mbox{ as } \; T \rightarrow \infty,
\enq
for some constant $c$.

The paper is organized as follows. Section 2  introduces some notations, formulates the dissipativity conditions on $b,\sigma$, and assumptions on $f$.   We then  state ergo\-dicity properties on the regime switching process $(X,I)$ in \reff{regimeX} as well as on the controlled diffusion process $X^\alpha$ in \reff{controlX}.  In Section 3, we  prove the existence and uniqueness of a solution to the fully nonlinear elliptic  HJB equation \reff{GHJBE} and its relation to BSDE with nonpositive jumps over an infinite horizon. Section 4 is concerned with the ergodic equation \reff{GEHJB} and its probabilistic representation in terms of ergodic BSDE with nonpositive jumps.  Convergence results \reff{VTlambda} and \reff{VTPHI} are studied in Section 5.  We collect in the Appendix some proofs and technical estimates needed in the paper.

\section{Ergodicity properties}

\setcounter{equation}{0} 
\setcounter{Theorem}{0} \setcounter{Proposition}{0}
\setcounter{Corollary}{0} \setcounter{Lemma}{0}
\setcounter{Definition}{0} \setcounter{Remark}{0}

\subsection{Notations and  assumptions}

Let $(\Omega,\Fc,\P)$ be a complete probability space on which are defined  a $d$-dimensional
Brownian motion $W$ $=$ $(W_t)_{t\geq 0}$  and an independent Poisson random measure $\mu$ on $\R_+\times A$, where $A$ is a compact   subset of $\R^q$, endowed with its Borel $\sigma$-field $\Bc(A)$.  We assume that the random measure $\mu$ has the intensity measure $\vartheta(da)dt$ for some finite measure $\vartheta$  on $(A,\Bc(A))$.  We set $\tilde\mu(dt,da)$ $=$ $\mu(dt,da)-\vartheta(da)dt$ the compensated martingale measure
associated to $\mu$, and denote by $\F$ $=$ $(\Fc_t)_{t\geq 0}$ the completion of the natural filtration generated by $W$ and $\mu$. We also denote, for any $T>0$, $\Pc_T$ the $\sigma$-field  of $\F$-predictable subsets of $[0,T]\times\Omega$. Let us introduce some additional notations. We denote by:
\begin{itemize}
\item ${\bf L^p(}\Fc_t{\bf)}$, $p$ $\geq$ $1$, $t\geq0$, the set of $\Fc_t$-measurable random variables $X$ such that $\E[|X|^p]$ $<$ $\infty$.
\item ${\bf S^2(t,T)}$, $0\leq t<T<\infty$, the set  of real-valued c\`adl\`ag adapted processes $Y$ $=$
$(Y_s)_{t\leq s\leq T}$ satisfying
\[
\E\Big[\sup_{t\leq s\leq T}|Y_s|^2\Big] \ < \ \infty.
\]
We also define ${\bf S_{\text{loc}}^2}:=\cap_{T>0}{\bf S^2(0,T)}$.
\item ${\bf L^p(W;t,T)}$, $p$ $\geq$ $1$, $0\leq t<T<\infty$, the set of  $\R^d$-valued predictable processes
$Z=(Z_s)_{t\leq s\leq T}$ such that
\[
\E\bigg[\bigg(\int_t^T |Z_s|^2 ds\bigg)^{\frac{p}{2}}\bigg] <\infty.
\]
We also define ${\bf L_{\text{loc}}^p(W)}:=\cap_{T>0}{\bf L^p(W;0,T)}$.
\item ${\bf L^p(\tilde\mu;t,T)}$, $p$ $\geq$ $1$, $0\leq t<T<\infty$, the set of
$\Pc_T\otimes\Bc(A)$-measurable maps $U\colon[t,T]\times\Omega\times A\rightarrow \R$ such that
\[
\E\bigg[ \int_t^T\bigg(\int_A  |U_s(a)|^2 \vartheta(da)\bigg)^{\frac{p}{2}}ds\bigg] \ < \ \infty.
\]
We also define ${\bf L_{\text{loc}}^p(\tilde\mu)}:=\cap_{T>0}{\bf L^p(\tilde\mu;0,T)}$.
\item ${\bf K^2(t,T)}$, $0\leq t<T<\infty$, the  set of  nondecreasing c\`adl\`ag predictable processes $K$ $=$ $(K_s)_{t\leq s\leq T}$ such that $\E[|K_T|^2]$ $<$ $\infty$ and $K_t$ $=$ $0$. We also define ${\bf K_{\text{loc}}^2}:=\cap_{T>0}{\bf K^2(0,T)}$.
\end{itemize}

We are given some  continuous functions $b\colon\R^d\times \R^q \rightarrow \R^d$, $\sigma\colon\R^d\times \R^q \rightarrow \R^{d\times d}$, and 
consider the forward regime switching process $(X,I)$ governed by the stochastic differential equation in $\R^d\times\R^q$:
\begin{equation}
\label{forward}
\begin{cases}
dX_t= b(X_t,I_t)\,dt +   \sigma(X_t,I_t)\,dW_t, \\
dI_t=\int_A(a-I_{t^-})\,\mu(dt,da).
\end{cases}
\end{equation}
We note that the fact that $\sigma$ is a square matrix does not involve any loss of generality, since we are not going to assume any nondegeneracy condition.  In particular, some rows or columns of 
$\sigma$ may be equal to zero.  
In the following we use the notation $M\trans$ for the transpose of any matrix $M$, and $\|M\|^2$ $=$ $\text{tr}(MM\trans)$ for the Hilbert-Schmidt norm. We shall make  the following assumptions on the coefficients $b$ and $\sigma$. 

\vspace{2mm}

{\bf (H1)} 
\begin{itemize}
\item[(i)] There exists a positive constant $L_1$  such that for all $x,x'\in\R^d$, $a,a'\in\R^q$,
\beqs
|b(x,a)-b(x',a')| + \|\sigma(x,a)-\sigma(x',a')\|  &\le&   L_1 \,(|x-x'|+|a-a'|). 
\enqs
\item[(ii)] There exists a constant $\gamma>0$ such that for all  $x,x'\in\R^d$, $a\in A$,
\beq \label{dissipative}
(x-x').(b(x,a)-b(x',a)) + \frac12 \|\sigma(x,a)-\sigma(x',a)\|^2 & \le &  \ -\gamma \,|x-x'|^2.
\enq
\end{itemize} 
It is well-known that under {\bf (H1)}(i), there exists a unique solution $(X_t^{x,a},I_t^a)_{t\geq0}$ to \eqref{forward} starting from $(x,a)\in\R^d\times\R^q$ at time $t=0$.  Notice that when $a$ $\in$ $A$, then  $I_t^a\in A$ for all $t$ $\geq$ $0$.   
Condition {\bf (H1)}(ii) is called {\it dissipativity condition} and will ensure the ergodicity of the process $X$, as stated  in the next paragraph. 


\begin{Example}
{\rm  Let  $b(x,a)$ $=$ $B(a)x + D(a)$, $\sigma(x,a)$ $=$ $\Sigma(a)$ for some vector valued Lipschitz function $D$, and matrix valued Lipschitz functions $B$, $\Sigma$ on $A$, such that $B$ is uniformly stable:
\beqs
x.B(a)x & \leq & - \gamma |x|^2, \;\;\; \forall x \in \R^d, \;  a \in A.
\enqs
In this case, {\bf (H1)} is satisfied, and this example corresponds to a controlled Ornstein-Uhlenbeck process with uncertain mean-reversion and volatility. 
\ep
}
\end{Example}

\vspace{1mm}

We also consider some real-valued continuous function $f$ on $\R^d\times\R^q\times\R$ satisfying the following assumption:

\vspace{2mm}

{\bf (H2)} 
\begin{itemize}
\item[(i)] There exists a positive constant $L_2$  such that for all $x,x'\in\R^d$, $a,a'\in \R^q$, $y,y'$ $\in$ $\R$,
\beqs
|f(x,a,y)-f(x',a',y')|   &\le&   L_2 \,(|x-x'|+ |a-a'| + |y-y'|). 
\enqs
\item[(ii)] The function $y$ $\in$ $\R$ $\longmapsto$ $f(x,a,y)$ is nonincreasing for all $(x,a)$ $\in$ $\R^d\times \R^q$.
\end{itemize}

\vspace{1mm}

We end this paragraph of notations by  introducing the following set  of probability measures, which shall play an important role in the sequel for establishing estimates. Let $\Vc_n$ be the set of 
$\Pc\otimes\Bc(A)$-measurable maps  valued in $[1,n+1]$,  $\Vc$ $=$ $\cup_{n\in\N} \Vc_n$,  
and consider for $\nu$ $\in$ $\Vc$, the probability measure $\P^\nu$ equivalent to 
$\P$ on $(\Omega,\Fc_T)$, for any $T$ $>$ $0$, with Radon-Nikodym density: 
\beq \label{Pnu}
\frac{d\P^\nu}{d\P}\Big|_{\Fc_t} &=& \zeta_t^\nu \; := \; \Ec_t \Big( \int_0^. \int_A (\nu_s(a)-1) \tilde\mu(ds,da) \Big), 
\enq
for $0\leq t\leq T$, where $\Ec(.)$ is the Dol\'eans-Dade exponential. Actually, since $\nu$ $\in$ $\Vc$ is essentially bounded, it is shown in Lemma 2.4  in \cite{khapha12} that $(\zeta^\nu)_{0\leq t\leq T}$ is a uniformly integrable $\P$-martingale, with $\zeta_T^\nu$ $\in$ ${\bf L^2}(\Fc_T)$,  for any $T$ $>$ $0$,  and so it  defines a probability measure $\P^\nu$ via \reff{Pnu}. 
We shall denote by $\E^\nu$ the expectation under $\P^\nu$.  Moreover, by Girsanov theorem,   the compensator of $\mu$ under $\P^\nu$  is $\nu_t(a)\vartheta(da)dt$, while $W$ remains a Brownian motion independent of $\mu$ under $\P^\nu$. We denote by $\tilde\mu^\nu(dt,da)$ $=$ $\mu(dt,da)-\nu_t(a)\vartheta(da)dt$ the compensated martingale measure of $\mu$ under $\P^\nu$.

\subsection{Ergodicity}

We now use the dissipativity condition in {\bf (H1)}(ii) to state moment estimates and stability results on the component  $X$ of  \reff{forward}.

\begin{Lemma} \label{estimX}
Let Assumption {\bf (H1)} hold.

\noindent \textup{(i)}  There exists a positive constant $C$ $=$ $C_{b,\sigma}$ depending only on $b,\sigma$  such that 
for all  $x\in\R^d$ and $a\in A$, 
\beq \label{estimate_x^2}
\sup_{t\geq0} \sup_{\nu \in \Vc} \E^\nu\big[|X_t^{x,a}|^2\big]  &\le&   C (1+|x|^2),
\enq

\noindent \textup{(ii)} For all $t$ $\geq$ $0$, $x,x'$ $\in$ $\R^d$, $a$ $\in$ $A$, 
\beq \label{ergoX}
\sup_{\nu\in\Vc} \E^\nu\big[|X_t^{x,a}-X_t^{x',a}|^2\big] & \leq & |x-x'|^2e^{-2\gamma t}. 
\enq
\end{Lemma}

The proof relies on rather standard arguments based on It\^o's formula and Gronwall's lemma, and is reported in the Appendix.

\begin{Remark}
{\rm
We shall need the following generalization of estimate \eqref{estimate_x^2}, for all $x\in\R^d$ and $a\in A$,
\beq \label{estimate_x^2_conditional}
\sup_{s\geq t} \sup_{\nu \in \Vc} \E^\nu\big[|X_s^{x,a}|^2\big|\Fc_t\big]  &\le&   C (1+|x|^2),
\enq
which is valid with the same constant $C=C_{b,\sigma}$, independent of $t$, as in \eqref{estimate_x^2}, when Assumption {\bf (H1)} holds.
\ep
}
\end{Remark}

\vspace{3mm}

Let $\underline{\alpha}\colon\R^d\rightarrow A$  be a feedback control and let $X$ $=$ $X^{\underline{\alpha}}$ be the associated diffusion process governed by 
\beq \label{xmarkov}
dX_t &=&  b(X_t,\underline{\alpha}(X_t))\,dt + \sigma(X_t,\underline{\alpha}(X_t))\,dW_t.
\enq
Suppose that the functions
\begin{equation}
\label{underline_b_sigma}
\underline b(x) \ := \ b(x,\underline\alpha(x)), \qquad \underline\sigma(x) \ := \ \sigma(x,\underline\alpha(x))
\end{equation}
are Lipschitz.  Then,  equation \eqref{xmarkov} defines a time-homogeneous Markov process $\{X_t^{\underline{\alpha}}\,,\,t\geq 0\}$, and we denote by $(P_t^{\underline{\alpha}})_{t\geq0}$ the associated semigroup, which acts  on  $B(\R^d)$, the set of bounded measurable functions $\varphi$, by 
\beqs
P_t^{\underline{\alpha}}\varphi(x) & = & \E_x\big[\varphi(X_t^{\underline{\alpha}} )\big], \qquad t \geq 0, \; x\in\R^d.
\enqs
Notice that $(P_t^{\underline{\alpha}})_{t\geq0}$ has the Feller property, i.e., for any $f\in C_b(\R^d)$, the space of continuous and bounded functions on $\R^d$, we have that 
$P_t^{\underline{\alpha}} f\in C_b(\R^d)$.  The next result shows the ergodicity  of $X^{\underline{\alpha}}$.

\begin{Proposition}
\label{P:LawInvariant} 
Let $\underline{\alpha}\colon\R^d\rightarrow A$ be  a feedback control such that $\underline b,\underline\sigma$ in \eqref{underline_b_sigma} satisfy Assumption {\bf (H1)}.  Then $X^{\underline{\alpha}}$ is ergodic, i.e., the following assertions are valid: 
\begin{itemize}
\item[\textup{(i)}] There exists a unique invariant probability measure $\rho$ $=$ $\rho^{\underline{\alpha}}$ on $\R^d$$:$ 
\beqs
\int P_t^{\underline{\alpha}} \varphi(x) \rho(dx)  &=& \int \varphi(x) \rho(dx),  \;\;\;  \forall \; t \geq 0, \;  \varphi \in B(\R^d). 
\enqs
\item[\textup{(ii)}] $X_t^{\underline{\alpha}}$  converges weakly to  $\rho$ as $t\to\infty$$:$ 
\beq \label{psiergodic}
P_t^{\underline{\alpha}}\varphi(x)  & \longrightarrow & \int \varphi(x) \rho(dx), \;\;\; \mbox{ as } \; t \rightarrow \infty, \;\;\; \forall\,x \in \R^d,\,\varphi\in C_b(\R^d). 
\enq
\end{itemize}
Moreover,  $\int |x|^2\rho(dx)<\infty$, and the convergence \reff{psiergodic}  holds for all continuous $\varphi$ satisfying a linear growth condition.
\end{Proposition}

The proof is based on the ``pullback'' method (see, e.g., Theorem 6.3.2 in \cite{daprato_zabczyk96}) and is detailed in the Appendix.

\section{Elliptic HJB equation}

\setcounter{equation}{0} 
\setcounter{Theorem}{0} \setcounter{Proposition}{0}
\setcounter{Corollary}{0} \setcounter{Lemma}{0}
\setcounter{Definition}{0} \setcounter{Remark}{0}

For any $\beta$ $>$ $0$, let us  consider the  fully nonlinear elliptic equation of HJB type:
\beq \label{ellipticbetapde}
\beta \,v^{\beta} - \sup_{a\in A}\big[\Lc^av^{\beta} + f(x,a,\beta v^\beta)\big] & = &  0, \;\;\; \mbox{ on } \; \R^d,
\enq
where
\beqs
\Lc^a\varphi  &= &  b(x,a). D_x\varphi  + \frac{1}{2}\text{tr}\big(\sigma\sigma\trans(x,a)D_x^2\varphi\big). 
\enqs
 
Notice that in the particular case where $f$ $=$ $f(x,a)$ does not depend on $y$, the equation \reff{ellipticbetapde} is the dynamic programming equation associated to the stochastic control problem on infinite horizon:
\beqs
v^\beta(x) & := & \sup_{\alpha\in\Ac} \E_x \Big[ \int_0^\infty e^{-\beta t} f(X_t^\alpha,\alpha_t) dt \Big].  
\enqs
In this case, it is easy to see from the Lipschitz and growth condition on $f$ in {\bf (H2)}(i), and the estimates \reff{estimate_x^2}, \reff{ergoX} that 
the sequence of functions  $(v^\beta)_\beta$ satisfies the uniform estimates:
\beq \label{estimvbeta}
v^\beta(x) \; \leq \; \frac{C}{\beta}(1 + |x|),  & &  |v^\beta(x)-v^\beta(x')| \; \leq \; C |x-x'|, \;\;\; \forall x,x' \in \R^d, 
\enq
for some positive constant $C$ independent of $\beta$. 

In the general case $f$ $=$ $f(x,a,y)$, this section is devoted to the existence and uniqueness of a viscosity solution $v^\beta$ to \reff{ellipticbetapde}, and to uniform  estimate on $(v^\beta)_\beta$ as in \reff{estimvbeta}. 
To this purpose, we introduce  the following  class of BSDE with nonpositive jumps over an infinite horizon, for any $\beta$ $>$ $0$: 
\begin{align}
\label{BSDE}
Y_t^\beta \ &= \ Y_T^\beta - \beta\int_t^T Y_s^\beta ds + \int_t^T f(X_s,I_s,\beta Y_s^\beta) ds + K_T^\beta-K_t^\beta \notag \\
&\quad \ - \int_t^T Z_s^\beta dW_s - \int_t^T \int_A U_s^\beta(a)\tilde\mu(ds,da), \qquad 0 \leq t \leq T, \;  \forall T \in (0,\infty),  
\end{align}
and
\begin{equation}
\label{JumpCon}
U_t^\beta(a) \ \leq \ 0, \qquad dt\otimes d\P\otimes\vartheta(da)\text{-a.e.}
\end{equation}

BSDEs driven by Brownian motion over an infinite horizon  have been introduced in \cite{darling95}, \cite{darpardoux97}, studied also in \cite{brihu98} and extended in \cite{roy04}, and related to elliptic semi-linear PDEs.  
Here, we extend this definition to BSDEs driven by Brownian motion and Poisson random measure, and with the nonpositivity constraint on the jump component.  
A \emph{minimal solution} to the elliptic BSDE with nonpositive jumps \eqref{BSDE}-\eqref{JumpCon}  is a quadruple $(Y^\beta,Z^\beta,U^\beta,K^\beta)\in\mathbf{S_{\textup{loc}}^2}\times\mathbf{L_{\text{loc}}^2(W)}\times\mathbf{L_{\text{loc}}^2(\tilde\mu)}\times\mathbf{K_{\text{loc}}^2}$ satisfying  \eqref{BSDE}-\eqref{JumpCon}, with $|Y_t^\beta|\leq C(1+|X_t|)$, for all $t\geq0$ and for some constant $C$, such that for any other solution 
$(\bar Y^\beta,\bar Z^\beta,\bar U^\beta, \bar K^\beta)\in\mathbf{S_{\textup{loc}}^2}\times\mathbf{L_{\text{loc}}^2(W)}\times\mathbf{L_{\text{loc}}^2(\tilde\mu)}\times\mathbf{K_{\text{loc}}^2}$ to \eqref{BSDE}-\eqref{JumpCon}, satisfying $|\bar Y_t^\beta|\leq C'(1+|X_t|)$ for some constant $C'$, we have $Y_t^\beta \leq \bar Y_t^\beta$, $\P$-a.s., for all  $t\geq0$.

\begin{Remark}
\label{P:Uniqueness}
{\rm There exists at most one minimal solution to \eqref{BSDE}-\eqref{JumpCon}. Indeed, let $(Y,Z,U,K)$ and $(\tilde Y,\tilde Z,\tilde U,\tilde K)$ be two minimal solutions to \eqref{BSDE}-\eqref{JumpCon}. The uniqueness of the $Y$ component is clear by definition. Regarding the other components, taking the difference between the two backward equations we obtain
\begin{equation}
\label{BSDEdifference}
\int_0^t \big(Z_s - \tilde Z_s\big)dW_s \ = \ K_t - \tilde K_t - \int_0^t\int_A \big(U_s(a) - \tilde U_s(a)\big)\tilde\mu(ds,da),
\end{equation}
for all $0\leq t\leq T$, $\P$-almost surely. Then, we see that the right-hand side is a finite variation process, while the left-hand side has not finite variation, unless $Z = \tilde Z$. Now, from \eqref{BSDEdifference}, we obtain the identity
\[
\int_0^t\int_A \big(U_s(a) - \tilde U_s(a)\big)\mu(ds,da) \ = \ \int_0^t\int_A \big(U_s(a) - \tilde U_s(a)\big)\vartheta(da)ds + K_t - \tilde K_t,
\]
where the right-hand side is a predictable process, therefore it has no totally inaccessible jumps (see, e.g., Proposition 2.24, Chapter I, in \cite{jacod_shiryaev03}); on the other hand, the left-hand side is a pure-jump process with totally inaccessible jumps, unless $U=\tilde U$. As a consequence, we must have $U=\tilde U$, from which it follows that $K=\tilde K$.
\ep
}
\end{Remark}

\vspace{3mm}

In the sequel, we prove by a penalization approach  the existence of  the minimal solution to \eqref{BSDE}-\eqref{JumpCon}, which shall provide  the solution to the elliptic nonlinear HJB equation \reff{ellipticbetapde}.  Then, by using this probabilistic representation of $v^\beta$, we shall state uniform Lipschitz estimate on $(v^\beta)_\beta$.

\subsection{Elliptic penalized  BSDE}

For any $\beta>0$ and $n\in\N$, we consider the penalized BSDE on $[0,\infty)$, $\P$-a.s.,
\begin{align}
Y_t^{\beta,n} \ &= \ Y_T^{\beta,n} - \beta \int_t^T Y_s^{\beta,n}\,ds + \int_t^T f(X_s,I_s,\beta Y_s^{\beta,n})\,ds + n  \int_t^T \int_A (U_s^{\beta,n}(a))_+\vartheta(da)\,ds \notag \\
&\quad \ - \int_t^TZ_s^{\beta,n}\,dW_s - \int_t^T\int_AU_s^{\beta,n}(a)\,\tilde\mu (ds,da), \qquad 0\leq t\leq T<\infty, \label{penalizedbetabsde}
\end{align}
where $h_+=\max(h,0)$ denotes the positive part of the function $h$.

\vspace{2mm}

We first state an a priori estimate on the above elliptic penalized BSDE. 

\begin{Lemma} \label{estimapriori}
Suppose that \textup{\bf (H1)} holds. Let $(x,a),(x',a')\in\R^d\times A$ and $(Y^{1,\beta,n},Z^{1,\beta,n},U^{1,\beta,n})$ $($resp. $(Y^{2,\beta,n},Z^{2,\beta,n},U^{2,\beta,n})$$)$ be a solution in $\mathbf{S_{\textup{loc}}^2}\times\mathbf{L_{\textup{loc}}^2(W)}\times\mathbf{L_{\textup{loc}}^2(\tilde\mu)}$ to \eqref{penalizedbetabsde}, with $(X,I)=(X^{x,a},I^a)$ $($resp. $(X,I)=(X^{x',a'},I^{a'})$$)$ and $f=f_1$ $($resp. $f=f_2$$)$ satisfying assumption \textup{\bf (H2)}. Set $\Delta_t Y$ $=$ $Y_t^{1,\beta,n}-Y_t^{2,\beta,n}$, $\Delta_t Z$ $=$ $Z_t^{1,\beta,n}-Z_t^{2,\beta,n}$, $\Delta_t U(a'')$ $=$ $U_t^{1,\beta,n}(a'')-U_t^{2,\beta,n}(a'')$,  
$\Delta_t'  f_1$ $=$ $f_1(X_t^{x,a},I_t^a,\beta Y_t^{1,\beta,n})$ $-$ $f_1(X_t^{x',a'},I_t^{a'},\beta Y_t^{1,\beta,n})$, and  $\Delta_t f$ $=$ 
$f_1(X_t^{x',a'},I_t^{a'},\beta Y_t^{2,\beta,n})-f_2(X_t^{x',a'},I_t^{a'},\beta Y_t^{2,\beta,n})$, $t\geq 0$, $a''\in A$. Then, there exists 
$\nu$ $\in$ $\Vc_n$ such that for all $T$ $\in$ $(0,\infty)$,
\beq
|\Delta_t Y|^2 & \leq &  \E^\nu \bigg[ e^{-2\beta(T-t)} |\Delta_T Y|^2 + 2 \int_t^T e^{-2\beta(s-t)} \Delta_s Y 
(\Delta_s' f_1 + \Delta_s f) ds \bigg| \Fc_t \bigg], \label{EstimateDeltaYZU}
\enq
for all $0\leq t\leq T$.
\end{Lemma} 
{\bf Proof.} See Appendix.
\ep

\vspace{2mm}

The next result states the existence and uniqueness of a solution to  \eqref{penalizedbetabsde}, and uniform estimate on the solution.

\begin{Proposition}
\label{P:ExistUniqPen}
Let Assumptions {\bf (H1)} and {\bf (H2)}  hold. Then, for any $(x,a,\beta,n)\in\R^d\times A\times(0,\infty)\times\N$,  there exists a solution $(Y^{x,a,\beta,n},Z^{x,a,\beta,n},U^{x,a,\beta,n})\in\mathbf{S_{\textup{loc}}^2}\times\mathbf{L_{\textup{loc}}^2(W)}\times\mathbf{L_{\textup{loc}}^2(\tilde\mu)}$ to \eqref{penalizedbetabsde}, with $(X,I)$ $=$ $(X^{x,a},I^a)$, and satisfying:
\beq \label{estimYbetan}
|Y_t^{x,a,\beta,n}| &\leq& \frac{C_{b,\sigma,f}}{\beta}(1+|X_t^{x,a}| ), \;\;\; \forall t \geq 0,
\enq
for some positive constant  $C_{b,\sigma,f}$ depending only on $b,\sigma,f$. Moreover, this solution is unique in the class of triplets $(Y,Z,U)\in\mathbf{S_{\textup{loc}}^2}\times\mathbf{L_{\textup{loc}}^2(W)}\times\mathbf{L_{\textup{loc}}^2(\tilde\mu)}$ satisfying the condition $|Y_t|\leq C(1+|X_t^{x,a}|)$, for all $t\geq0$ and for some positive constant $C$ (possibly depending on $x$, $a$, $\beta$ and $n$).
\end{Proposition}
\textbf{Proof.} See Appendix.
\ep

\vspace{3mm}

For any $(x,a,\beta,n)\in\R^d\times A\times(0,\infty)\times\N$, we notice that $Y_0^{x,a,\beta,n}$ is a constant since it is $\Fc_0$-measurable. Therefore, for each $\beta$ $>$ $0$, $n$ $\in$ $\N$, we introduce the function $v^{\beta,n}\colon\R^d\times A\rightarrow\R$ defined as
\begin{equation}
\label{E:v^beta,n}
v^{\beta,n}(x,a) \ := \ Y_0^{x,a,\beta,n}, \qquad  \,(x,a)\in\R^d\times A.
\end{equation}

Let us now investigate some  key properties of the function $v^{\beta,n}$. We first state a uniform Lipschitz estimate on $(v^{\beta,n})$.

\begin{Lemma}
\label{L:vbeta,n_Bdd_Lip}
Let Assumptions {\bf (H1)} and {\bf (H2)}  hold.  For any $(\beta,n)\in(0,\infty)\times\N$, the function $v^{\beta,n}$ is such that: 
$Y^{x,a,\beta,n}_t=v^{\beta,n}(X^{x,a}_t,I^{a}_t)$, for all $t\geq0$, and  $(x,a)\in\R^d\times A$. 
Moreover, there exists some positive constant $C$ depending only on $b,\sigma,f$, and independent of $\beta,n$  such that 
\beq
v^{\beta,n}(x,a) & \leq &  \frac{C}{\beta}\big( 1 + |x| \big), \\
\big|v^{\beta,n}(x,a)-v^{\beta,n}(x',a)\big| & \leq & C  |x-x'|, \label{v^beta,n_Lipschitz}
\enq
for all $x,x'\in\R^d$ and $a\in A$. 
\end{Lemma}
\textbf{Proof.} See Appendix. 
\ep

\vspace{3mm}

As expected, for fixed $(\beta,n)$,  the function $v^{\beta,n}$ is related to the elliptic integro-differential equation:
\beq
\beta \,v^{\beta,n}(x,a) - \Lc^av^{\beta,n}(x,a) - \Mc^av^{\beta,n}(x,a) - f\big(x,a,\beta v^{\beta,n}(x,a)\big) &  &  \label{ellipticpdepenalizedbeta} \\
- n\int_A(v^{\beta,n}(x,a') -v^{\beta,n}(x,a))_+\,\vartheta(da') &  = & 0, \;\; \mbox{ on } \; \R^d\times A, \nonumber 
\enq
where
\[
\Mc^a\varphi(a) \ = \ \int_A \big(\varphi(a')-\varphi(a)\big) \vartheta(da'),
\]
for any $\varphi\in C(A)$. More precisely, we have the following result.

\begin{Proposition}
\label{P:ViscPropv^beta,n}
Let Assumptions {\bf (H1)} and {\bf (H2)}  hold. Then, the function $v^{\beta,n}$ in \eqref{E:v^beta,n} is a continuous viscosity solution to \eqref{ellipticpdepenalizedbeta}, i.e., it is continuous on 
$\R^d\times A$ and it is a viscosity supersolution $($resp. subsolution$)$ to \eqref{ellipticpdepenalizedbeta}, namely
\beqs
\beta \,\varphi(x,a) \ &\geq \ (resp.\text{ }\leq) &  \Lc^a\varphi(x,a) + \Mc^a\varphi(x,a) + f\big(x,a,\beta\varphi(x,a)\big) \\
&  & \quad \ + \; n  \int_A(\varphi(x,a') -\varphi(x,a))_+\,\vartheta(da')
\enqs
for any $(x,a)\in\R^d\times A$ and any $\varphi\in C^2(\R^d\times\R^q)$ such that
\beqs
0 \;\, = \;\, (v^{\beta,n}-\varphi)(x,a) & = &  \min_{\R^d\times A}(v^{\beta,n}-\varphi) \quad (resp.\text{ }\max_{\R^d\times A}(v^{\beta,n}-\varphi)). 
\enqs
\end{Proposition}
\textbf{Proof.} See Appendix. 
\ep

\subsection{Elliptic BSDE with nonpositive jumps}

We can now prove the existence of the minimal solution to the elliptic BSDE with nonpo\-sitive jumps  \eqref{BSDE}-\eqref{JumpCon}.

\begin{Proposition}
\label{T:Exist}
Let Assumptions {\bf (H1)} and {\bf (H2)}  hold.  Then, for any $\beta>0$ and $(x,a)\in\R^d\times A$ there exists a solution $(Y^{x,a,\beta},Z^{x,a,\beta},U^{x,a,\beta},K^{x,a,\beta})\in\mathbf{S_{\textup{loc}}^2}\times\mathbf{L_{\textup{loc}}^2(W)}\times\mathbf{L_{\textup{loc}}^2(\tilde\mu)}\times\mathbf{K_{\textup{loc}}^2}$ to \eqref{BSDE}-\eqref{JumpCon}, with $(X,I)=(X^{x,a},I^a)$. Moreover
\begin{enumerate}
\item[\textup{(i)}] $Y^{x,a,\beta}$ is the increasing limit of $(Y^{x,a,\beta,n})_n$ and satisfies 
\beq \label{growthYbeta}
|Y_t^{x,a,\beta}| &\leq& \frac{C}{\beta}(1+|X_t^{x,a}|), \;\;\; \forall  t\geq0,
\enq
for some positive constant $C$ independent of $\beta,x,a,t$. 
\item[\textup{(ii)}] $(Z^{x,a,\beta}_{|[0,T]},U^{x,a,\beta}_{|[0,T]})$, for any $T>0$, is the strong $($resp. weak$)$ limit of $(Z^{x,a,\beta,n}_{|[0,T]},U^{x,a,\beta,n}_{|[0,T]})_n$  in ${\bf L^p(W;0,T)}\times{\bf L^p(\tilde \mu;0,T)}$, with $p\in[1,2)$, $($resp. in ${\bf L^2(W;0,T)}\times{\bf L^2(\tilde \mu;0,T)}$$)$.
\item[\textup{(iii)}] $K_t^{x,a,\beta}$ is the weak limit of $(K_t^{x,a,\beta,n})_n$ in ${\bf L^2}(\Fc_t)$, for any $t\geq0$.
\end{enumerate}
Furthermore, this solution is minimal in the class of quadruplets $(Y,Z,U,K)\in\mathbf{S_{\textup{loc}}^2}\times\mathbf{L_{\textup{loc}}^2(W)}\times\mathbf{L_{\textup{loc}}^2(\tilde\mu)}\times\mathbf{K_{\textup{loc}}^2}$ satisfying the condition $|Y_t|\leq C(1+|X_t^{x,a}|)$, for all $t\geq0$ and for some positive constant $C$ (possibly depending on $x$, $a$, and $\beta$).
\end{Proposition}
\textbf{Proof.} See Appendix. 
\ep

\vspace{3mm}

For any $\beta>0$, let us introduce the deterministic function $v^\beta\colon\R^d\times A\rightarrow\R$ defined as follows
\begin{equation}
\label{vbeta}
v^\beta(x,a) \ := \ Y_0^{x,a,\beta}, \qquad \forall\,(x,a)\in\R^d\times A.
\end{equation}
From point (i) of Proposition  \ref{T:Exist}, it follows that $(v^{\beta,n})_n$ converges increasingly to $v^\beta$. Then, the identification $Y_t^{x,a,\beta,n}=v^{\beta,n}(X_t^{x,a},I_t^a)$ implies that $Y_t^{x,a,\beta}=v^\beta(X_t^{x,a},I_t^a)$. We shall now investigate the relation between $v^\beta$ and the fully nonlinear elliptic PDE of HJB type \eqref{ellipticbetapde}. More precisely, we shall prove that $v^\beta$ solves in the viscosity sense equation \eqref{ellipticbetapde}. The main issue is to prove that $v^\beta$ does not depend actually on $a$. However, as we do not know a priori that the function $v^\beta$ is continuous in both arguments, we shall rely on discontinuous viscosity solutions arguments as in \cite{khapha12},
and make the following assumptions on the set $A$ and the intensity measure $\vartheta$:

\vspace{3mm}

\begin{tabular}{cl}
\hspace{-2mm}{\bf (H$A$)}  &\hspace{-2mm} The interior $\mathring{A}$ of $A$ is connected, and $A$ $=$ Cl$(\mathring{A})$, the closure of its interior. \\
& \\
\hspace{-2mm}{\bf (H$\vartheta$)} &\hspace{-2mm} The measure $\vartheta$ supports the whole set $\mathring{A}$, and the boundary of $A$: $\partial A$ $=$ $A\backslash\mathring{A}$,\\
&\hspace{-2mm}  is negligible with respect to $\vartheta$. 
\end{tabular}

\vspace{2mm}

Notice that equation \reff{ellipticbetapde} does not depend on $\vartheta$, and this intensity measure only appears in order to give a probabilistic representation of $v^\beta$. Therefore, we have the choice to fix  an intensity measure $\vartheta$ satisfying condition {\bf (H$\vartheta$)}, which is anyway a fairly general condition easy to  realize.  In the sequel, we shall make the standing assumption that  {\bf (H$\vartheta$)} holds.

\vspace{3mm}

We can now state the main result of this section.

\begin{Theorem}
\label{T:EllipticVisc}
Let Assumptions {\bf (H1)},  {\bf (H2)}, {\bf (H$A$)}, and {\bf (H$\vartheta$)} hold. Then, for any $\beta>0$, the function $v^\beta$ in \eqref{vbeta} does not depend on the variable $a$ on $\R^d\times\mathring A$:
\begin{equation}
\label{va}
v^\beta(x,a) \ = \ v^\beta(x,a'), \qquad  a,a' \in \mathring A,
\end{equation}
for all $x$ $\in$ $\R^d$, and we set  by misuse of notation: $v^\beta(x)$ $=$ $v^\beta(x,a)$, $x\in\R^d$,  for any $a\in\mathring A$. Then, $v^\beta$ is the unique continuous viscosity solution to equation \eqref{ellipticbetapde}, i.e., it is continuous on $\R^d$ and it is a viscosity supersolution (resp. subsolution) to \eqref{ellipticbetapde}, namely:
\beqs
\beta \,\varphi(x) & \geq \ (resp.\text{ }\leq) &  \sup_{a\in A} \big[\Lc^a\varphi(x) + f(x,a,\beta\,\varphi(x))\big],
\enqs
for any $x\in\R^d$ and any $\varphi\in C^2(\R^d)$ such that
\beqs
0 \;\, = \;\, (v^\beta-\varphi)(x) \ = \ \min_{\R^d}(v^\beta-\varphi) \quad (\text{resp. }\max_{\R^d}(v^\beta-\varphi)).
\enqs
Moreover, there exists some positive constant $C$ independent of $\beta$ such that:
\beq
|\beta v^\beta(x)|  &\leq &  C (1+|x|), \;\;\; \forall x \in \R^d,   \label{vlin}\\
|v^\beta(x)-v^\beta(x')| & \leq &  C|x-x'|,  \qquad \forall\,x,x'\in\R^d. \label{E:vbeta_Lip}
\enq
\end{Theorem}
\textbf{Proof.} We use the corresponding result for the parabolic case in  \cite{khapha12} to prove the non dependence of $v^\beta$ on $a$, and then the 
viscosity property to the elliptic equation. More precisely, we start by  observing that, for any $\beta>0$ and $T>0$, 
$v^{\beta,n}$ is a viscosity solution to the parabolic PDE on $[0,T]\times\R^d\times A$:
\beqs
\beta \,w(t,x,a)  & =  &  \Dt{w}(t,x,a) + \Lc^a w(t,x,a) + \Mc^a w(t,x,a) + f(x,a,\beta\,w(t,x,a)) \\
& &  \ + n\int_A(w(t,x,a') -w(t,x,a))_+\,\vartheta(da'),
\enqs
with terminal condition $w(T,x,a)=v^{\beta,n}(x,a)$, for all $(x,a)\in\R^d\times A$, i.e. 
\beq
\beta \,\varphi(t,x,a)  &\geq \ (\text{resp. }\leq) & \Dt{\varphi}(t,x,a) + \Lc^a\varphi(t,x,a) + \Mc^a\varphi(t,x,a) + f(x,a,\beta\,\varphi(t,x,a)) \nonumber \\
& &  \ + n\int_A(\varphi(t,x,a') -\varphi(t,x,a))_+\,\vartheta(da'),  \label{E:v^beta,n-varphi_time}
\enq
for any $(t,x,a)\in[0,T)\times\R^d\times A$ and any $\varphi\in C^{1,2}([0,T]\times(\R^d\times\R^q))$ such that
\begin{equation}
\label{E:v^beta,n-varphi_time2}
0 \ = \ (v^{\beta,n}-\varphi)(t,x,a) \ = \ \min_{[0,T]\times\R^d\times A}(v^{\beta,n}-\varphi) \quad (\text{resp. }\max_{[0,T]\times\R^d\times A}(v^{\beta,n}-\varphi)).
\end{equation}
To prove this, for any $(t,x,a)\in[0,T)\times\R^d\times A$ and any $\varphi\in C^{1,2}([0,T]\times(\R^d\times A))$ satisfying 
\eqref{E:v^beta,n-varphi_time}-\eqref{E:v^beta,n-varphi_time2}, set $\psi_t(x',a'):=\varphi(t,x',a')$, for all $(x',a')\in\R^d\times\R^q$. 
Then $\psi_t\in C^2(\R^d\times\R^q)$ and
\[
0 \ = \ (v^{\beta,n}-\psi_t)(x,a) \ = \ \min_{\R^d\times A}(v^{\beta,n}-\psi_t) \quad (\text{resp. }\max_{\R^d\times A}(v^{\beta,n}-\psi_t)).
\]
Observing that  $\partial_t\varphi(t,x,a)\leq0$ (resp. $\partial_t\varphi(t,x,a)\geq0$) since $\varphi(t,x,a)=\max_{t'\in[0,T)}\varphi(t',x,a)$ (resp. $\varphi(t,x,a)=\min_{t'\in[0,T)}\varphi(t',x,a)$), it follows from Proposition \ref{P:ViscPropv^beta,n} that \eqref{E:v^beta,n-varphi_time} is satisfied. 
Then, from Theorem 3.1  in \cite{khapha12}, we deduce by sending $n$ to infinity  that the function $v^\beta$ does not depend on the variable $a$, and so \eqref{va} holds. We should point out that in \cite{khapha12} the terminal condition of the parabolic PDE solved by $v^{\beta,n}$ does not depend on $n$, contrary to our case. However, the part of the proof of Theorem 3.1 in \cite{khapha12} regarding the independence with respect to the variable $a$ does not involve the terminal condition, so the result still holds. Next, we obtain again 
from Theorem 3.1 in \cite{khapha12} that $v^\beta$ solves in the viscosity sense the parabolic PDE on  $[0,T)\times\R^d\times A$ (as before, we do not call in the terminal condition):
\beqs
\beta \,w(t,x,a) - \Dt{w}(t,x,a) - \sup_{a\in A}\big[ \Lc^a w(t,x,a) + f(x,a,\beta\,w(t,x,a)) \big]   &= & 0.
\enqs
Since this equation holds for any $T$, and $v^\beta$ does not depend on $t$,  we obtain that $v^\beta$ is a viscosity solution to the elliptic equation \eqref{ellipticbetapde}. The uniqueness follows  from Theorem 7.4 in \cite{ishii89}. 
Finally, the linear growth and Lipschitz properties \reff{vlin}-\eqref{E:vbeta_Lip} of $v^\beta$ are direct consequences of the corresponding properties 
\reff{growthYbeta} for  $Y_0^{x,a,\beta}$ and \eqref{v^beta,n_Lipschitz} for $v^{\beta,n}$, respectively. 
\ep

\begin{Remark}
{\rm
Notice that, from the identification $Y^{x,a,\beta}=v^\beta(X^{x,a})$, for any $x\in\R^d$ and for some $a\in\mathring A$, and the Lipschitz property \eqref{E:vbeta_Lip}, it follows that $Y^{x,a,\beta}$ is a continuous process, so that $U^{x,a,\beta}\equiv0$, while $K^{x,a,\beta}$ is also a continuous process.
\ep
}
\end{Remark}

\section{Ergodic HJB equation and ergodic BSDE with nonpositive jumps}

\setcounter{equation}{0} 
\setcounter{Theorem}{0} \setcounter{Proposition}{0}
\setcounter{Corollary}{0} \setcounter{Lemma}{0}
\setcounter{Definition}{0} \setcounter{Remark}{0}

This section is devoted to the existence of a solution pair $(\lambda,\phi)$ to the ergodic HJB equation
\beq \label{ergodicHJB}
\lambda - \sup_{a\in A}\big[\Lc^a\phi  + f(x,a,\lambda)\big] & = & 0, \;\;\; \mbox{ on } \; \R^d,
\enq
and to its probabilistic representation in terms of ergodic BSDE with nonpositive jumps.  We first give  the definition of viscosity solution to equation \eqref{ergodicHJB}.

\begin{Definition}
\textup{(i)} A pair $(\lambda,\phi)$, with $\lambda$ a real number and $\phi\colon\R^d\rightarrow\R$ a lower $($resp. upper$)$ semicontinuous function,  is called a \emph{viscosity supersolution} 
$($resp. \emph{viscosity subsolution}$)$ to equation \eqref{ergodicHJB} if
\[
\lambda \ \geq \ (resp.\text{ }\leq) \ \sup_{a\in A}\{\Lc^a\varphi(x) + f(x,a,\lambda)\},
\]
for any $x\in\R^d$ and any $\varphi\in C^2(\R^d)$ such that
\[
(\phi-\varphi)(x) \ = \ \min_{\R^d}(\phi-\varphi) \quad (\text{resp. }\max_{\R^d}(\phi-\varphi)).
\]
\textup{(ii)} A pair $(\lambda, \phi)$, with $\lambda$ a real number and $\phi\colon\R^q\rightarrow\R$ a continuous function, 
 is called a \emph{viscosity solution} to equation \eqref{ergodicHJB} if it is both a viscosity supersolution and a viscosity subsolution to \eqref{ergodicHJB}.
\end{Definition}

\begin{Theorem}
\label{T:ErgodicVisc}
Let  Assumptions {\bf (H1)}, {\bf (H2)}, {\bf (H$A$)}, and {\bf (H$\vartheta$)} hold. Then, there exists a visco\-sity solution pair $(\lambda,\phi)$ to \reff{ergodicHJB} with 
$\phi$ Lipschitz on $\R^d$. 
\end{Theorem}
{\bf Proof.}  
We follow the approximation procedure as in  \cite{ichihara12} or \cite{fuhrman_hu_tess09}: for any $\beta$ $>$ $0$, let $v^\beta$ be the solution of \reff{ellipticbetapde} given in Theorem 
\ref{T:EllipticVisc},  and define  $\lambda_\beta$ $\in$ $\R$, and the function $\phi^\beta\colon\R^d\rightarrow\R$ by:
\beqs
\lambda_\beta \; := \; \beta v^\beta(0), & & \phi^\beta(x) \;  := \;    v^\beta(x)-v^\beta(0), \qquad  \; x\in\R^d. 
\enqs
By \reff{vlin}-\reff{E:vbeta_Lip},  we see that there exists some positive constant $C$ independent of $\beta$ such that 
\beqs
\sup_{\beta > 0} |\lambda_\beta| & < & \infty, \\
\sup_{\beta > 0} |\phi^\beta(x)| & \leq &  C |x|, \qquad \forall\,x\in\R^d.  
\enqs
Then, the family of functions $(\phi^\beta)_\beta$ is equicontinuous and uniformly bounded on every compact subset of $\R^d$. Thus, by means of Ascoli-Arzel\`a theorem (for more details, see (4.4) and (4.5) in \cite{fuhrman_hu_tess09}), we can construct a sequence $(\beta_k)_k$ decreasing monotonically to zero such that, for all $x\in\R^d$,
\beq \label{lambda}
\lambda_{\beta_k}  \;  \overset{k\rightarrow\infty}{\longrightarrow}  \; \lambda, & & 
\phi^{\beta_k}(x) \ \overset{k\rightarrow\infty}{\longrightarrow} \ \phi(x), 
\enq
for some  real constant $\lambda$, and a Lipschitz function $\phi\colon \R^d\rightarrow\R$. Moreover, the convergence of $\phi^{\beta_k}$ towards $\phi$ is uniform on compact sets. 
 
Now, from the elliptic equation \eqref{ellipticbetapde} satisfied by $v^\beta$,  and by definition of $(\lambda_\beta, \phi^\beta)$, we see that $\phi^\beta$ is a viscosity solution to: 
\beq
\label{ellipticbetapde2}
\lambda_\beta +  \beta \phi^\beta  -  \sup_{a\in A} \big[ \Lc^a \phi^{\beta}  +  f(x,a,\lambda_\beta  + \beta \phi^\beta) \big] &=& 0, \;\;\; \mbox{ on } \; \R^d. 
\enq
Let us denote by 
\begin{align*}
F_k(x,r,p,M) \ &:= \ \lambda_{\beta_k} + \beta_k\,r \\
&\quad \ - \sup_{a\in A}\big[ b(x,a).p  + \frac{1}{2}\text{tr}\big(\sigma\sigma\trans(x,a)M\big) + f(x,a,\lambda_{\beta_k} + \beta_k\,r)\big], 
\end{align*}
so that by \reff{ellipticbetapde2},  $\phi^{\beta_k}$ is a viscosity solution to: 
\beqs
F_k(x,\phi^{\beta_k},D_x\phi^{\beta_k},D_x^2\phi^{\beta^k}) &=& 0,
\enqs
and set 
\[
G(x,r,p,M) \ := \ \lambda - \sup_{a\in A}\big[ b(x,a).p  + \frac{1}{2}\text{tr}\big(\sigma\sigma\trans(x,a)M\big) + f(x,a,\lambda)\big],
\]
for all $(x,r,p,M)\in\R^d\times\R\times\R^d\times\R^{d\times d}$. Then, it follows from \eqref{lambda} that $\lim_{k\rightarrow\infty} F_k(x,r,p,M)$ $=$ $G(x,r,p,M)$. As a consequence, from the method of half-relaxed limits of Barles and Perthame, see Remark 6.3 in \cite{crandallishiilions92}, we deduce by \reff{lambda} that  $\phi$ is a viscosity solution to:
\beqs
G(x,\phi,D_x\phi,D_x^2\phi) &=& 0,
\enqs 
i.e., $(\lambda,\phi)$ is a viscosity solution to equation \eqref{ergodicHJB}.
\ep

\vspace{3mm}

We postpone the uniqueness problem of the ergodic equation \reff{ergodicHJB} to the next section, and conclude this section by  providing a probabilistic representation formula for a solution 
to the ergodic equation.   Let us introduce  the ergodic BSDE with nonpositive jumps, $\P$-a.s.,
\begin{align}
\label{ergodicbsde}
Y_t \ &= \ Y_T + \int_t^T \big(f(X_s,I_s,\lambda)-\lambda\big) ds + K_T - K_t \notag \\
&\quad \ - \int_t^T Z_s dW_s - \int_t^T\int_A U_s(a) \tilde\mu(ds,da), \qquad 0 \leq t \leq T < \infty
\end{align}
and
\begin{equation}
\label{JumpConErgodic}
U_t(a) \ \leq \ 0, \qquad dt\otimes d\P\otimes\vartheta(da)\text{-a.e.}
\end{equation}
Here, in addition to the components $(Y,Z,U,K)$, the real number $\lambda$ is  part of the unknowns of the ergodic BSDE.  We recall that ergodic BSDEs driven by Brownian motion have been defined  in 
\cite{fuhrman_hu_tess09}  for the study of optimal control problems on the drift of a diffusion process,  which are related to ergodic semilinear HJB equations.  In this paper, we extend this definition by considering the jump 
constraint \reff{JumpConErgodic}, and our first purpose  is to introduce a notion of ``minimal'' solution to \eqref{ergodicbsde}-\eqref{JumpConErgodic}.  However,  we notice that the ``natural"  definition of minimal solution as the solution $(\bar Y,\bar Z,\bar U,\bar K,\lambda)\in\mathbf{S_{\text{loc}}^2}\times\mathbf{L_{\textup{loc}}^2(W)}\times\mathbf{L_{\textup{loc}}^2(\tilde\mu)}\times\mathbf{K_{\textup{loc}}^2}\times\R$ to \eqref{ergodicbsde}-\eqref{JumpConErgodic} such that for any other solution $(\tilde Y,\tilde Z,\tilde U,\tilde K,\tilde \lambda)\in\mathbf{S_{\text{loc}}^2}\times\mathbf{L_{\textup{loc}}^2(W)}\times\mathbf{L_{\textup{loc}}^2(\tilde\mu)}\times\mathbf{K_{\textup{loc}}^2}\times\R$ to \eqref{ergodicbsde}-\eqref{JumpConErgodic} we have $\bar Y\leq\tilde Y$,  is not relevant  in this case, since $(\bar Y-c,\bar Z,\bar U,\bar K,\lambda)$, with $c>0$, would be another solution to \eqref{ergodicbsde}-\eqref{JumpConErgodic}, contradicting the minimality of $(\bar Y,\bar Z,\bar U,\bar K,\lambda)$. For this reason, we give the following definition of minimal solution to 
the ergodic BSDE with nonpositive jumps \eqref{ergodicbsde}-\eqref{JumpConErgodic}.

\begin{Definition}
\label{D:ErgMinimality}
A quintuple $(\underline Y,\underline Z,\underline U,\underline K,\lambda)\in\mathbf{S_{\textup{loc}}^2}\times\mathbf{L_{\textup{loc}}^2(W)}\times\mathbf{L_{\textup{loc}}^2(\tilde\mu)}\times\mathbf{K_{\textup{loc}}^2}\times\R$ is called a \emph{minimal solution} to \eqref{ergodicbsde}-\eqref{JumpConErgodic} if, for any $T>0$, $(\underline Y_{|[0,T]},\underline Z_{|[0,T]},\underline U_{|[0,T]},\underline K_{|[0,T]})$ is a minimal solution to the BSDE with nonpositive jumps on $[0,T]$ with terminal condition $\underline Y_T$, $\P$-a.s.,
\begin{align}
\label{BSDE_delta}
Y_t \ &= \ \underline Y_T + \int_t^T \big(f(X_s,I_s,\lambda)-\lambda\big) ds + K_T - K_t \notag \\
&\quad \ - \int_t^T Z_s dW_s - \int_t^T \int_A U_s(a)\tilde\mu(ds,da), \qquad 0 \leq t \leq T,
\end{align}
together with the jump constraint
\begin{equation}
\label{JumpCon_delta}
U_t(a) \ \leq \ 0, \qquad dt\otimes d\P\otimes\vartheta(da)\text{-a.e.}
\end{equation}
In other words, for any other solution $(\tilde Y,\tilde Z,\tilde U,\tilde K)\in\mathbf{S^2(0,T)}\times\mathbf{L^2(W;0,T)}\times\mathbf{L^2(\tilde\mu;0,T)}\times\mathbf{K^2(0,T)}$ to \eqref{BSDE_delta}-\eqref{JumpCon_delta} we have $\underline Y_t \leq \tilde Y_t$, $\P$-a.s., for all $0\leq t\leq T$.
\end{Definition}

\begin{Remark}
{\rm We do not know a priori if there is uniqueness of a minimal solution to \eqref{ergodicbsde}-\eqref{JumpConErgodic}.  This will be discussed later in Remark \ref{remunilambda}. 
\ep
}
\end{Remark}

We shall now prove that any Lipschitz continuous viscosity solution to the ergodic HJB equation \eqref{ergodicHJB},   admits a probabilistic representation in terms of a minimal solution to the ergodic BSDE \eqref{ergodicbsde}.

\begin{Theorem}
\label{T:FeynmanKac_phi}
Let  Assumptions {\bf (H1)}, {\bf (H2)}, {\bf (H$A$)}, and {\bf (H$\vartheta$)} hold.  Let $(\lambda,\phi)$, with $\phi$ Lipschitz, be a viscosity solution to the ergodic equation \eqref{ergodicHJB}. 
Then, $(\lambda,\phi)$ can be represented by means of a minimal solution to the ergodic BSDE \eqref{ergodicbsde}-\eqref{JumpConErgodic}, namely:
\begin{itemize}
\item[\textup{(i)}] For any $a\in\mathring{A}$, there exists a minimal solution $(\underline Y^{x,a},\underline Z^{x,a},\underline U^{x,a},\underline K^{x,a},\lambda)\in\mathbf{S_{\textup{loc}}^2}\times\mathbf{L_{\textup{loc}}^2(W)}\times\mathbf{L_{\textup{loc}}^2(\tilde\mu)}\times \mathbf{K_{\textup{loc}}^2}\times\R$ to the ergodic BSDE with nonpositive jumps, $\P$-a.s.,
\begin{align}
\underline Y_t^{x,a} \ &= \ \underline Y_T^{x,a} + \int_t^T \big(f(X_s^{x,a},I_s^a,\lambda)-\lambda\big)ds + \underline K_T^{x,a} - \underline K_t^{x,a} - \int_t^T \underline Z_s^{x,a} dW_s \nonumber \\
&\quad \ - \int_t^T\int_A \underline U_s^{x,a}(a') \tilde\mu(ds,da'), \qquad\qquad 0\leq t\leq T<\infty, \label{BSDEunderlineY}
\end{align}
and
\beq
\underline U_t^{x,a}(a') \ \leq \ 0, \qquad dt\otimes d\P\otimes \vartheta(da')\text{-a.e.} \label{underlineconstraint}
\enq
\item[\textup{(ii)}] $\underline Y_t^{x,a}=\phi(X_t^{x,a})$ for $t\geq0$, and, in particular,
\[
\phi(x) \ = \ \underline Y_0^{x,a}, \qquad \text{for all }x\in\R^d,
\]
for some $a\in\mathring{A}$.
\end{itemize} 
\end{Theorem}
\textbf{Proof.}
We start by observing that for any $T>0$, $\phi$ is a viscosity solution to the following parabolic equation, in the unknown $\psi$,
\begin{equation}
\label{elliptic_pde_phi}
\begin{cases}
-\dfrac{\partial\psi(t,x)}{\partial t} - \ \sup_{a\in A}\big[ \Lc^a\psi(t,x) + f(x,a,\lambda) - \lambda\big] = 0,  \qquad\qquad &(t,x)\in[0,T)\times\R^d, \\
\psi(T,x) = \phi(x), &x\in\R^d.
\end{cases}
\end{equation}
It follows from Theorem 7.4 in \cite{ishii89} that $\phi$ is the unique uniformly continuous viscosity solution to \eqref{elliptic_pde_phi}. From Theorem 3.1 in \cite{khapha12}, we see that this viscosity solution admits a representation in terms of the minimal solution to a BSDE with nonpositive jumps on $[0,T]$. More precisely, for any $T>0$, let $(Y^{x,a,T},Z^{x,a,T},U^{x,a,T},K^{x,a,T})\in\mathbf{S^2(0,T)}\times\mathbf{L^2(W;0,T)}\times\mathbf{L^2(\tilde\mu;0,T)}\times\mathbf{K^2(0,T)}$ be the minimal solution to the BSDE with nonpositive jumps on $[0,T]$, $\P$-a.s.,
\begin{align}
\label{ergodicbsde_proof}
Y_t \ &= \ \phi(X_T^{x,a}) + \int_t^T \big(f(X_s^{x,a},I_s^a,\lambda) - \lambda\big) ds + K_T-K_t \notag \\
&\quad \ - \int_t^T Z_s dW_s - \int_t^T \int_A U_s(a')\tilde\mu(ds,da'), \qquad 0 \leq t \leq T,
\end{align}
and
\begin{equation}
\label{ergodicbsde_proof_jump}
U_t(a') \ \leq \ 0, \qquad dt\otimes d\P\otimes\vartheta(da')\text{-a.e.}
\end{equation}
Then, from Theorem 3.1 in \cite{khapha12} we see that $Y_t^{x,a,T}=\phi(X_t^{x,a})$, for all $0\leq t\leq T$. The identification $Y_t^{x,a,T}=\phi(X_t^{x,a})$ also implies that $Y^{x,a,T}$ does not depend on $T$. Moreover, using the fact that the $Y$ component remains the same, reasoning as in Remark \ref{P:Uniqueness}, we can prove that $Z^{x,a,T}$, $U^{x,a,T}$, and $K^{x,a,T}$ do not depend on $T$. Let  $\underline Y^{x,a}_{|[0,T]}=Y^{x,a,T}$, $\underline Z^{x,a}_{|[0,T]}=Z^{x,a,T}$, $\underline U^{x,a}_{|[0,T]}=U^{x,a,T}$, and $\underline K^{x,a}_{|[0,T]}=K^{x,a,T}$, for any $T>0$. Then, we see that, for any $T>0$, the quadruple $(\underline Y^{x,a}_{|[0,T]},\underline Z^{x,a}_{|[0,T]},\underline U^{x,a}_{|[0,T]},\underline K^{x,a}_{|[0,T]})$ is the minimal solution to \eqref{ergodicbsde_proof}-\eqref{ergodicbsde_proof_jump}, from which we deduce that $(\underline Y^{x,a},\underline Z^{x,a},\underline U^{x,a},\underline K^{x,a},\lambda)$ is a minimal solution to the ergodic BSDE \eqref{ergodicbsde}-\eqref{JumpConErgodic}.
\ep

\begin{Remark}
\label{R:MinSol-->ErgHJB}
{\rm
Notice that, a minimal solution to the ergodic BSDE with nonpositive jumps \eqref{ergodicbsde}-\eqref{JumpConErgodic} provides a viscosity solution to the ergodic HJB equation \eqref{ergodicHJB}. In other words, the converse of the result stated in Theorem \ref{T:FeynmanKac_phi} is also true. More precisely, let $\phi\colon\R^d\times\R^q\rightarrow\R$ be a Lipschitz function and $\lambda\in\R$. For any $x\in\R^d$, consider $a\in\mathring{A}$ and a minimal solution $(\underline Y^{x,a},\underline Z^{x,a},\underline U^{x,a},\underline K^{x,a},\lambda)\in\mathbf{S_{\textup{loc}}^2}\times\mathbf{L_{\textup{loc}}^2(W)}\times\mathbf{L_{\textup{loc}}^2(\tilde\mu)}\times \mathbf{K_{\textup{loc}}^2}\times\R$ to the ergodic BSDE with nonpositive jumps \eqref{BSDEunderlineY}-\eqref{underlineconstraint}, such that $\underline Y_t^{x,a}=\phi(X_t^{x,a},I_t^a)$, $t\geq0$. Then, $\phi$ does not depend on the variable $a$ and $(\lambda,\phi)$ is a viscosity solution to the ergodic equation \eqref{ergodicHJB}. Indeed, fix $T>0$, then by Definition \ref{D:ErgMinimality} we know that $(\underline Y^{x,a},\underline Z^{x,a},\underline U^{x,a},\underline K^{x,a})$ is the unique minimal solution to the BSDE with nonpositive jumps on $[0,T]$ in \eqref{BSDE_delta}-\eqref{JumpCon_delta}. It follows from Theorem 3.1 in \cite{khapha12} that $\phi$ does not depend on $a$ (this last property follows from Proposition 3.2 in \cite{khapha12}, which does not call in the terminal condition; therefore, even if in our case, as opposite to \cite{khapha12}, the terminal condition depends on $a$, the result is still valid) and it is a viscosity solution to equation \eqref{elliptic_pde_phi}. As a consequence, $(\lambda,\phi)$ is a viscosity solution to the ergodic equation \eqref{ergodicHJB}.
\ep
}
\end{Remark}

\section{Convergence of solutions}

\setcounter{equation}{0} 
\setcounter{Theorem}{0} \setcounter{Proposition}{0}
\setcounter{Corollary}{0} \setcounter{Lemma}{0}
\setcounter{Definition}{0} \setcounter{Remark}{0}

Let us consider the parabolic fully nonlinear equation of HJB type: 
 \begin{equation}
 \label{parabolicpde}
\begin{cases}
\DT{v}  -  \sup_{a\in A}\big[ \Lc^a v  + f\big(x,a,\frac{v}{T+1}\big)\big], \; = \; 0,  \;\;\; \mbox{ on } (0,\infty)\times\R^d, \\
v(0,.) \; = \; h, \;\;\; \mbox{ on } \; \R^d,
\end{cases}
\end{equation}
where $h$ is a Lipschitz function on $\R^d$. Existence and uniqueness of a solution $v(T,x)$ to \reff{parabolicpde} is studied in the next paragraph.  Now, let us consider a solution pair $(\lambda,\phi)$, with $\phi$ Lipschitz,  
to the ergodic equation \eqref{ergodicHJB}.  The main result of this section is to prove that
\begin{equation}
\label{convergone}
\frac{v(T,x)}{T} \ \overset{T\rightarrow\infty}{\longrightarrow} \ \lambda.
\end{equation}
Consequently, this will show the uniqueness of the component $\lambda$ in \reff{ergodicHJB}.  We shall end this part of the paper by proving a stronger  result than \eqref{convergone} under additional assumptions. More precisely, suppose that $\phi$ belongs to $C^2(\R^d)$ so that  $(\lambda,\phi)$ is a classical solution to the ergodic equation \eqref{ergodicHJB}, and 
assume that in the ergodic equation,  the supremum is attained at $a=\underline\alpha(x)$, for every $x\in\R^d$, for some locally Lipschitz function $\underline\alpha\colon A\rightarrow\R^d$. 
Then, the following convergence holds
\begin{equation}
\label{convergtwo}
v(T,x)-(\lambda{T}+\phi(x)) \ \overset{T\rightarrow\infty}{\longrightarrow} \ c,
\end{equation}
for some constant $c$. In particular, $\phi$ is uniquely determined up to a constant.

\subsection{Wellposedness of the parabolic equation \eqref{parabolicpde}}

We shall build a solution to equation \eqref{parabolicpde} through BSDE methods, as this construction will be useful in the sequel. More precisely, from Theorem 3.1 in \cite{khapha12}, under {\bf (H1)}(i), {\bf (H2)}(i),  {\bf (H$A$)}, and {\bf (H$\vartheta$)}, there exists a uniformly continuous viscosity solution $v$ to equation \eqref{parabolicpde}, which admits the following probabilistic representation formula
\begin{equation}
\label{v}
v(T,x) \ = \ Y_0^{x,a,T}, \qquad (T,x)\in[0,\infty)\times\R^d,
\end{equation}
for any $a\in\mathring{A}$, where $(Y^{x,a,T},Z^{x,a,T},U^{x,a,T},K^{x,a,T})\in\mathbf{S^2(0,T)}\times\mathbf{L^2(W;0,T)}\times\mathbf{L^2(\tilde\mu;0,T)}\times\mathbf{K^2(0,T)}$, with $v(T-t,X_t^{x,a})=Y_t^{x,a,T}$ for all $0\leq t\leq T$, is the unique minimal solution to the BSDE with nonpositive jumps on $[0,T]$, $\P$-a.s.,
\beq
\label{BSDE_v}
Y_t^{x,a,T} \ &= & h(X_T^{x,a}) + \int_t^T f\big(X_s^{x,a},I_s^a,\frac{Y_s^{x,a,T}}{T-s+1}\big) ds + K_T^{x,a,T} - K_t^{x,a,T} \\
& & \ - \int_t^T Z_s^{x,a,T} dW_s - \int_t^T\int_A U_s^{x,a,T}(a') \tilde\mu(ds,da'), \qquad 0\leq t\leq T \notag
\enq
and
\begin{equation}
\label{JumpConstraint_v}
U_t^{x,a,T}(a') \ \leq \ 0, \qquad dt\otimes d\P\otimes \vartheta(da')\text{-a.e.}
\end{equation}
Moreover, from Theorem 7.4 in \cite{ishii89}, we have that $v$ is the unique uniformly continuous viscosity solution to \eqref{parabolicpde} (observe that Theorem 7.4 in \cite{ishii89} applies to elliptic equations on unbounded domains; interpreting $t$ as a space variable, \eqref{parabolicpde} can be seen as an elliptic equation on the space domain $[0,\infty)\times\R^d$, so that we can now apply Theorem 7.4 in \cite{ishii89}).

\begin{Remark}
\label{R:v^T}
{\rm
Notice that Theorem 3.1 in \cite{khapha12} is designed for \emph{backward} parabolic PDEs, while \eqref{parabolicpde} is a forward parabolic equation. However, we can exploit Theorem 3.1 in \cite{khapha12} by proceeding as follows. For any $T>0$, we consider the HJB equation on $[0,T]\times\R^d$:
\begin{equation}
\label{parabolicpdehjb}
\begin{cases}
-\Dt{v^T} - \sup_{a\in A}\big[ \Lc^a v^T(t,x)   +  f\big(x,a,\frac{v^T}{T-t+1}\big)\big] = 0,   \quad (t,x)\in[0,T) \times\R^d, \\
v^T(T,x) \; = \;  h(x),  \quad x\in\R^d.
\end{cases}
\end{equation}
Under {\bf (H1)}(i), {\bf (H2)}(i), {\bf (H$A$)}, and {\bf (H$\vartheta$)}, it follows from Theorem 3.1 in \cite{khapha12} that there exists a uniformly continuous viscosity solution $v^T$ to equation \eqref{parabolicpdehjb}, which admits a probabilistic representation formula in terms of the unique minimal solution to a certain BSDE with nonpositive jumps. Define the function $v(t,x):=v^T(T-t,x)$, for all $(t,x)\in[0,T]\times\R^d$ and for any $T>0$. Notice that $v$ is well-defined, thanks to  uniqueness results (see, e.g., Theorem 7.4 in \cite{ishii89}) of viscosity solutions to equation \eqref{parabolicpdehjb} (in particular, $v^T(T-t,x)=v^{T'}(T'-t,x)$, for any $(t,x)\in[0,T]\times\R^d$ and $0\leq T\leq T'<\infty$). Moreover, from the viscosity properties of $v^T$ it follows that $v$ is the unique uniformly continuous viscosity solution to equation \eqref{parabolicpde}. Then, from the probabilistic representation formula for $v^T$ we deduce the representation formula \eqref{v} for $v$.
\ep
}
\end{Remark}

In conclusion, we have proved the following result.

\begin{Proposition}
Let  Assumptions {\bf (H1)}, {\bf (H2)}, {\bf (H$A$)}, and {\bf (H$\vartheta$)} hold.  
Then, the function $v$ given by \eqref{v} is the unique uniformly continuous viscosity solution to \eqref{parabolicpde} on $[0,\infty)\times\R^d$. 
\end{Proposition}

\subsection{First convergence result: the proof of (\ref{convergone})}

Let $(\lambda,\phi)$, with $\phi$ Lipschitz, be a viscosity solution to the ergodic equation \eqref{ergodicHJB}. Let us introduce the function $w\colon[0,\infty)\times\R^d\rightarrow\R$ given by
\begin{equation}
\label{w}
w(T,x) \ := \ v(T,x) - \big(\lambda T+\phi(x)\big),
\qquad (T,x)\in[0,\infty)\times\R^d.
\end{equation}
The aim is to state an upper and lower  estimate  for $w$, uniformly in time $T$, so that by dividing by $T$, we obtain the convergence of the long run average $v(T,.)/T$ to $\lambda$. 
Classical PDE arguments (in the case where $f$ does not depend on $y$) rely on the smoothness of  $v$ and $\phi$ in order to prove that $w$ is a sub and supersolution to some PDE without cost or gain function. 
Then  by comparison principle, and under ergodicity conditions, one would obtain for $w$ a lower and upper bound  function which does not depend on time.  In our general framework, the major difficulty is due to the non-regularity in general of $v$ and $\phi$, especially when there is singularity of the diffusion coefficient. In this case, it is not clear, even with the notion of viscosity solution, how to derive an equation 
 for $w$ involving the difference of  $v$ and $\phi$. 
We circumvent this issue by adopting an alternative approach where we use probabilistic representations formulae for $v$ and $\phi$. We are also interested in the case where $f(x,a,y)$ depends on $y$, that we shall actually tackle by using the nondecreasing feature of $f$ in $y$ and imposing the following additional assumption. 

\vspace{2mm}

{\bf (H3)} \hspace{7mm} The function $f$ can be written as $f(x,a,y)$ $=$ $f_0(x,a)+f_1(x,a,y)$, where $f_1$ can be either the zero function or it satisfies, for all $x\in\R^d$, $a\in A$, $y,y'\in\R$,
\beqs
y \ > \ y' \quad & \Longrightarrow & \quad f_1(x,a,y) - f_1(x,a,y') \ \leq \ -\kappa (y-y'),
\enqs
for some constant $\kappa>0$.

\vspace{1mm}

\begin{Theorem}
\label{T:Convergence}
Let Assumptions {\bf (H1)}, {\bf (H2)}, {\bf (H3)}, {\bf (H$A$)}, and {\bf (H$\vartheta$)} hold. Then, there exists a positive constant $C$ such that, the function $w$ defined in \eqref{w} satisfies
\beq
\label{estimate_w}
-C(1+|x|) \;  \leq \;   w(T,x) & \leq & C(1+|x|), \;\;\; (T,x) \in [0,\infty)\times\R^d.
\enq
In particular, we have
\begin{equation}
\label{convergonebis}
\frac{v(T,x)}{T} \ \overset{T\rightarrow\infty}{\longrightarrow} \ \lambda.
\end{equation}
\end{Theorem}

{\begin{Remark}
\label{R:Convergence}
{\rm
We report here the \emph{proof of Theorem \ref{T:Convergence} when $f=f(x,a)$ does not depend on $y$,} since it is much easier. Recall from Remark \ref{R:v^T} that $v(t,x)=v^T(T-t,x)$, where $v^T$ is the unique uniformly continuous viscosity solution to the Hamilton-Jacobi-Bellman equation \eqref{parabolicpdehjb}. Therefore, $v^T$ admits a stochastic control representation, which in terms of $v$ reads
\beqs
v(T,x) &= &  \sup_{\alpha \in \Ac}\,\E\bigg[\int_0^T f\big(X_s^{x,\alpha},\alpha_s\big) ds 
 + h(X_T^{x,\alpha})\bigg], \qquad \forall\,(T,x)\in[0,\infty)\times\R^d,
\enqs
where $\Ac$ is the set of adapted control processes valued in $A$,  and $(X_t^{x,\alpha})_{t\geq0}$ is the unique solution to the controlled equation \eqref{controlX} starting from $x$ at time $0$. Similarly, we know from the proof of Theorem \ref{T:FeynmanKac_phi} that $\phi$ is the unique uniformly continuous viscosity solution to equation \eqref{elliptic_pde_phi}, so that $\phi$ is given by
\beqs
\phi(x) &= &  \sup_{\alpha \in \Ac}\,\E\bigg[\int_0^T \big[f\big(X_s^{x,\alpha},\alpha_s\big) - \lambda\big] ds 
 + \phi(X_T^{x,\alpha})\bigg], \qquad \forall\,(T,x)\in[0,\infty)\times\R^d.
\enqs
From the definition of $w$, we have
\begin{align*}
w(T,x) \ &= \ \sup_{\alpha \in \Ac}\,\E\bigg[\int_0^T f\big(X_s^{x,\alpha},\alpha_s\big) ds 
 + h(X_T^{x,\alpha})\bigg] - \sup_{\alpha \in \Ac}\,\E\bigg[\int_0^T f\big(X_s^{x,\alpha},\alpha_s\big) ds 
 + \phi(X_T^{x,\alpha})\bigg] \\
&\leq \ \sup_{\alpha \in \Ac}\,\E\big[(h-\phi)(X_T^{x,\alpha})\big].
\end{align*}
Proceeding in a similar way, we obtain
\begin{align*}
w(T,x) \ &\geq \ \inf_{\alpha \in \Ac}\E\big[(h-\phi)(X_T^{x,\alpha})\big].
\end{align*}
Since $h$ and $\phi$ are Lipschitz, from estimate \eqref{estimate_x^2} we deduce \eqref{estimate_w}.
\ep
}
\end{Remark}

\noindent\textbf{Proof of Theorem \ref{T:Convergence}.}
We recall from \eqref{v} and Theorem \ref{T:FeynmanKac_phi} the nonlinear Feynman-Kac formulae
\[
\phi(X_t^{x,a}) \ = \ \underline Y_t^{x,a}, \qquad v(T-t,X_t^{x,a}) \ = \ Y_t^{x,a,T}, \qquad  \; 0 \leq t \leq T, 
\]
for all $(T,x)$ $\in$ $[0,\infty)\times\R^d$, and any $a\in\mathring{A}$. Fix then $a$ $\in$ $\mathring{A}$, and define, for $(T,x)$ $\in$ $[0,\infty)\times\R^d$, the process:
\beqs
\tilde Y_t^{x,a,T} & := & Y_t^{x,a,T} - \lambda(T-t) - \underline Y_t^{x,a}, \;\;\;\;\;\;\;  0 \leq t \leq T.
\enqs   
Then, by definition of $w$ in \reff{w}, we have
\beq \label{wtildeY}
w(T,x) & = & \tilde Y_0^{x,a,T}, \qquad \forall\,(T,x)\in[0,\infty)\times\R^d.
\enq
Moreover, from the BSDE \eqref{BSDE_v}-\eqref{JumpConstraint_v} for $Y^{x,a,T}$ and \reff{BSDEunderlineY}-\reff{underlineconstraint} for $\underline Y^{x,a}$, 
we derive the BSDE for $\tilde Y^{x,a,T}$: 
\begin{align*}
\tilde Y_t^{x,a,T} \ &= \ (h - \phi)(X_T^{x,a}) + \int_t^T\big(f\big(X_s^{x,a},I_s^a,\frac{Y_s^{x,a,T}}{T-s+1}\big) - f(X_s^{x,a},I_s^a,\lambda)\big)ds \\
&\quad \ + K_T^{x,a,T} - K_t^{x,a,T} - \big(\underline K_T^{x,a}-\underline K_t^{x,a}\big) - \int_t^T\big(Z_s^{x,a,T}-\underline Z_s^{x,a}\big)dW_s \\
&\quad \ - \int_t^T\int_A \big(U_s^{x,a,T}(a')-\underline U_s^{x,a}(a')\big)\tilde\mu(ds,da'), \qquad 0\leq t\leq T,
\end{align*}
and
\begin{align*}
U_t^{x,a,T}(a') \ &\leq \ 0, \qquad dt\otimes d\P\otimes\vartheta(da')\text{-a.e.} \\
\underline U_t^{x,a}(a') \ &\leq \ 0, \qquad dt\otimes d\P\otimes\vartheta(da')\text{-a.e.}
\end{align*}
We shall now prove  suitable upper and lower bounds for $\tilde Y^{x,a,T}$, and thus for $w(T,x)$.

$\bullet$ \emph{Step 1. Upper bound: $w(T,x)\leq C(1+|x|)$.} Let us consider the BSDE with nonpositive jumps on $[0,T]$:
\begin{align}
\label{BSDE1}
\hat Y_t^{x,a,T} \ &= \ (h - \phi)(X_T^{x,a}) 
+ \int_t^T\big(f\big(X_s^{x,a},I_s^a,\frac{\hat Y_s^{x,a,T}+\lambda(T-s)+\phi(X_s^{x,a})}{T-s+1}\big) \notag \\ 
&\quad \ - f(X_s^{x,a},I_s^a,\lambda)\big)ds + \hat K_T^{x,a,T} - \hat K_t^{x,a,T} \\
 & \quad - \int_t^T\hat Z_s^{x,a,T}dW_s  - \int_t^T\int_A \hat U_s^{x,a,T}(a')\tilde\mu(ds,da') \nonumber
\end{align}
and
\begin{equation}
\label{JumpComp1}
\hat U_t^{x,a,T}(a') \ \leq \ 0, \qquad dt\otimes d\P\otimes\vartheta(da')\text{-a.e.}
\end{equation}
We know  from Theorem 2.1 in \cite{khapha12} that there exists a unique minimal solution
$(\hat Y^{x,a,T}$, $\hat Z^{x,a,T}$, $\hat U^{x,a,T}$, $\hat K^{x,a,T})$  in $\mathbf{S^2(0,T)}\times\mathbf{L^2(W;0,T)}\times\mathbf{L^2(\tilde\mu;0,T)}\times\mathbf{K^2(0,T)}$ to equation \eqref{BSDE1}-\eqref{JumpComp1}. Set $\bar Y_t^{x,a,T}=\hat Y_t^{x,a,T}+\lambda(T-t)+\underline Y_t^{x,a}$, $t\in[0,T]$, and recall that $\underline Y_t^{x,a}$ $=$ $\phi(X_t^{x,a})$. Then, from the BSDEs satisfied by $\hat Y^{x,a,T}$ and $\underline Y^{x,a}$, we easily see that  
$(\bar Y^{x,a,T},\hat Z^{x,a,T}+\underline Z^{x,a},\hat U^{x,a,T}+\underline U^{x,a},\hat K^{x,a,T}+\underline K^{x,a})$ is a solution to \eqref{BSDE_v}-\eqref{JumpConstraint_v}. From the minimality of $(Y^{x,a,T},Z^{x,a,T},U^{x,a,T},K^{x,a,T})$, we get:  $Y_t^{x,a,T}\leq\bar Y_t^{x,a,T}$, and by subtracting to both sides $\lambda(T-t)+\underline Y_t^{x,a}$ we end up with 
\beq \label{YtildeYhat}
\tilde Y_t^{x,a,T} &\leq & \hat Y_t^{x,a,T}, \;\;\; t\in[0,T].
\enq
Let us now derive an upper bound for $\hat Y_0^{x,a,T}$. To this end, we introduce the associated penalized BSDE with jumps on $[0,T]$, for $n$ $\in$ $\N$:
\begin{align}
\label{BSDE1_penalized}
\hat Y_t^{x,a,T,n} \ &= \ (h - \phi)(X_T^{x,a}) 
+ \int_t^T\big(f\big(X_s^{x,a},I_s^a,\frac{\hat Y_s^{x,a,T,n}+\lambda(T-s)+\phi(X_s^{x,a})}{T-s+1}\big) \notag \\ 
&\quad \ - f(X_s^{x,a},I_s^a,\lambda)\big)ds + n\int_t^T\int_A (\hat U_s^{x,a,T,n}(a'))_+ \vartheta(da')ds \\
 & \quad - \int_t^T\hat Z_s^{x,a,T,n}dW_s - \int_t^T\int_A \hat U_s^{x,a,T,n}(a')\tilde\mu(ds,da').  \nonumber
\end{align}
From the uniform Lipschitz condition on $f(x,a,y)$ with respect to $y$, together with Assumptions {\bf (H2)}(ii) and {\bf (H3)}, we have
\beq
\label{Ineq}
& & f\big(X_s^{x,a},I_s^a,\frac{\hat Y_s^{x,a,T,n}+\lambda(T-s) +\phi(X_s^{x,a})}{T-s+1}\big)   \ - f(X_s^{x,a},I_s^a,\lambda)   \label{flipalpha} \notag \\
&\leq&  \rho_s^n \big( \hat Y_s^{x,a,T,n} + \phi(X_s^{x,a}) - \lambda \big),   \qquad  0 \leq s \leq T,
\enq
with (we suppose here that $f_1$ in {\bf (H3)} is not the zero function; otherwise, $\rho^n$ can be taken equal to zero everywhere and the proof becomes easier)
\begin{equation}
\label{rho^n}
\rho_s^n \ = \ - \frac{\kappa}{T+1} 1_{\{\hat Y_s^{x,a,T,n} + \phi(X_s^{x,a}) - \lambda>0\}} - L_2 1_{\{\hat Y_s^{x,a,T,n} + \phi(X_s^{x,a}) - \lambda\leq 0\}}, \qquad  0 \leq s \leq T.
\end{equation}
Then, applying  It\^o's formula to $e^{\int_0^t\rho_r^n dr}\hat Y_t^{x,a,T,n}$ between $0$ and $T$, we get from \reff{BSDE1_penalized}:
\begin{align*}
\hat Y_0^{x,a,T,n} \ &= \  e^{\int_0^T \rho_r^n dr}(h-\phi)(X_T^{x,a}) - \int_0^T \rho_s^n e^{\int_0^s\rho_r^n dr} \hat Y_s^{x,a,T,n} ds \\
&\quad \ + \int_0^T e^{\int_0^s\rho_r^n dr}\big( f\big(X_s^{x,a},I_s^a,\frac{\hat Y_s^{x,a,T,n}+\lambda(T-s) +\phi(X_s^{x,a})}{T-s+1}\big) \\
&\quad \ - f(X_s^{x,a},I_s^a,\lambda)\big)ds + n\int_0^T \int_A e^{\int_0^s\rho_r^n dr}(\hat U_s^{x,a,T,n}(a'))_+\vartheta(da')ds \\
&\quad \ - \int_0^T e^{\int_0^s\rho_r^n dr}\hat Z_s^{x,a,T,n}dW_s - \int_0^T \int_A e^{\int_0^s\rho_r^n dr}\hat U_s^{x,a,T,n}(a')\tilde\mu(ds,da').
\end{align*}
Using \eqref{Ineq}, we obtain
\begin{align}
\label{Ineq2}
\hat Y_0^{x,a,T,n} \ &\leq \  e^{\int_0^T \rho_r^n dr}(h-\phi)(X_T^{x,a}) + \int_0^T \rho_s^n e^{\int_0^s\rho_r^n dr}\big(\phi(X_s^{x,a})-\lambda\big) ds \notag \\
&\quad \ + n\int_0^T \int_A e^{\int_0^s\rho_r^n dr}(\hat U_s^{x,a,T,n}(a'))_+\vartheta(da')ds \\
&\quad \ - \int_0^T e^{\int_0^s\rho_r^n dr}\hat Z_s^{x,a,T,n}dW_s - \int_0^T \int_A e^{\int_0^s\rho_r^n dr}\hat U_s^{x,a,T,n}(a')\tilde\mu(ds,da'). \notag
\end{align}
Now, from Proposition 2.1 in \cite{khapha12} we have the following dual representation formula for the right-hand side of \eqref{Ineq2}:
\begin{align*}
&e^{\int_0^T \rho_r^n dr}(h-\phi)(X_T^{x,a}) + \int_0^T \rho_s^n e^{\int_0^s\rho_r^n dr}\big(\phi(X_s^{x,a})-\lambda\big) ds \notag \\
& + n\int_0^T \int_A e^{\int_0^s\rho_r^n dr}(\hat U_s^{x,a,T,n}(a'))_+\vartheta(da')ds \\
& - \int_0^T e^{\int_0^s\rho_r^n dr}\hat Z_s^{x,a,T,n}dW_s - \int_0^T \int_A e^{\int_0^s\rho_r^n dr}\hat U_s^{x,a,T,n}(a')\tilde\mu(ds,da') \\
&= \ \sup_{\nu\in\Vc} \E^\nu\bigg[ e^{\int_0^T \rho_r^n dr} (h-\phi)(X_T^{x,a}) + \int_0^T \rho_s^n e^{\int_0^s \rho_r^n dr}\big(\phi(X_s^{x,a})-\lambda\big) ds\bigg].
\end{align*}
Therefore, we get
\beqs
\hat Y_0^{x,a,T,n}  &\leq &  \sup_{\nu\in\Vc} \E^\nu\bigg[ e^{\int_0^T \rho_r^n dr} (h-\phi)(X_T^{x,a}) + \int_0^T \rho_s^n e^{\int_0^s \rho_r^n dr}\big(\phi(X_s^{x,a})-\lambda\big) ds\bigg].
\enqs
From the definition of $\rho^n$ in \eqref{rho^n}, we find
\beqs
\hat Y_0^{x,a,T,n}  &\leq & \sup_{\nu\in\Vc} \E^\nu \bigg[|h-\phi|(X_T^{x,a}) 
+ \max\Big(\frac{\kappa}{T+1},L_2\Big) \int_0^T e^{-\min(\frac{\kappa}{T+1},L_2)s}|\phi(X_s^{x,a})-\lambda| ds\bigg].
\enqs
Recalling that $h$ and $\phi$ are Lipschitz, from estimate \eqref{estimate_x^2} we obtain
\[
\hat Y_0^{x,a,T,n} \ \leq \ C (1 + |x|),
\]
for some positive constant $C$, independent of $x$, $a$, $T$, and $n$. Since from Theorem 2.1 in \cite{khapha12} we have that  $\hat Y_0^{x,a,T,n}$ converges to  $\hat Y_0^{x,a,T}$, as $n$ goes to infinity, 
we get the same estimate: $\hat Y_0^{x,a,T}$ $\leq$ $C(1+|x|)$, and  therefore, from \reff{wtildeY} and \reff{YtildeYhat}, we deduce that  
\beqs
w(T,x) &\leq& C(1+|x|). 
\enqs

$\bullet$ \emph{Step 2. Lower bound: $w(T,x)\geq -C(1+|x|)$.} As in step 1, where we built an upper bound for $\tilde Y^{x,a,T}$ using the minimality of $Y^{x,a,T}$, here we shall construct a lower bound for $\tilde Y^{x,a,T}$ exploiting the minimality of $\underline Y^{x,a}$ in the sense of Definition \ref{D:ErgMinimality}. In particular, we fix $T>0$ and we recall that $(\underline Y^{x,a}_{|[0,T]},\underline Z^{x,a}_{|[0,T]},\underline U^{x,a}_{|[0,T]},\underline K^{x,a}_{|[0,T]})$ is the minimal solution to \eqref{BSDE_delta}-\eqref{JumpCon_delta} on $[0,T]$ with terminal condition $\phi(X_T^{x,a})$. Now, let us  consider the BSDE with nonnegative jumps on $[0,T]$:
\begin{align}
\label{BSDE2}
\check Y_t^{x,a,T} \ &= \ (h - \phi)(X_T^{x,a}) + \int_t^T\big(f\big(X_s^{x,a},I_s^a,\frac{v(T-s,X_s^{x,a})}{T-s+1}\big) - f(X_s^{x,a},I_s^a,\lambda)\big)ds \\
&\quad \  - \big(\check K_T^{x,a,T}-\check K_t^{x,a,T}\big) - \int_t^T\check Z_s^{x,a,T}dW_s   - \int_t^T\int_A \check U_s^{x,a,T}(a')\tilde\mu(ds,da') \notag
\end{align}
and
\begin{equation}
\label{JumpComp2}
\check U_t^{x,a,T}(a') \ \geq \ 0, \qquad dt\otimes d\P\otimes\vartheta(da')\text{-a.e.}
\end{equation}
Theorem 2.1 in \cite{khapha12} gives the existence of  a unique maximal solution $(\check Y^{x,a,T}$, $\check Z^{x,a,T}$, $\check U^{x,a,T},\check K^{x,a,T})\in\mathbf{S^2(0,T)}\times\mathbf{L^2(W;0,T)}\times\mathbf{L^2(\tilde\mu;0,T)}\times\mathbf{K^2(0,T)}$ to \eqref{BSDE2}-\eqref{JumpComp2}. Actually, Theorem 2.1 in \cite{khapha12} is designed for minimal solutions, while here we deal with  maximal solutions; however, it is easy to show that $-\check Y^{x,a,T}$ is a minimal solution to a certain BSDE with nonpositive jumps, therefore we can apply Theorem 2.1 to $-\check Y^{x,a,T}$. Set $\overline Y_t^{x,a,T}=-\check Y_t^{x,a,T}+Y_t^{x,a,T}-\lambda(T-t)$, $t\in[0,T]$, then
\begin{align*}
&(\overline Y^{x,a,T},\overline Z^{x,a,T},\overline U^{x,a,T},\overline K^{x,a,T}) \\
&:= \ (\overline Y^{x,a,T},-\check Z^{x,a,T}+Z^{x,a,T},-\check U^{x,a,T}+U^{x,a,T},\check K^{x,a,T}+K^{x,a,T})
\end{align*}
is a solution to \eqref{BSDE_delta}-\eqref{JumpCon_delta} on $[0,T]$ with 
terminal condition $\phi(X_T^{x,a})$. From the already mentioned minimality of $(\underline Y^{x,a}_{|[0,T]},\underline Z^{x,a}_{|[0,T]},\underline U^{x,a}_{|[0,T]},\underline K^{x,a}_{|[0,T]})$ to 
\eqref{BSDE_delta}-\eqref{JumpCon_delta}, we see that $\underline Y_t^{x,a}\leq\overline Y_t^{x,a,T}$,  
and by subtracting  to both sides $Y_t^{x,a,T}-\lambda(T-t)$, we end up with
\beq \label{YtildeYcheck}
\check Y_t^{x,a,T} &\leq & \tilde Y_t^{x,a,T}, \;\;\; t\in[0,T].
\enq
We now derive a lower bound for $\check Y_0^{x,a,T}$ by means of a dual representation formula. In particular, we see from Theorem 2.2  in \cite{khapha12} that $\check Y_0^{x,a,T}$ admits the dual representation formula (observe that, as for Theorem 2.1 in \cite{khapha12}, Theorem 2.2 in \cite{khapha12} is designed for minimal solutions, while here we deal with  maximal solutions; however, since $-\check Y^{x,a,T}$ is a minimal solution to a certain BSDE with nonpositive jumps, we can apply Theorem 2.2 to $-\check Y^{x,a,T}$)
\[
\check Y_0^{x,a,T} \ = \ \inf_{\nu\in\Vc} \E^\nu\bigg[ (h-\phi)(X_T^{x,a}) + \int_0^T\big(f\big(X_s^{x,a},I_s^a,\frac{v(T-s,X_s^{x,a})}{T-s+1}\big) - f(X_s^{x,a},I_s^a,\lambda)\big)ds \bigg].
\]
From the Lipschitz property of $h$ and $\phi$, and estimate \eqref{estimate_x^2}, we have
\begin{equation}
\label{E:LowerBound}
\inf_{\nu\in\Vc} \E^\nu\big[ (h-\phi)(X_T^{x,a})\big] \ \geq \ -C(1+|x|).
\end{equation}
Moreover, from the uniform Lipschitz condition on $f(x,a,y)$ with respect to $y$, and the nondecreasing property of $y$ $\mapsto$ $f(x,a,y)$  in {\bf (H2)}, there exists some adapted, nonpositive, and bounded process 
$\zeta$  such that
\beqs
f\big(X_s^{x,a},I_s^a,\frac{v(T-s,X_s^{x,a})}{T-s+1}\big) - f(X_s^{x,a},I_s^a,\lambda) &=& \zeta_s \Big( \frac{v(T-s,X_s^{x,a})}{T-s+1} - \lambda \Big),
\enqs
for all $0\leq s\leq T$. Therefore, we have
\begin{align}
\label{E:LowerBound2}
&\inf_{\nu\in\Vc} \E^\nu\bigg[\int_0^T\big(f\big(X_s^{x,a},I_s^a,\frac{v(T-s,X_s^{x,a})}{T-s+1}\big) - f(X_s^{x,a},I_s^a,\lambda)\big)ds \bigg] \notag \\
&= \ \inf_{\nu\in\Vc} \E^\nu\bigg[ \int_0^T\zeta_s \Big( \frac{v(T-s,X_s^{x,a})}{T-s+1} - \lambda \Big) ds \bigg].
\end{align}
From step 1 and the Lipschitz property of $\phi$, it follows that $v(T,x)-\lambda T\leq C(1+|x|)$, and consequently $\frac{v(T,x)}{T+1}-\lambda\leq C(1+|x|)$. Hence, since $\zeta$ is nonpositive,
\[
\zeta_s\Big(\frac{v(T-s,X_s^{x,a})}{T-s+1} - \lambda\Big) \ \geq \ \zeta_s C(1+|X_s^{x,a}|) \ \geq \ -L_2C(1+|X_s^{x,a}|),
\]
where we used the fact that $\zeta$ is bounded by $L_2$, the Lipschitz constant of $f$. Plugging the above estimate into \eqref{E:LowerBound2} combined with \eqref{estimate_x^2}, and recalling \eqref{E:LowerBound}, we find
\[
\check Y_0^{x,a,T} \ \geq \ -C(1+|x|).
\]
From \reff{wtildeY} and \reff{YtildeYcheck}, we conclude that $w(T,x)\geq -C(1+|x|)$.
\ep

\begin{Remark} \label{remunilambda}
{\rm
(i) From Theorem \ref{T:Convergence} and Remark \ref{R:MinSol-->ErgHJB}, we deduce a uniqueness result for the component $\lambda$ of a minimal solution to the ergodic BSDE with nonpositive jumps \eqref{ergodicbsde}-\eqref{JumpConErgodic}. Indeed, consider a family of minimal solutions to \eqref{ergodicbsde}-\eqref{JumpConErgodic} as in Remark \ref{R:MinSol-->ErgHJB}. Then, from Theorem \ref{T:Convergence} we see that $\lambda$ is given by \eqref{convergonebis}.

\noindent (ii) Let $f\in{\bf L_{\textup{loc}}^1}([0,\infty);\R)$ be such that $\int_0^\infty e^{-\beta t}f(t)dt$ exists for $\beta>0$. A theorem which states that, under certain conditions on $f$, if $\lim_{\beta\rightarrow0^+}\beta\int_0^\infty e^{-\beta t}f(t)dt=f_\infty\in\R$ then $\lim_{T\rightarrow\infty}\frac{1}{T}\int_0^Tf(t)dt=f_\infty$, is called a \emph{Tauberian theorem}, see, e.g., \cite{arendtbatty95}. In our paper, Theorem \ref{T:ErgodicVisc} together with Theorem \ref{T:Convergence} can be thought as a ``robust'' Tauberian theorem. Indeed, in Theorem \ref{T:ErgodicVisc} we proved that the convergence of $\beta_k v^{\beta_k}(x)=\beta_k\sup_\alpha\E_x[\int_0^\infty e^{-\beta_k t}f(X_t^\alpha,\alpha_t)dt]$ towards $\lambda$, and also of $v^{\beta_k}(x)-v^{\beta_k}(0)$ towards $\phi$, allows to construct a viscosity solution $(\lambda,\phi)$ to the ergodic equation \eqref{ergodicHJB}. Then, Theorem \ref{T:Convergence} implies the convergence of $\frac{v(T,x)}{T}=\frac{1}{T}\sup_\alpha\E_x[\int_0^T f(X_t^\alpha,\alpha_t)dt + h(X_T^\alpha)]$ towards $\lambda$.

\ep
}
\end{Remark}

\subsection{Further convergence result via verification: the proof of (\ref{convergtwo})}

We conclude this section by presenting, in the form of a verification theorem, the following result, which shows the validity of the convergence \eqref{convergtwo}.

\begin{Theorem}
\label{T:Verification}
Let  Assumptions {\bf (H1)}, {\bf (H2)}, {\bf (H$A$)}, and {\bf (H$\vartheta$)} hold.  Suppose that:
\begin{enumerate}
\item[\textup{(i)}] $(\lambda,\phi)$, with  $\phi\in C^2(\R^d)$ and  Lipschitz,  is a classical solution to the ergodic equation \eqref{ergodicHJB}.
\item[\textup{(ii)}] In the ergodic equation \eqref{ergodicHJB},  the supremum is attained at $a=\underline\alpha(x)$, for every $x\in\R^d$, for some function $\underline\alpha\colon A\rightarrow\R^d$ such that $\underline b,\underline\sigma$ in \eqref{underline_b_sigma} satisfy Assumption {\bf (H1)}.
\end{enumerate}
Consider the unique (viscosity) solution $v$ to \eqref{v}. Then, there exists a real  constant $c$ such that
\[
v(T,x)-(\lambda{T}+\phi(x)) \ \overset{T\rightarrow\infty}{\longrightarrow} \ c,
\]
for all $x\in\textup{supp}\,\rho$, the support of the invariant  measure $\rho$ given by  Proposition \ref{P:LawInvariant}. In particular, when $\textup{supp}\,\rho=\R^d$ we deduce that $\phi$  is uniquely determined up to a constant. 
\end{Theorem}

\begin{Remark}
{\rm  The existence of a smooth solution $(\lambda,\phi)$ to the ergodic equation \eqref{ergodicHJB} is ensured  under a uniform ellipticity condition, see, e.g., Theorem 1.1 in \cite{safonov88}. More precisely, suppose that assumptions {\bf (H1)}, {\bf (H2)}, {\bf (H$A$)}, and {\bf (H$\vartheta$)} hold and let $(\lambda,\phi)$ be a viscosity solution to \eqref{ergodicHJB} with $\phi$ Lipschitz, whose existence follows from Theorem \ref{T:ErgodicVisc}. Then, to exploit Theorem 1.1 in \cite{safonov88}, we fix $\delta>0$ and we consider the elliptic equation in the unknown $\psi$ on the bounded domain $B_R\subset\R^d$ (the open ball of radius $R>0$ centered at the origin)
\beq
\delta\psi(x) - \sup_{a\in A}\big[\Lc^a\psi + f(x,a,\lambda) - \lambda + \delta\phi(x)\big] & = & 0, \;\;\; \mbox{ on } \; B_R, \label{ergodicHJB_Ball} \\
\psi & = & \phi, \;\;\; \mbox{ on } \; \partial B_R. \label{ergodicHJB_Ball_Boundary}
\enq
Notice that, thanks to the presence of the term ``$\delta\psi(x)$'' in equation \eqref{ergodicHJB_Ball}, we can apply comparison Theorem 3.3 in \cite{ishii89}, from which it follows that $\phi$ is the unique uniformly continuous viscosity solution to equation \eqref{ergodicHJB_Ball}-\eqref{ergodicHJB_Ball_Boundary}. Let us now impose the following uniform ellipticity condition: there exists $\nu\in(0,1]$, possibly depending on $R$, such that
\[
\nu|\xi|^2 \ \leq \ \sum_{i,j=1}^d (\sigma\sigma\trans)_{ij}(x,a)\xi_i\xi_j \ \leq \ \nu^{-1}|\xi|^2, \qquad \forall\,\xi\in\R^d,
\]
for all $x\in\overline B_R$ and $a\in A$. Then, as explained in Remark 1.1 of \cite{safonov88}, under the above assumption,  Theorem 1.1 in \cite{safonov88} holds, and   there exists a unique solution 
$\psi\in C^2(B_R)\cap C(\overline B_R)$ to equation \eqref{ergodicHJB_Ball}-\eqref{ergodicHJB_Ball_Boundary}. Theorem 3.3 in \cite{ishii89} implies that $\psi$ coincides with our function $\phi$, so that 
$\phi\in C^2(B_R)$. Since $R$ is arbitrary, we conclude that $\phi\in C^2(\R^d)$.
\ep
}
\end{Remark}

\noindent \textbf{Proof.}
\emph{Step 1.} Notice that, for any $T,S>0$ and for all $x\in\R^d$, we have
\beq
\label{E:StochControlRepr}
v(T+S,x)  &= & \sup_{\alpha\in\Ac}\,\E\bigg[\int_0^{T} f\big(X_s^{x,\alpha},\alpha_s,\frac{v(T+S-s,X_s^{x,\alpha})}{T+S-s+1}\big) ds 
 + v(S,X_{T}^{x,\alpha})\bigg],
\enq
where $\Ac$ is the set of adapted control processes valued in $A$,  and $(X_t^{x,\alpha})_{t\geq0}$ is the unique solution to the controlled equation \eqref{controlX} starting from $x$ at time $0$. As a matter of fact, to prove \eqref{E:StochControlRepr} we recall from Remark \ref{R:v^T} that $v(T+S,x)=v^{T+S}(0,x)$, for all $x\in\R^d$, where $v^{T+S}$ is the unique uniformly continuous viscosity solution to the following Hamilton-Jacobi-Bellman equation in the unknown 
$\tilde v^{T+S}$:
\beqs
- \Dt{\tilde v^{T+S}}  - \sup_{a\in A}\big[ \Lc^a \tilde v^{T+S}(t,x)  + \; f \big(x,a,\frac{v^{T+S}}{T+S-t+1}\big)\big] &=& 0,  \qquad (t,x)\in[0,T]\times\R^d, \\
\tilde v^{T+S}(T,x) &=& v^{T+S}(T,x),  \qquad x\in\R^d.
\enqs
As a consequence, $v^{T+S}$ is given, for all $(t,x)\in[0,T]\times\R^d$, by the stochastic control representation:
\beqs
v^{T+S}(t,x) &= &  \sup_{\alpha \in \Ac}\,\E\bigg[\int_0^{T-t}f\big(X_s^{x,\alpha},\alpha_s,\frac{ v^{T+S}(s+t,X_s^{x,\alpha})}{T+S-s-t+1}\big) ds 
 + v^{T+S}(T,X_{T-t}^{x,\alpha})\bigg],
\enqs
which implies \eqref{E:StochControlRepr}. In particular, we have
\beqs
v(T+S,x) &\geq &  \E\bigg[\int_0^{T} f\big(X_s^{x,\underline\alpha},\underline\alpha(X_s^{x,\underline\alpha}),
\frac{v(T+S-s,X_s^{x,\underline\alpha})}{T+S-s+1}\big) ds  + v(S,X_{T}^{x,\underline\alpha})\bigg].
\enqs
On the other hand, applying It\^o's formula to $\phi(X_s^{x,\underline\alpha})$ between $0$ and $T$, and using the optimality of $\underline\alpha$ in the ergodic equation \eqref{ergodicHJB}, 
we obtain
\beqs
\phi(x) &=&  \E\bigg[ \int_0^{T} f(X_s^{x,\underline\alpha},\underline\alpha(X_s^{x,\underline\alpha}),\lambda) ds -  \lambda T +  \phi(X_{T}^{x,\underline\alpha}) \bigg].
\enqs
Therefore, $w$ in \reff{w} satisfies
\beq
w(T+S,x) &\geq &  \E\bigg[w(S,X_T^{x,\underline\alpha}) \label{wstep1}   \\
& & \;\;\; + \;  \int_0^T f\big(X_s^{x,\underline\alpha},\underline\alpha(X_s^{x,\underline\alpha}),\frac{ v(T+S-s,X_s^{x,\underline\alpha})}{T+S-s+1} \big) 
- f(X_s^{x,\underline\alpha},\underline\alpha(X_s^{x,\underline\alpha}),\lambda) ds\bigg]. \nonumber
\enq
\emph{Step 2.} Let us prove that there exists a positive constant $C$ such that
\[
|w(T,x) - w(T,x')| \ \leq \ C|x-x'|,
\]
for all $T\geq0$ and $x,x'\in\R^d$. Recalling \reff{w} and since $\phi$ is Lipschitz, it is therefore enough to prove that  the function $v$ satisfies 
\begin{equation}
\label{vsipTlip}
|v(T,x)-v(T,x')| \ \leq \ C |x-x'|, \qquad T\geq 0,\,x\in\R^d,
\end{equation}
for some positive constant $C$. We know that $v(T,x)=Y_0^{x,a,T}$ is represented by the minimal solution to the BSDE with nonpositive jumps \eqref{BSDE_v}-\eqref{JumpConstraint_v} on $[0,T]$. We recall from Theorem 2.1 in \cite{khapha12} that
 $Y^{x,a,T,n}\uparrow Y^{x,a,T}$, where
$(Y^{x,a,T,n},Z^{x,a,T,n},U^{x,a,T,n})\in\mathbf{S^2(0,T)}\times\mathbf{L^2(W;0,T)}\times\mathbf{L^2(\tilde\mu;0,T)}$ is the solution to the penalized BSDE on $[0,T]$:
\begin{align*}
Y_t^{x,a,T,n} \ &= \ h(X_T^{x,a}) + \int_t^T f(X_s^{x,a},I_s^a,\frac{Y_s^{x,a,T,n}}{T-s+1}) ds + n  \int_t^T \int_A \big(U_s^{x,a,T,n}(a')\big)_+\vartheta(da')ds \notag \\
&\quad \ - \int_t^T Z_s^{x,a,T,n} dW_s - \int_t^T\int_A U_s^{x,a,T,n}(a') \tilde\mu(ds,da'), \qquad 0\leq t\leq T
\end{align*}
and
\[
U_t^{x,a,T,n}(a') \ \leq \ 0, \qquad dt\otimes d\P\otimes \vartheta(da')\text{-a.e.}
\]
Then, (\ref{vsipTlip}) follows once we get: 
\[
|Y^{x,a,T,n}_0-Y^{x',a,T,n}_0| \ \leq \ C |x-x'|,
\]
for a constant $C$ that does not depend on $x$, $a$, $T$, and $n$. This
can be done using Girsanov theorem and the dissipativity condition \eqref{dissipative} in the same way as for \eqref{v^beta,n_Lipschitz}. \\
\emph{Step 3.} Now, we proceed as in \cite{ichihara12} and we introduce the set $\Gamma$ which contains all the $\omega$-limits of the family $\{w(T,\cdot)\}_{T>1}$ in $C(\R^d)$ (we endow $C(\R^d)$ with the topology for which $f_j\rightarrow f$ in $C(\R^d)$ if and only if $f_j$ converges uniformly to $f$ on any compact subset of $\R^d$). In other words, $\Gamma$ is given by
\[
\Gamma \ := \ \big\{w_\infty\in C(\R^d)\colon\,w(T_j,\cdot)\rightarrow w_\infty\text{ in }C(\R^d)\text{ for some }(T_j)_{j\in\N}\text{ with }T_j\rightarrow\infty\big\}.
\]
It follows from step 2 that the family $\{w(T,\cdot)\}_{T>1}$ is relatively compact in $C(\R^d)$. In particular, $\Gamma\neq\emptyset$, and any $w_\infty$ in $\Gamma$ is Lipschitz. 
To conclude, it suffices to prove that every $w_\infty\in\Gamma$ is equal to the same constant $c\in\R$ on $\text{supp}\,\rho$. 

\emph{- Step 3a.} We first show that any element of $\Gamma$ is constant on $\text{supp}\,\rho$. Let $w_\infty\in\Gamma$, therefore there exists a sequence $(T_j)_{j\in\N}$, with $T_j\rightarrow\infty$, such that $w(T_j,\cdot)\rightarrow w_\infty$ in $C(\R^d)$ as $j\rightarrow\infty$. From \reff{wstep1} with $S$ $=$ $T_j-T$, we have
\beq
w(T_j,x) &\geq &  \E\bigg[w(T_j-T,X_T^{x,\underline\alpha}) \label{E:Limit_w_infty_proof}   \\
& & \;\;\; + \;  \int_0^T f\big(X_s^{x,\underline\alpha},\underline\alpha(X_s^{x,\underline\alpha}),\frac{ v(T_j-s,X_s^{x,\underline\alpha})}{T_j-s+1} \big) 
- f(X_s^{x,\underline\alpha},\underline\alpha(X_s^{x,\underline\alpha}),\lambda) ds\bigg]. \nonumber
\enq 
From \eqref{convergonebis} we have for all $s$ $\in$ $[0,T]$, 
\beqs
\label{E:Limit_v-lambda_proof}
\frac{v(T_j-s,X_s^{x,\underline\alpha})}{T_j-s+1} & \overset{j\rightarrow\infty}{\longrightarrow}  & \lambda, \qquad \P\text{-a.s.}
\enqs
Therefore, sending $j\rightarrow\infty$ in \eqref{E:Limit_w_infty_proof}, and by the dominated convergence theorem,  we obtain
\beqs
w_\infty(x)  &\geq & \E\big[ w_\infty(X_T^{x,\underline\alpha})\big].
\enqs
Moreover, choosing $T:=T_j$ and letting $j\rightarrow\infty$ in the above inequality, we obtain, from Proposition \ref{P:LawInvariant}:
\[
w_\infty(x) \ \geq \ \int_{\R^d} w_\infty(x')\rho(dx').
\]
Now, taking the infimum with respect to $x\in\R^d$, we end up with
\[
0 \ \geq \ \int_{\R^d}(w_\infty(x')-\inf_{\R^d}w_\infty)\rho(dx') \ \geq \ 0.
\]
As a consequence $w_\infty=\inf_{\R^d}w_\infty$, $\rho$-a.s.,
therefore $w_\infty$ is constant on $\text{supp}\,\rho$.

\emph{- Step 3b.} We next prove that every $w_\infty\in\Gamma$ is equal to the same constant $c$ on $\text{supp}\,\rho$. Proceeding as in the derivation of \eqref{wstep1}, we have:
\beq 
& & \E\big[ w(R+S,X_{T-R}^{x,\underline\alpha}\big]   \label{E:Limit_w_infty_proof_bis}  \\
&\geq &  \E\bigg[w(S,X_T^{x,\underline\alpha})    + \;  \int_R^T f\big(X_s^{x,\underline\alpha},\underline\alpha(X_s^{x,\underline\alpha}),\frac{ v(T+S-s,X_s^{x,\underline\alpha})}{T+S-s+1} \big) 
- f(X_s^{x,\underline\alpha},\underline\alpha(X_s^{x,\underline\alpha}),\lambda) ds\bigg]. \nonumber
\enq 
for any $T,S,R>0$ with $R\leq T$. Suppose that there exist two real constants $c_1,c_2$ and two diverging sequences $(T_j)_{j\in\N}$ and $(S_j)_{j\in\N}$ such that $w(T_j,\cdot)\rightarrow c_1$ and $w(S_j,\cdot)\rightarrow c_2$ on $\text{supp}\,\rho$ as $j\rightarrow\infty$. Let us take $T:=T_j$, $R=T_j-S$, and $S:=S_k$ in \eqref{E:Limit_w_infty_proof_bis}. Then, letting $j\rightarrow\infty$ we obtain (notice that, by \eqref{E:Limit_v-lambda_proof} and Lebesgue's dominated convergence theorem, the two integral terms in \eqref{E:Limit_w_infty_proof_bis} simplify one with the other as $j\rightarrow\infty$)
\beqs
c_1 &\geq&  \int_{\R^d}w(S_k,x')\rho(dx').
\enqs
Now, sending $k\rightarrow\infty$, we find
\beqs
c_1 & \geq & \lim_{k\rightarrow\infty}\int_{\R^d}w(S_k,x')\rho(dx') \ = \ \int_{\R^d}C_2\rho(dx') \ = \ c_2.
\enqs
Therefore $c_1\geq c_2$. Changing the role of $(T_j)_{j\in\N}$ and $(S_j)_{j\in\N}$, we also find $c_2\geq c_1$. Hence, $c_1=c_2$ and the claim follows.
\ep

\begin{Remark}
{\rm
Under the hypotheses of Theorem \ref{T:Verification} and when $f=f(x,a)$ does not depend on $y$, $\lambda$ can be interpreted as value of an \emph{ergodic control problem} with gain functional
\beqs
J(x,\alpha) &:=& \limsup_{T\rightarrow\infty}\E\Big[ \int_0^T f(X_t^{x,\alpha},\alpha_t) dt + h(X_T^{x,\alpha}) \Big],
\enqs
where $X^{x,\alpha}$ is the controlled diffusion process satisfying \eqref{controlX}, starting from $x\in\R^d$ at time $0$, and $\alpha$ $\in$ $\Ac$ is a control process, i.e., an $A$-valued adapted process. More precisely, it is clear that $J(x,\alpha)$ depends only on the asymptotic behavior of the trajectories of $X^{x,\alpha}$. Therefore, from the ergodicity of $X^{x,\alpha}$, we expect that there exists a real number $\lambda^*$, independent of $x\in\R^d$, such that
\beqs
\lambda^* &:=& \sup_{\alpha\in\Ac}J(x,\alpha), \qquad \forall\,x\in\R^d,
\enqs
namely, $\lambda^*$ is the value of the ergodic control problem. Let us prove that $\lambda^*=\lambda$. Firtsly, observe that, since $f$ does not depend on $y$, the function $v$ in \eqref{v} admits the stochastic control representation
\beqs
v(T,x) &=& \sup_{\alpha\in\Ac}\E\Big[ \int_0^T f(X_t^{x,\alpha},\alpha_t) dt + h(X_T^{x,\alpha}) \Big].
\enqs
From \eqref{convergonebis} we know that, for any $x\in\R^d$,
\beq
\label{lambda=lambda^*}
\lambda \;\,=\;\, \lim_{T\rightarrow\infty}\frac{v(T,x)}{T} &=& \lim_{T\rightarrow\infty}\sup_{\alpha\in\Ac}\frac{1}{T}\E\Big[ \int_0^T f(X_t^{x,\alpha},\alpha_t) dt + h(X_T^{x,\alpha}) \Big] \notag \\
&=& \lim_{T\rightarrow\infty}\sup_{\alpha\in\Ac}\frac{1}{T}\E\Big[ \int_0^T f(X_t^{x,\alpha},\alpha_t) dt \Big],
\enq
where the last equality follows from the fact that $\lim_{T\rightarrow\infty}\sup_{\alpha\in\Ac}\frac{1}{T}\E[h(X_T^{x,\alpha})]=0$, which is a consequence of the Lipschitz property of $h$ and estimate \eqref{estimate_x^2}.
From \eqref{lambda=lambda^*} we see that $\lambda^*\leq\lambda$. To prove the reverse inequality, fix $x\in\R^d$, then, applying It\^o's formula to $\phi(X_t^{x,\underline\alpha})$ between $0$ and $T$, and using the optimality of $\underline\alpha$ in the ergodic equation \eqref{ergodicHJB}, we obtain
\beqs
\lambda &=& \frac{1}{T}\E\Big[ \int_0^{T} f(X_t^{x,\underline\alpha},\underline\alpha(X_t^{x,\underline\alpha})) dt +  \phi(X_{T}^{x,\underline\alpha}) - \phi(x)\Big].
\enqs
From the Lipschitz property of $\phi$ and estimate \eqref{estimate_x^2}, we have $\frac{1}{T}\E[\phi(X_{T}^{x,\underline\alpha}) - \phi(x)]\rightarrow0$ as $T\rightarrow\infty$, therefore
\beqs
\lambda &=& \lim_{T\rightarrow\infty}\frac{1}{T}\E\Big[ \int_0^{T} f(X_t^{x,\underline\alpha},\underline\alpha(X_t^{x,\underline\alpha})) dt\Big] \\
&\leq & \sup_{\alpha\in\Ac}\Big\{\limsup_{T\rightarrow\infty}\frac{1}{T}\E\Big[ \int_0^{T} f(X_t^{x,\alpha},\alpha_t) dt\Big]\Big\} \;\,=\;\, \lambda^*,
\enqs
which implies that $\lambda^*=\lambda$.
\ep
}
\end{Remark}

\appendix

\setcounter{equation}{0} \setcounter{Assumption}{0}
\setcounter{Theorem}{0} \setcounter{Proposition}{0}
\setcounter{Corollary}{0} \setcounter{Lemma}{0}
\setcounter{Definition}{0} \setcounter{Remark}{0}

\renewcommand\thesection{Appendix}

\section{}

\renewcommand\thesection{\Alph{subsection}}

\renewcommand\thesubsection{\Alph{subsection}}

\subsection{Ergodicity proofs}

\subsubsection{Proof of Lemma \ref{estimX}}

$\bullet$ {\bf Proof of (i)}\\
Let $t\geq0$, then an application of It\^o's formula to $e^{\gamma s}|X_s^{x,a}|^2$ between $0$ and $t$ yields
\begin{align*}
e^{\gamma t}|X_t^{x,a}|^2 \ &= \ |x|^2 + \gamma\int_0^t e^{\gamma s}|X_s^{x,a}|^2 ds + 2\int_0^t e^{\gamma s}X_s^{x,a}.b(X_s^{x,a},I_s^a) ds \\
&\quad \ + \int_0^t e^{\gamma s}\|\sigma(X_s^{x,a},I_s^a)\|^2ds + 2\int_0^t e^{\gamma s}(X_s^{x,a})\trans\sigma(X_s^{x,a},I_s^a)dW_s. \notag
\end{align*}
Rearranging the terms in a suitable way so to exploit the dissipativity condition {\bf (H1)}(ii), we obtain
\begin{align*}
&e^{\gamma t}|X_t^{x,a}|^2 \ = \ |x|^2 + \gamma\int_0^t e^{\gamma s}|X_s^{x,a}|^2 ds + 2\int_0^t e^{\gamma s}X_s^{x,a}.(b(X_s^{x,a},I_s^a)-b(0,I_s^a)) ds \notag \\
& + \int_0^t e^{\gamma s}\text{tr}\big[\big(\sigma(X_s^{x,a},I_s^a)-\sigma(0,I_s^a)\big)\big(\sigma(X_s^{x,a},I_s^a)-\sigma(0,I_s^a)\big)\trans\big]ds \notag \\
& + 2\int_0^t e^{\gamma s}X_s^{x,a}.b(0,I_s^a) ds +  2\int_0^t e^{\gamma s}\text{tr}\big[\sigma(0,I_s^a)\big(\sigma(X_s^{x,a},I_s^a)-\sigma(0,I_s^a)\big)\trans\big]ds \notag \\
& + \int_0^t e^{\gamma s}\text{tr}\big[\sigma(0,I_s^a)\sigma(0,I_s^a)\trans\big]ds + 2\int_0^t e^{\gamma s}(X_s^{x,a})\trans\sigma(X_s^{x,a},I_s^a)dW_s.
\end{align*}
Using {\bf (H1)}, we find
\begin{align*}
e^{\gamma t}|X_t^{x,a}|^2 \ &\leq \ |x|^2 + \gamma\int_0^t e^{\gamma s}|X_s^{x,a}|^2 ds - 2\gamma\int_0^t e^{\gamma s}|X_s^{x,a}|^2 ds + 2M_1\int_0^t e^{\gamma s}|X_s^{x,a}|ds \\
&\quad \ + 2M_1L_1\int_0^t e^{\gamma s}|X_s^{x,a}|ds + M_1^2\int_0^t e^{\gamma s}ds + 2\int_0^t e^{\gamma s}(X_s^{x,a})\trans\sigma(X_s^{x,a},I_s^a)dW_s,
\end{align*}
where $M_1:=\sup_{a\in A}(|b(0,a)|+\|\sigma(0,a)\|)$. From the inequality $|X_s^{x,a}|\leq\eps|X_s^{x,a}|^2+1/(4\eps)$, for any $\eps>0$, we obtain
\begin{align}
\label{E:ItoProof4bis}
e^{\gamma t}|X_t^{x,a}|^2 \ &\leq \ |x|^2 - (\gamma-2M_1\eps-2M_1L_1\eps)\int_0^t e^{\gamma s}|X_s^{x,a}|^2ds \notag \\
&\quad \ + \bigg(\frac{M_1+M_1L_1}{2\eps}+M_1^2\bigg)\int_0^t e^{\gamma s}ds + 2\int_0^t e^{\gamma s}(X_s^{x,a})\trans\sigma(X_s^{x,a},I_s^a)dW_s.
\end{align}
We can find $\eps$ such that $\gamma-2M_1\eps-2M_1L_1\eps\geq0$ (more precisely, if $M_1=0$ then $\eps$ can be any positive real number; otherwise we take $\eps\leq \gamma/(2M_1+2M_1L_1)$), therefore  (also multiplying both sides in \eqref{E:ItoProof4bis} by $e^{-\gamma t}$)
\begin{align}
\label{E:ItoProof2}
|X_t^{x,a}|^2 \ &\leq \ e^{-\gamma t}|x|^2 + \bigg(\frac{M_1+M_1L_1}{2\eps}+M_1^2\bigg)\frac{1-e^{-\gamma t}}{\gamma} \notag \\
&\quad \ + 2e^{-\gamma t}\int_0^t e^{\gamma s}(X_s^{x,a})\trans\sigma(X_s^{x,a},I_s^a)dW_s.
\end{align}
Now, consider $\nu\in\Vc$ and recall that $W$ remains a Brownian motion under $\P^\nu$. Then, the following well-known estimate holds under {\bf (H1)}(i): for all $T>0$ and $p\geq1$, there exists some positive constant $\bar C_{T,p}$ such that
\begin{equation}
\label{estimate2}
\sup_{\nu\in\Vc}\E^\nu\Big[\sup_{0\leq s\leq T}|X_s^{x,a}|^p\Big] \ \leq \ \bar C_{T,p}\big(1+|x|^p\big), \qquad \forall\,(x,a)\in\R^d\times A.
\end{equation}
Estimate \eqref{estimate2} implies that the local martingale
\beqs
(M_T)_{T\geq0} &:=& \bigg(\int_0^T e^{\gamma s}(X_s^{x,a})\trans\sigma(X_s^{x,a},I_s^a)dW_s\bigg)_{T\geq0}
\enqs
is indeed a $\P^\nu$-martingale. Then, we have $\E^\nu[e^{-\gamma t}M_t]=0$. Therefore, taking the expectation $\E^\nu$ in \eqref{E:ItoProof2}, we find
\begin{align*}
\E^\nu\big[|X_t^{x,a}|^2\big] \ &\leq \ |x|^2 + \bigg(\frac{M_1+M_1L_1}{2\eps}+M_1^2\bigg)\frac{1}{\gamma},
\end{align*}
from which we deduce \eqref{estimate_x^2} with $C:=\sqrt{\max\{1,[(M_1+M_1L_1)/(2\eps)+M_1^2]/\gamma\}}$.

\vspace{3mm}

\noindent $\bullet$ {\bf Proof of (ii)}\\
Applying It\^o's formula to $|X_t^{x,a}-X_t^{x',a}|^2$ we find
\begin{align}
\label{E:ItoProof1}
|X_t^{x,a}-X_t^{x',a}|^2 \ &= \ |x-x'|^2 + 2\int_0^t (X_s^{x,a}-X_s^{x',a}).(b(X_s^{x,a},I_s^a)-b(X_s^{x',a},I_s^a)) ds \notag \\
&\quad \ + \int_0^t \|\sigma(X_s^{x,a},I_s^a)-\sigma(X_s^{x',a},I_s^a)\|^2 ds \\
&\quad \ + 2\int_0^t (X_s^{x,a}-X_s^{x',a})\trans(\sigma(X_s^{x,a},I_s^a)-\sigma(X_s^{x',a},I_s^a))dW_s. \notag
\end{align}
Now, consider $\nu\in\Vc$ and recall that $W$ remains a Brownian motion under $\P^\nu$. Using estimate \eqref{estimate2}, we see that the local martingale
\[
\bigg(\int_0^t (X_s^{x,a}-X_s^{x',a})\trans(\sigma(X_s^{x,a},I_s^a)-\sigma(X_s^{x',a},I_s^a))dW_s\bigg)_{t\geq0}
\]
is indeed a $\P^\nu$-martingale. Therefore, taking the expectation $\E^\nu$ with respect to $\P^\nu$ in \eqref{E:ItoProof1} and using the dissipativity condition \eqref{dissipative}, we obtain
\begin{align*}
\E^\nu\big[|X_t^{x,a}-X_t^{x',a}|^2\big] \ &\leq \ |x-x'|^2 - 2\gamma\int_0^t \E^\nu\big[|X_s^{x,a}-X_s^{x',a}|^2\big] ds,
\end{align*}
which implies
\[
\E^\nu\big[|X_t^{x,a}-X_t^{x',a}|^2\big] \ \leq \ |x-x'|^2e^{-2\gamma t}.
\]
\ep

\subsubsection{Proof of Proposition \ref{P:LawInvariant}}

\emph{Step 1. Existence and uniqueness of $\rho$.} Let $\tilde W=(\tilde W_t)_{t\geq0}$ be a $d$-dimensional Brownian motion, independent of $W$ and $\mu$. Then, we define
\[
\bar W_t \ = \
\begin{cases}
W_t, \qquad &t\geq0, \\
\tilde W_{-t}, &t<0.
\end{cases}
\]
For any $T\in\R$ and $x\in\R^d$, we denote $X^{T,x}=(X_t^{T,x})_{t\geq T}$ the unique solution to the equation on $[T,\infty)$:
\begin{equation}\label{xmarkov2}
dX_t \ = \ \underline b(X_t)\,dt + \underline \sigma(X_t)\,d\bar W_t, \quad t\geq T, \qquad X_T \ = \ x.
\end{equation}
From the time-homogeneity of equation \eqref{xmarkov2}, it follows the law invariance property $\Lc(X_t^{T,x})$ $=$ $\Lc(X_{t-T}^x)$, for $t\geq T$, where $X^x$ is the solution to \eqref{xmarkov2} starting from $x$ at time $0$.\\
Let $S>T>0$ and $x\in\R^d$, then, applying It\^o's formula to the difference $|X_t^{-S,x}-X_t^{-T,x}|^2$ from $-T$ to $t\in[-T,0]$, we obtain
\begin{align}
\label{E:ErgodicityProof}
|X_t^{-S,x}-X_t^{-T,x}|^2 \ &= \ |X_{-T}^{-S,x}-x|^2 + 2\int_{-T}^t (X_s^{-S,x}-X_s^{-T,x}).(\underline b(X_s^{-S,x})-\underline b(X_s^{-T,x})) ds \notag \\
&\quad \ + \int_{-T}^t \|\underline\sigma(X_s^{-S,x})-\underline\sigma(X_s^{-T,x})\|^2 ds \\
&\quad \ + 2\int_{-T}^t (X_s^{-S,x}-X_s^{-T,x})\trans(\underline\sigma(X_s^{-S,x})-\underline\sigma(X_s^{-T,x}))d\bar W_s. \notag
\end{align}
Taking the expectation and using the dissipativity condition \eqref{dissipative}, we find
\begin{align*}
\E\big[|X_t^{-S,x}-X_t^{-T,x}|^2\big] \ &\leq \ \E\big[|X_{-T}^{-S,x}-x|^2\big] - 2\gamma\int_{-T}^t \E\big[|X_s^{-S,x}-X_s^{-T,x}|^2\big] ds,
\end{align*}
which implies
\begin{align}
\label{X_-S_-T_Cauchy}
\E\big[|X_0^{-S,x}-X_0^{-T,x}|^2\big] \ &\leq \ \E\big[|X_{-T}^{-S,x}-x|^2\big] e^{-2\gamma T}.
\end{align}
Similar to \eqref{estimate_x^2}, we can prove that there exists a positive constant $\bar C$, depending only on the $L_1$, $M_1$, and $\gamma$, such that
\begin{equation}
\label{X_-S_-T_Cauchy2}
\E\big[|X_{-T}^{-S,x}-x|^2\big] \ \leq \ \bar C(1+|x|^2).
\end{equation}
Plugging \eqref{X_-S_-T_Cauchy2} into \eqref{X_-S_-T_Cauchy}, we obtain
\begin{equation}
\label{X_-S_-T_Cauchy4}
\E\big[|X_0^{-S,x}-X_0^{-T,x}|^2\big] \ \leq \ \bar C(1+|x|^2) e^{-2\gamma T}.
\end{equation}
It follows from \eqref{X_-S_-T_Cauchy4} that $(X_0^{-T,x})_{T>0}$ converges, as $T\rightarrow\infty$, to some square integrable random variable $\eta^x$, which a priori depends on $x$. Let $x'\in\R^d$, then applying It\^o's formula to $|X_s^{-T,x}-X_s^{-T,x'}|^2$ between $-T$ and $t\in[-T,0]$, we find
\begin{align*}
|X_t^{-T,x}-X_t^{-T,x'}|^2 \ &= \ |x-x'|^2 + 2\int_{-T}^t (X_s^{-T,x}-X_s^{-T,x'}).(\underline b(X_s^{-T,x})-\underline b(X_s^{-T,x'})) ds \\
&\quad \ + \int_{-T}^t \|\underline\sigma(X_s^{-T,x})-\underline\sigma(X_s^{-T,x'})\|^2 ds \\
&\quad \ + 2\int_{-T}^t (X_s^{-T,x}-X_s^{-T,x'})\trans(\underline\sigma(X_s^{-T,x})-\underline\sigma(X_s^{-T,x'}))d\bar W_s.
\end{align*}
Taking the expectation and using the dissipativity condition \eqref{dissipative}, we obtain
\begin{align*}
\E\big[|X_t^{-T,x}-X_t^{-T,x'}|^2\big] \ &\leq \ |x-x'|^2 - 2\gamma\int_{-T}^t \E\big[|X_s^{-T,x}-X_s^{-T,x'}|^2\big] ds,
\end{align*}
which implies
\begin{align*}
&\E\big[|X_0^{-T,x}-X_0^{-T,x'}|^2\big] \ \leq \ |x-x'|^2 e^{-2\gamma T} \ \overset{T\rightarrow\infty}{\longrightarrow} \ 0.
\end{align*}
As a consequence, $\eta^x=\eta^{x'}=:\eta$. We denote $\rho:=\Lc(\eta)$. Finally, using the law invariance property already recalled, and the fact that convergence in $\mathbf{L^2(}\P\mathbf{)}$ implies convergence in law, we deduce
\begin{equation}
\label{E:LawConv}
\Lc(X_T^x) \ = \ \Lc(X_0^{-T,x}) \ \longrightarrow \ \rho,
\end{equation}
weakly as $T\rightarrow\infty$. From the square integrability of $\eta$ we see that $\int_{\R^d}|x|^2\rho(dx)<\infty$. Let us now prove the invariance property. Let $\varphi\in C_b(\R^d)$, then, from the Markov property we have
\[
P_{t+s}^{\underline{\alpha}}\varphi(x) \ = \ P_t^{\underline{\alpha}}(P_s^{\underline{\alpha}}\varphi)(x), \qquad \forall\,t,s\geq0.
\]
Sending $t\rightarrow\infty$, using \eqref{E:LawConv} and the Feller property, we obtain
\begin{equation}
\label{E:InvarianceProperty_phi_Continuous}
\int_{\R^d}\varphi(x)\rho(dx) \ = \ \int_{\R^d} P_s^{\underline{\alpha}}\varphi(x)\rho(dx), \qquad \forall\,s\geq0,
\end{equation}
for all $\varphi\in C_b(\R^d)$. By a monotone class argument, we see that \eqref{E:InvarianceProperty_phi_Continuous} remains true for all $\varphi\in B(\R^d)$, which implies the invariant property of $\rho$. Concerning the uniqueness of $\rho$, let us consider another invariance probability measure $\nu$ and take $\varphi\in C_b(\R^d)$, then
\[
\int_{\R^d}\varphi(x)\nu(dx) \ = \ \int_{\R^d}P_s^{\underline{\alpha}}\varphi(x)\nu(dx) \ \overset{s\rightarrow\infty}{\longrightarrow} \ \int_{\R^d}\bigg(\int_{\R^d}\varphi(x)\rho(dx)\bigg) \nu(dx) \ = \ \int_{\R^d}\varphi(x)\rho(dx).
\]
Since the result holds for any $\varphi\in C_b(\R^d)$, we deduce the uniqueness of $\rho$.\\
\emph{Step 2. \eqref{psiergodic} is valid for any continuous $\varphi$ satisfying a linear growth condition.} For any $R>0$, consider a continuous function $\chi_R\colon\R^d\rightarrow[0,1]$ which is equal to $1$ on $B_R\subset\R^d$ (the open ball of radius $R$ centered at the origin) and is equal to $0$ on $\R^d\backslash B_{2R}$. Then
\begin{align}
\label{E:LawConv_LinearGrowth1}
\bigg|\E[\varphi(X_t^x)] - \int_{\R^d}\varphi(x)\rho(dx)\bigg| \ &\leq \ \bigg|\E[\varphi(X_t^x)\chi_R(X_t^x)] - \int_{\R^d}\varphi(x)\chi_R(x)\rho(dx)\bigg| \\
&\quad \ + \bigg|\E[\varphi(X_t^x)(1-\chi_R(X_t^x))] - \int_{\R^d}\varphi(x)(1-\chi_R(x))\rho(dx)\bigg|. \notag
\end{align}
Since $\varphi\chi_R\in C_b(\R^d)$, the first term on the right-hand side of \eqref{E:LawConv_LinearGrowth1} goes to zero as $t\rightarrow\infty$. If
\begin{equation}
\label{E:LawConv_LinearGrowth2}
\lim_{R\rightarrow\infty}\limsup_{t\rightarrow\infty}\bigg|\E[\varphi(X_t^x)(1-\chi_R(X_t^x))] - \int_{\R^d}\varphi(x)(1-\chi_R(x))\rho(dx)\bigg| \ = \ 0,
\end{equation}
then, taking first $\limsup_{t\rightarrow\infty}$ and then $\lim_{R\rightarrow\infty}$ in \eqref{E:LawConv_LinearGrowth1}, we get the thesis. Therefore it remains to prove \eqref{E:LawConv_LinearGrowth2}. From Cauchy-Schwarz inequality, the linear growth property of $\varphi$, and estimate \eqref{X_-S_-T_Cauchy2}, we have that there exists a positive constant $C$ such that
\begin{align*}
\big|\E[\varphi(X_t^x)(1-\chi_R(X_t^x))]\big| \ &\leq \ \sqrt{\E[|\varphi(X_t^x)|^2]}\sqrt{\E_x[|1-\chi_R(X_t^x)|^2]} \\
&\leq \ C(1+|x|)\sqrt{\E[|1-\chi_R(X_t^x)|^2]}.
\end{align*}
Since the function $|1-\chi_R|^2\in C_b(\R^d)$, we find
\[
\limsup_{t\rightarrow\infty}\big|\E[\varphi(X_t^x)(1-\chi_R(X_t^x))]\big| \ \leq \ C(1+|x|)\sqrt{\int_{\R^d}|1-\chi_R(x)|^2\rho(dx)}.
\]
Similarly, we have
\begin{align*}
\bigg|\int_{\R^d}\varphi(x)(1-\chi_R(x))\rho(dx)\bigg| \ &\leq \ \sqrt{\int_{\R^d}|\varphi(x)|^2\rho(dx)}\sqrt{\int_{\R^d}|1-\chi_R(x)|^2\rho(dx)} \\
&\leq \ C\bigg(1+\sqrt{\int_{\R^d}|x|^2\rho(dx)}\bigg)\sqrt{\int_{\R^d}|1-\chi_R(x)|^2\rho(dx)},
\end{align*}
where we recall that $\int_{\R^d}|x|^2\rho(dx)<\infty$. In conclusion, we obtain
\begin{align*}
&\limsup_{t\rightarrow\infty}\bigg|\E[\varphi(X_t^x)(1-\chi_R(X_t^x))] - \int_{\R^d}\varphi(x)(1-\chi_R(x))\rho(dx)\bigg| \\
&\leq \ C\bigg(1+|x|+\sqrt{\int_{\R^d}|x|^2\rho(dx)}\bigg)\sqrt{\int_{\R^d}|1-\chi_R(x)|^2\rho(dx)}.
\end{align*}
Notice that
\[
\int_{\R^d}|1-\chi_R(x)|^2\rho(dx) \ \leq \ \int_{\R^d\backslash B_R}\rho(dx) \ = \ \rho(\R^d\backslash B_R) \ \overset{R\rightarrow\infty}{\longrightarrow} \ 0,
\]
which implies \eqref{E:LawConv_LinearGrowth2}.
\ep

\subsection{Elliptic BSDEs}

\subsubsection{Proof of Lemma \ref{estimapriori}}

Let $0\leq t\leq T<\infty$ and apply It\^o's formula to $e^{-2\beta s}|\Delta Y_s|^2$ between $t$ and $T$, then
\begin{align}
\label{penalizedBSDE_Delta}
e^{-2\beta t}|\Delta Y_t|^2 \ &= \ e^{-2\beta T}|\Delta Y_T|^2 + 2n  \int_t^T \int_A e^{-2\beta s}\Delta Y_s\big[(U_s^{1,\beta,n}(a''))_+ - (U_s^{2,\beta,n}(a''))_+\big]\vartheta(da'')ds  \notag \\
&\quad \ + 2 \int_t^T e^{-2\beta s}\Delta Y_s\big(f_1(X_s^{x,a},I_s^a,\beta Y_s^{1,\beta,n}) - f_2(X_s^{x',a'},I_s^{a'},\beta Y_s^{2,\beta,n})\big) ds \notag \\
&\quad \ - 2\int_t^Te^{-2\beta s}\Delta Y_s\Delta Z_s\,dW_s - 2\int_t^T\int_Ae^{-2\beta s}\Delta Y_s\Delta U_s(a'')\,\tilde\mu (ds,da'') \notag \\
&\quad \ - \int_t^T e^{-2\beta s}|\Delta Z_s|^2 ds - \int_t^T\int_A e^{-2\beta s}|\Delta U_s(a'')|^2 \mu(ds,da'').
\end{align}
Notice that, using the nonincreasing property of $f_1$ in $y$, we have
\begin{align*}
\int_t^T e^{-2\beta s}\Delta Y_s\big(f_1(X_s^{x,a},I_s^a,\beta Y_s^{1,\beta,n}) - f_2(X_s^{x',a'},I_s^{a'},\beta Y_s^{2,\beta,n})\big) ds \\
\leq \ \int_t^T e^{-2\beta s}\Delta Y_s(\Delta_s'f_1 + \Delta_s f) ds.
\end{align*}
Now, define the $[1,n+1]$-valued map $\nu$ as follows
\[
\nu_t(a'') \ = \ 1 + n\frac{(U_t^{1,\beta,n}(a''))_+ - (U_t^{2,\beta,n}(a''))_+}{\Delta U_t(a'')}1_{\{\Delta U_t(a'')\neq0\}}, \qquad t\geq0,\,a''\in A.
\]
Observe that $\nu$ is a $\Pc\otimes\Bc(A)$-measurable map satisfying $1\leq \nu_s(a)\leq n+1$, $ds\otimes d\P\otimes\vartheta(da)$-a.e., then $\nu\in\Vc_n$. Let us consider the probability measure $\P^\nu$ equivalent to $\P$ on $(\Omega,\Fc_T)$ with Radon-Nikodym density given by \eqref{Pnu}. Recalling that $\tilde\mu^\nu$ denotes the compensated martingale measure associated to $\mu$ under $\P^\nu$, equation \eqref{penalizedBSDE_Delta} can be rewritten as follows
\begin{align}
\label{penalizedBSDE_Delta2}
&e^{-2\beta t}|\Delta Y_t|^2 + \int_t^T e^{-2\beta s}|\Delta Z_s|^2 ds + \int_t^T\int_A e^{-2\beta s}|\Delta U_s(a'')|^2 \mu(ds,da'') \notag \\
&\leq \ e^{-2\beta T}|\Delta Y_T|^2 + 2 \int_t^T e^{-2\beta s}\Delta Y_s(\Delta_s'f_1 + \Delta_s f) ds - 2\int_t^Te^{-2\beta s}\Delta Y_s\Delta Z_s\,dW_s \notag \\
&\quad \ - 2\int_t^T\int_Ae^{-2\beta s}\Delta Y_s\Delta U_s(a'')\,\tilde\mu^\nu (ds,da'').
\end{align}
From Lemma 2.5  in \cite{khapha12}, we see that the two stochastic integrals on the right-hand side of \eqref{penalizedBSDE_Delta2} are martingales. Hence, taking the expectation $\E^\nu$, conditional on $\Fc_t$, with respect to $\P^\nu$ in \eqref{penalizedBSDE_Delta2}, we end up with estimate \eqref{EstimateDeltaYZU}.
\ep

\subsubsection{Proof of Proposition \ref{P:ExistUniqPen}} 

\emph{Uniqueness.} Fix $(\beta,n)\in(0,\infty)\times\N$ and consider two solutions $(Y^{1,\beta,n},Z^{1,\beta,n},U^{1,\beta,n})$, $(Y^{2,\beta,n},Z^{2,\beta,n},U^{2,\beta,n})\in\mathbf{S_{\textup{loc}}^2}\times\mathbf{L_{\textup{loc}}^2(W)}\times\mathbf{L_{\textup{loc}}^2(\tilde\mu)}$ to \eqref{penalizedbetabsde}. Set $\Delta Y_t$ $=$ $Y_t^{1,\beta,n}-Y_t^{2,\beta,n}$, $\Delta Z_t$ $=$ $Z_t^{1,\beta,n}-Z_t^{2,\beta,n}$, and $\Delta U_t(a')$ $=$ $U_t^{1,\beta,n}(a')-U_t^{2,\beta,n}(a')$, $t\geq 0$, $a'\in A$. Let $0\leq t\leq T<\infty$. Then, from estimate \eqref{EstimateDeltaYZU} with $f_1=f_2=f$ (so that $\Delta' f_1=\Delta f=0$), there exists $\nu$ $\in$ $\Vc_n$ such that 
\begin{equation}
\label{penalizedBSDE_Uniq2}
e^{-2\beta t}|\Delta Y_t|^2 \ \leq \ \E^\nu\big[e^{-2\beta T}|\Delta Y_T|^2\big|\Fc_t\big].
\end{equation}
Moreover, recall from \eqref{estimate_x^2} that the following estimate holds
\[
\E^\nu\big[|X_T^{x,a}|^2\big] \ \leq \ C_{b,\sigma}(1+|x|^2), \qquad \forall\,T\geq 0.
\]
Since $|\Delta Y_T|\leq2C(1+|X_T^{x,a}|)$, we conclude that $\E^\nu[e^{-2\beta T}|\Delta Y_T|^2]\rightarrow0$, as $T\rightarrow\infty$. From \eqref{penalizedBSDE_Uniq2} it follows that $\Delta Y=0$. Finally, plugging $\Delta Y=0$ into \eqref{penalizedBSDE_Delta2}, we conclude that $\Delta Z=0$ and $\Delta U=0$.\\
\emph{Existence. Step 1. Approximating BSDE.} Fix $(x,a,\beta,n)\in\R^d\times\R^q\times(0,\infty)\times\N$, $T>0$, and consider the backward stochastic differential equation on $[0,T]$ given by, $\P$-a.s.,
\begin{align}
Y_t \ &= \ - \beta \int_t^T Y_s\,ds+ n  \int_t^T \int_A (U_s(a'))_+\vartheta(da')\,ds + \int_t^T f(X_s^{x,a},I_s^a,\beta Y_s)\,ds \notag \\
&\quad \ - \int_t^TZ_s\,dW_s - \int_t^T\int_AU_s(a')\,\tilde\mu (ds,da'), \qquad 0\leq t\leq T. \label{penalizedBSDE_T}
\end{align}
Notice that \eqref{penalizedBSDE_T} has a zero terminal condition at the final time $T$. It follows from Lemma 2.4 in \cite{tanli94} that there exists a unique solution $(Y^T,Z^T,U^T)\in{\bf S^2(0,T)}\times{\bf L^2(W;0,T)}\times{\bf L^2(\tilde\mu;0,T)}$ to \eqref{penalizedBSDE_T}.\\
\emph{Step 2. Estimate for $Y^T$.} Let $t\in[0,T]$, then, from estimate \eqref{EstimateDeltaYZU} with $(Y^{1,\beta,n},Z^{1,\beta,n},U^{1,\beta,n})=(Y^T,Z^T,U^T)$, $(Y^{2,\beta,n},Z^{2,\beta,n},U^{2,\beta,n})=(0,0,0)$, $f_1=f$, and $f_2=0$, there exists $\nu\in\Vc_n$:
\begin{align}
\label{estimateY_Proof0}
|Y_t^T|^2 \ &\leq \ 2\int_t^T e^{-2\beta(s-t)}\E^\nu\big[Y_s^Tf(X_s^{x,a},I_s^a,0)\big|\Fc_t\big] ds \notag \\
&\leq \ 2\int_t^T e^{-2\beta(s-t)} \sqrt{\E^\nu\big[|Y_s^T|^2\big|\Fc_t\big]} \sqrt{\E^\nu\big[|f(X_s^{x,a},I_s^a,0)|^2\big|\Fc_t\big]} ds.
\end{align}
Set $g(s)=e^{-2\beta(s-t)}\E^\nu[|Y_s^T|^2|\Fc_t]$ and $h(s)=2e^{-\beta(s-t)}\sqrt{\E^\nu[|f(X_s^{x,a},I_s^a,0)|^2]}$, for any $s\in[t,T]$. Then, recalling that $g(T)=0$, inequality \eqref{estimateY_Proof0} becomes
\[
g(t) \ \leq \ g(T) + \int_t^T \sqrt{g(s)}h(s)ds.
\]
Our aim is to derive a Gronwall type estimate for $g$. To this end, define
\[
\tilde g(t) \ := \ g(T) + \int_t^T \sqrt{g(s)}h(s)ds, \qquad 0\leq t\leq T.
\]
Notice that $\tilde g\in C^1([0,T])$. Moreover $g(t)\leq \tilde g(t)$, for any $0\leq t\leq T$, and
\[
\tilde g(T') \ = \ \tilde g(T'') + \int_{T'}^{T''} \sqrt{g(s)}h(s)ds \ \leq \ \tilde g(T'') + \int_{T'}^{T''} \sqrt{\tilde g(s)}h(s)ds, \quad t\leq T'<T''\leq T.
\]
Dividing by $T''-T'$ and letting $T''-T'\rightarrow0$, we deduce the differential inequality
\[
\tilde g'(s) \ \geq \ -h(s)\sqrt{\tilde g(s)}, \qquad t\leq s\leq T.
\]
We have
\[
\frac{d\sqrt{\tilde g(s)}}{ds} \ \geq \ -\frac{1}{2}h(s),
\]
which yields
\[
\sqrt{\tilde g(T)} - \sqrt{\tilde g(t)} \ \geq \ -\frac{1}{2}\int_t^T h(s) ds.
\]
Therefore, we find
\[
|Y_t^T| \ = \ \sqrt{g(t)} \ \leq \ \sqrt{\tilde g(t)} \ \leq \ \int_t^T e^{-\beta(s-t)}\sqrt{\E^\nu\big[|f(X_s^{x,a},I_s^a,0)|^2\big|\Fc_t\big]} ds.
\]
Recalling that $|f(x,a,0)|\leq L_2|x|+M_2$, with $M_2:=\sup_{a\in A}|f(0,a,0)|$, so that $|f(x,a,0)|^2\leq 2L_2^2|x|^2+2M_2^2$, and using the inequality $\sqrt{a+b}\leq\sqrt{a}+\sqrt{b}$, for any $a,b\in\R_+$, we find
\[
|Y_t^T| \ \leq \ \sqrt{2}L_2\int_t^\infty e^{-\beta(s-t)}\sqrt{\E^\nu\big[|X_s^{x,a}|^2\big|\Fc_t\big]} ds + \sqrt{2}M_2 \int_t^\infty e^{-\beta(s-t)} ds.
\]
From estimate \eqref{estimate_x^2_conditional}, we have
\begin{align}
\label{estimateY_Proof}
|Y_t^T| \ &\leq \ \sqrt{2}\big(L_2\sqrt{C_{b,\sigma}}\big(1+|X_t^{x,a}|\big)+M_2\big)\int_t^\infty e^{-\beta(s-t)} ds \notag \\
&= \ \sqrt{2}\frac{L_2\sqrt{C_{b,\sigma}}(1+|X_t^{x,a}|)+M_2}{\beta}.
\end{align}
\emph{Step 3. Convergence of $(Y^T)_{T>0}$.} Let $T,T'>0$, with $T<T'$, and denote $\Delta Y_t$ $=$ $Y_t^T-Y_t^{T'}$, $0\leq t\leq T$. Let $t\in[0,T]$, then estimate \eqref{EstimateDeltaYZU} reads
\begin{equation}
\label{Y_t^T'-Y_t^T}
|\Delta Y_t|^2 \ \leq \ e^{-2\beta(T-t)}\E^\nu\big[|\Delta Y_T|^2\big|\Fc_t\big] \ \overset{T\rightarrow\infty}{\longrightarrow} \ 0,
\end{equation}
where the convergence result follows from \eqref{estimateY_Proof} and \eqref{estimate_x^2_conditional}. Let us now consider the family of real-valued c\`adl\`ag adapted processes $(Y^T)_{T>0}$. It follows from \eqref{Y_t^T'-Y_t^T} that, for any $t\geq0$, the family $(Y_t^T(\omega))_{T>0}$ is Cauchy for almost every $\omega$, so that it converges $\P$-a.s. to some $\Fc_t$-measurable random variable $Y_t$, which is bounded from the right-hand side of \eqref{estimateY_Proof}. Moreover, using again \eqref{Y_t^T'-Y_t^T}, \eqref{estimateY_Proof}, and \eqref{estimate_x^2_conditional}, we see that, for any $0\leq S<T\wedge T'$, with $T,T'>0$, we have
\begin{equation}
\label{BSDE_penalized_sup-->0}
\sup_{0\leq t \leq S}|Y_t^{T'}-Y_t^T| \ \leq \ e^{-\beta(T\wedge T'-S)}C_0\sup_{0\leq t\leq S}(1+|X_t^{x,a}|) \ \overset{T,T'\rightarrow\infty}{\longrightarrow} \ 0, 
\end{equation}
where $C_0$ is a positive constant independent of $S,T,T'$. In other words, the family $(Y^T)_{T>0}$ converges $\P$-a.s. to $Y$ uniformly on compact subsets of $\R_+$. Since each $Y^T$ is a c\`adl\`ag process, it follows that $Y$ is c\`adl\`ag, as well. Finally, from estimate \eqref{estimateY_Proof} we see that $Y\in\mathbf{S_{\textup{loc}}^2}$ and
\beqs
|Y_t| &\leq & \frac{C}{\beta}\big(1+|X_t^{x,a}|\big), \qquad \forall\,t\geq0.
\enqs
\emph{Step 4. Convergence of $(Z^T,U^T)_{T>0}$.} Let $S,T,T'>0$, with $S<T<T'$. Then, applying It\^o's formula to $e^{-2\beta t}|Y_t^{T'}-Y_t^T|^2$ between $0$ and $S$, and taking the expectation, we find
\begin{align*}
&\E\int_0^S e^{-2\beta s}|Z_s^{T'}-Z_s^T|^2 ds + \E\int_0^S\int_A e^{-2\beta s}|U_s^{T'}(a')-U_s^T(a')|^2 \vartheta(da')ds \\
&= \ e^{-2\beta S}\E\big[|Y_S^{T'}-Y_S^T|^2\big] - |Y_0^{T'}-Y_0^T|^2 \\
&+ 2\E\int_0^S e^{-2\beta s}\big(Y_s^{T'}-Y_s^T\big)\big(f(X_s^{x,a},I_s^a,\beta Y_s^{T'}) -f(X_s^{x,a},I_s^a,\beta Y_s^T)\big)ds \\
&+ 2n\E\int_0^S\int_A e^{-2\beta s}\big(Y_s^{T'}-Y_s^T\big)\big((U_s^{T'}(a'))_+ - (U_s^T(a'))_+\big) \vartheta(da')ds.
\end{align*}
Since the map $y\mapsto f(x,a,y)$ is nonincreasing, we get (using also the inequality $ab\leq a^2/2+b^2/2$, for any $a,b\in\R$)
\begin{align*}
&\E\int_0^S e^{-2\beta s}|Z_s^{T'}-Z_s^T|^2 ds + \E\int_0^S\int_A e^{-2\beta s}|U_s^{T'}(a')-U_s^T(a')|^2 \vartheta(da')ds \\
&\leq \ e^{-2\beta S}\E\big[|Y_S^{T'}-Y_S^T|^2\big] + 2n\E\int_0^S\int_A e^{-2\beta s}|Y_s^{T'}-Y_s^T||(U_s^{T'}(a'))_+ - (U_s^T(a'))_+| \vartheta(da')ds \\
&\leq \ e^{-2\beta S}\E\big[|Y_S^{T'}-Y_S^T|^2\big] + 2n^2\vartheta(A)\E\int_0^S e^{-2\beta s}|Y_s^{T'} - Y_s^T|^2 ds \\
&\quad \ + \frac{1}{2}\E\int_0^S\int_A e^{-2\beta s}|(U_s^{T'}(a'))_+ - (U_s^T(a'))_+|^2 \vartheta(da')ds.
\end{align*}
Multiplying the previous inequality by $e^{2\beta S}$, we obtain
\begin{align*}
&\E\int_0^S |Z_s^{T'}-Z_s^T|^2 ds + \frac{1}{2}\E\int_0^S\int_A |U_s^{T'}(a')-U_s^T(a')|^2 \vartheta(da')ds \\
&\leq \ \E\int_0^S e^{2\beta(S-s)}|Z_s^{T'}-Z_s^T|^2 ds + \frac{1}{2}\E\int_0^S\int_A e^{2\beta(S-s)}|U_s^{T'}(a')-U_s^T(a')|^2 \vartheta(da')ds \\
&\leq \ \E\big[|Y_S^{T'}-Y_S^T|^2\big] + 2n^2\vartheta(A)\E\int_0^S e^{2\beta(S-s)}|Y_s^{T'} - Y_s^T|^2 ds \ \overset{T,T'\rightarrow\infty}{\longrightarrow} \ 0,
\end{align*}
where the convergence to zero follows from estimate \eqref{BSDE_penalized_sup-->0}. Then, for any $S>0$, we see that the family $(Z_{|[0,S]}^T,U_{|[0,S]}^T)_{T>S}$ is Cauchy in the Hilbert space $\mathbf{L^2(W;0,S)}\times\mathbf{L^2(\tilde\mu;0,S)}$. Therefore, we deduce that there exists $(\bar Z^S,\bar U^S)\in\mathbf{L^2(W;0,S)}\times\mathbf{L^2(\tilde\mu;0,S)}$ such that $(Z_{|[0,S]}^T,U_{|[0,S]}^T)_{T>S}$ converges to $(\bar Z^S,\bar U^S)$ in $\mathbf{L^2(W;0,S)}\times\mathbf{L^2(\tilde\mu;0,S)}$, i.e.,
\[
\E\int_0^S |Z_s^T-\bar Z_s^S|^2 ds + \E\int_0^S\int_A |U_s^T(a')-\bar U_s^S(a')|^2\vartheta(da')ds \ \overset{T\rightarrow\infty}{\longrightarrow} \ 0.
\]
Notice that $\bar Z_{|[0,S]}^{S'}=\bar Z^S$ and $\bar U_{|[0,S]}^{S'}=\bar U^S$, for any $0\leq S\leq S'<\infty$. Indeed, $(\bar Z_{|[0,S]}^{S'},\bar U_{|[0,S]}^{S'})$, as $(\bar Z^S,\bar U^S)$, is the limit in $\mathbf{L^2(W;0,S)}\times\mathbf{L^2(\tilde\mu;0,S)}$ of $(Z_{|[0,S]}^T,U_{|[0,S]}^T)_{T>S}$. Hence, we define $Z_s=\bar Z_s^S$ and $U_s=\bar U_s^S$, for all $s\in[0,S]$ and for any $S>0$. Observe that $(Z,U)\in\mathbf{L_{\text{loc}}^2(W)}\times\mathbf{L_{\text{loc}}^2(\tilde\mu)}$. Moreover, for any $S>0$, $(Z_{|[0,S]}^T,U_{|[0,S]}^T)_{T>S}$ converges to $(Z_{|[0,S]},U_{|[0,S]})$ in $\mathbf{L^2(W;0,S)}\times\mathbf{L^2(\tilde\mu;0,S)}$, i.e.,
\begin{equation}
\label{BSDE_penalized_sup-->0ZU_bis}
\E\int_0^S |Z_s^T-Z_s|^2 ds + \E\int_0^S\int_A |U_s^T(a')-U_s(a')|^2\vartheta(da')ds \ \overset{T\rightarrow\infty}{\longrightarrow} \ 0.
\end{equation}
Now, fix $S\in[0,T]$ and consider the BSDE satisfied by $(Y^T,Z^T,U^T)$ on $[0,S]$:
\begin{align*}
Y_t^T \ &= \ Y_S^T - \beta \int_t^S Y_s^T\,ds+ n  \int_t^S \int_A (U_s^T(a'))_+\vartheta(da')\,ds + \int_t^S f(X_s^{x,a},I_s^a,\beta Y_s^T)\,ds \notag \\
&\quad \ - \int_t^SZ_s^T\,dW_s - \int_t^S\int_AU_s^T(a')\,\tilde\mu (ds,da'), \qquad 0\leq t\leq S.
\end{align*}
From \eqref{BSDE_penalized_sup-->0} and \eqref{BSDE_penalized_sup-->0ZU_bis}, we can pass to the limit in the above BSDE by letting $T\rightarrow\infty$, keeping $S$ fixed. Then, we deduce that $(Y,Z,U)$ solves the penalized BSDE \eqref{penalizedbetabsde} on $[0,S]$. Since  $S$ is arbitrary, it follows that $(Y,Z,U)$ solves equation \eqref{penalizedbetabsde} on $[0,\infty)$.
\ep

\subsubsection{Proof of Lemma \ref{L:vbeta,n_Bdd_Lip}}

The linear growth of $v^{\beta,n}$ follows from \eqref{E:v^beta,n} and the estimate on $Y^{x,a,\beta,n}$ of Proposition \ref{P:ExistUniqPen}. Concerning the identification $Y_t^{x,a,\beta,n}=v^{\beta,n}(X_t^{x,a},I_t^a)$, it is a consequence, as usual, of the flow property $(X_T^{x,a},I_T^a)=(X_{T-t}^{X_t^{x,a},I_t^a},I_{T-t}^{I_t^a})$ $\P$-a.s., for any $0\leq t\leq T<\infty$, and from the uniqueness for the penalized BSDE. Finally, regarding the uniform Lipschitz condition \eqref{v^beta,n_Lipschitz} of $v^{\beta,n}$ with respect to $x$, consider $x,x'\in\R^d$ and set $\Delta Y_t=Y_t^{x,a,\beta,n}-Y_t^{x',a,\beta,n}$, $\Delta Z_t=Z_t^{x,a,\beta,n}-Z_t^{x',a,\beta,n}$, $\Delta U_t(a')=U_t^{x,a,\beta,n}(a')-U_t^{x',a,\beta,n}(a')$, $t\geq0$, $a'\in A$. Let $T\in(0,\infty)$, then from estimate \eqref{EstimateDeltaYZU} there exists $\nu\in\Vc_n$:
\begin{align}
\label{penalizedBSDE_x-x'2}
|\bar Y_0|^2 \ &\leq \ \E^\nu\big[e^{-2\beta T}|\bar Y_T|^2\big] \notag \\
&\quad \ + 2\int_0^T e^{-2\beta s}\E^\nu\big[\bar Y_s\big(f(X_s^{x,a},I_s^a,\beta Y_s^{x,a,\beta,n}) - f(X_s^{x',a},I_s^a,\beta Y_s^{x,a,\beta,n})\big)\big]ds \notag \\
&\leq \ \E^\nu\big[e^{-2\beta T}|\bar Y_T|^2\big] \\
&\quad \ + 2\int_0^T e^{-2\beta s}\sqrt{\E^\nu\big[|\bar Y_s|^2\big]}\sqrt{\E^\nu\big[\big|f(X_s^{x,a},I_s^a,\beta Y_s^{x,a,\beta,n}) - f(X_s^{x',a},I_s^a,\beta Y_s^{x,a,\beta,n})\big|^2\big]}ds. \notag
\end{align}
Set $g(s)=e^{-2\beta s}\E^\nu[|Y_s|^2]$ and $h(s)=2e^{-\beta s}\sqrt{\E^\nu[|f(X_s^{x,a},I_s^a,\beta Y_s^{x,a,\beta,n})-f(X_s^{x',a},I_s^a,\beta Y_s^{x,a,\beta,n})|^2]}$, for any $s\in[0,T]$, and proceed as in \eqref{estimateY_Proof0}. Then, we conclude that
\begin{align*}
|\bar Y_0|^2 \ &\leq \ \sqrt{\E^\nu\big[e^{-2\beta T}|\bar Y_T|^2\big]} + \int_0^T e^{-\beta s}\sqrt{\E^\nu\big[\big|f(X_s^{x,a},I_s^a,\beta Y_s^{x,a,\beta,n}) - f(X_s^{x',a},I_s^a,\beta Y_s^{x,a,\beta,n})\big|^2\big]} ds \\
&\leq \ \sqrt{\E^\nu\big[e^{-2\beta T}|\bar Y_T|^2\big]} + L_2\int_0^T e^{-\beta s}\sqrt{\E^\nu\big[\big|X_s^{x,a} - X_s^{x',a}\big|^2\big]} ds.
\end{align*}
Therefore, recalling that $|\bar Y_T|\leq2C_{b,\sigma,f}(1+|X_T^{x,a}|)/\beta$ and using estimate \eqref{estimate_x^2}, we obtain
\[
\E^\nu\big[e^{-2\beta T}|\bar Y_T|^2\big] \ \overset{T\rightarrow\infty}{\longrightarrow} \ 0.
\]
By Lemma \ref{estimX}, we find
\[
\E^\nu\big[|X_s^{x,a}-X_s^{x',a}|^2\big] \ \le \ e^{-2\gamma s}|x-x'|^2, \qquad \forall\,s\geq0.
\]
Then, we deduce that
\[
|\bar Y_0| \ \leq \ L_2\int_0^\infty e^{-(\beta+\gamma) s}|x-x'|ds \ = \ \frac{L_2}{\beta+\gamma}|x-x'| \ \leq \ \frac{L_2}{\gamma}|x-x'|,
\]
which implies \eqref{v^beta,n_Lipschitz}.
\ep

\subsubsection{Proof of Proposition \ref{P:ViscPropv^beta,n}}

\emph{Continuity.} Fix $(\beta,n)\in(0,\infty)\times\N$. Let $x,x'\in\R^d$ and $a,a'\in\R^q$. Set $\Delta Y_t=Y_t^{x,a,\beta,n}-Y_t^{x',a',\beta,n}$, $\Delta Z_t=Z_t^{x,a,\beta,n}-Z_t^{x',a',\beta,n}$, $\Delta U_t=U_t^{x,a,\beta,n}-U_t^{x',a',\beta,n}$. Then, from estimate \eqref{EstimateDeltaYZU} we find, for any $0\leq r\leq T$,
\begin{align*}
|\Delta Y_r|^2 \ &\leq \ e^{-2\beta T}\E^\nu\big[|\Delta Y_T|^2\big] \\
&\quad \ + 2\E^\nu\bigg[\int_r^T e^{-2\beta s}\Delta Y_s\big(f(X_s^{x,a},I_s^a,\beta Y_s^{x,a,\beta,n}) - f(X_s^{x',a'},I_s^{a'},\beta Y_s^{x,a,\beta,n})\big)ds \bigg] \\
&\leq \ e^{-2\beta T}\E^\nu\big[|\Delta Y_T|^2\big] + \E^\nu\int_r^T e^{-2\beta s}|\Delta Y_s|^2ds \notag \\
&\quad \ + \E^\nu\int_r^Te^{-2\beta s}\big|f(X_s^{x,a},I_s^a,\beta Y_s^{x,a,\beta,n}) - f(X_s^{x',a'},I_s^{a'},\beta Y_s^{x,a,\beta,n})\big|^2ds.
\end{align*}
From Gronwall's lemma applied to the map $s\mapsto\E^\nu[|\Delta Y_s|^2]$ we obtain
\begin{align*}
|\Delta Y_0|^2 \ &= \ e^{\frac{1-e^{-2\beta T}}{2\beta}}\bigg(e^{-2\beta T}\E^\nu\big[|\Delta Y_T|^2\big] \notag \\
&\quad \ + \E^\nu\int_0^Te^{-2\beta s}\big|f(X_s^{x,a},I_s^a,\beta Y_s^{x,a,\beta,n}) - f(X_s^{x',a'},I_s^{a'},\beta Y_s^{x,a,\beta,n})\big|^2ds\bigg).
\end{align*}
From the Lipschitz property of $f$ in {\bf (H2)}, we find
\begin{align}
\label{continuity_v^beta,n4}
|\Delta Y_0|^2 \ &= \ e^{\frac{1-e^{-2\beta T}}{2\beta}}\bigg(e^{-2\beta T}\E^\nu\big[|\Delta Y_T|^2\big] \notag \\
&\quad \ + 2L_2^2\int_0^Te^{-2\beta s}\big\{\E^\nu\big[|X_s^{x,a} - X_s^{x',a'}|^2\big] + \E^\nu\big[|I_s^a - I_s^{a'}|^2\big]\big\}ds\bigg).
\end{align}
Now, for any $\eps>0$, applying It\^o's formula to $e^{(2\gamma-\eps-\eps L_1^2)t}|X_t^{x,a}-X_t^{x',a'}|^2$ and proceeding as in the proof of estimate \eqref{ergoX}, we obtain
\[
\E^\nu\big[|X_t^{x,a}-X_t^{x',a'}|^2\big] \ \leq \ |x-x'|^2 + \Big(1+\frac{2}{\eps}\Big)L_1^2\int_0^t e^{(2\gamma-\eps-\eps L_1^2)(s-t)} \E^\nu\big[|I_s^a-I_s^{a'}|^2\big]ds.
\]
Denote by $T_1$ the first jump time of the marked point process $(T_n,\alpha_n)_{n\geq1}$ associated to the Poisson random measure $\mu$. Notice that the two processes $I^a$ and $I^{a'}$ coincide after $T_1$, while we have $I_s^a=a$ and $I_s^{a'}=a'$ before $T_1$. In other words, $|I_s^a-I_s^{a'}|=|a-a'|1_{\{s\leq T_1\}}\leq|a-a'|$. Therefore
\begin{align*}
\E^\nu\big[|X_t^{x,a}-X_t^{x',a'}|^2\big] \ &\leq \ |x-x'|^2 + \Big(1+\frac{2}{\eps}\Big)L_1^2|a-a'|^2\int_0^t e^{(2\gamma-\eps-\eps L_1^2)(s-t)} ds \\
&\leq \ |x-x'|^2 + \Big(1+\frac{2}{\eps}\Big)L_1^2|a-a'|^2\int_0^\infty e^{-(2\gamma-\eps-\eps L_1^2)s} ds \\
&= \ |x-x'|^2 + \Big(1+\frac{2}{\eps}\Big)L_1^2|a-a'|^2\frac{1}{2\gamma-\eps-\eps L_1^2}.
\end{align*}
Therefore, \eqref{continuity_v^beta,n4} becomes
\begin{align}
\label{continuity_v^beta,n4_bis}
|\Delta Y_0|^2 \ &\leq \ e^{\frac{1-e^{-2\beta T}}{2\beta}}\bigg(e^{-2\beta T}\E^\nu\big[|\Delta Y_T|^2\big] + C_0\big(|x-x'|^2 + |a-a'|^2\big)\int_0^Te^{-2\beta s}ds\bigg) \notag \\
&\leq \ e^{\frac{1-e^{-2\beta T}}{2\beta}}\Big(e^{-2\beta T}\E^\nu\big[|\Delta Y_T|^2\big] + C_0\big(|x-x'|^2 + |a-a'|^2\big)\frac{1}{2\beta}\Big),
\end{align}
for some positive constant $C_0$, possibly depending on $L_1,L_2,\eps$, but independent of $T$. Since $|\Delta Y_T|\leq2C_{b,\sigma,f}(1+|X_T^{x,a}|)/\beta$, using estimate \eqref{estimate_x^2} we see that $\E[e^{-2\beta T}|\Delta Y_T|^2]\rightarrow0$ as $T\rightarrow\infty$. Therefore, letting $T\rightarrow\infty$ in \eqref{continuity_v^beta,n4_bis}, it follows that $|\Delta Y_0|^2\rightarrow0$ as $(x',a')\rightarrow(x,a)$. Since $\Delta Y_0=v^{\beta,n}(x,a)-v^{\beta,n}(x',a')$, then $v^{\beta,n}$ is continuous in both arguments.\\
\emph{Viscosity property.} We shall now prove the viscosity supersolution property of $v^{\beta,n}$. A similar argument would show that $v^{\beta,n}$ it is a viscosity subsolution to equation \reff{ellipticpdepenalizedbeta}. Let $(\bar x,\bar a)\in\R^d\times\R^q$ and $\varphi\in C^2(\R^d\times\R^q)$ such that
\begin{align}
\label{E:min_v^beta,n-phi}
0 \ = \ (v^{\beta,n}-\varphi)(\bar x,\bar a) \ = \ \min_{\R^d\times\R^q}(v^{\beta,n}-\varphi).
\end{align}
Let us proceed by contradiction, assuming that
\begin{align*}
\beta \,\varphi(\bar x,\bar a) - \Lc^{\bar a}\varphi(\bar x,\bar a) - \Mc^{\bar a}\varphi(\bar x,\bar a) - f(\bar x,\bar a,\beta v^{\beta,n}(\bar x,\bar a)) & \\
- n\int_A[\varphi(\bar x,a') -\varphi(\bar x,\bar a)]_+\,\vartheta(da')& \ =: \ -2\eps \ < \ 0.
\end{align*}
Using the continuity of $b$, $\sigma$, $f$, and $v^{\beta,n}$, we find $\delta>0$ such that
\begin{align}
\label{E:eps}
\beta \,\varphi(x,a) - \Lc^a\varphi(x,a) - \Mc^a\varphi(x,a) - f(x,a,\beta v^{\beta,n}(x,a)) & \notag \\
- n\int_A[\varphi(x,a')-\varphi(x,a)]_+\,\vartheta(da')& \ \leq \ -\eps,
\end{align}
for any $(x,a)\in\R^d\times\R^q$, with $|x-\bar x|,|a-\bar a|<\delta$. Define
\[
\tau \ := \ \inf\big\{t\geq0 \colon |X_t^{\bar x,\bar a}-\bar x| > \delta,\, |I_t^{\bar a}-\bar a| > \delta\big\}\wedge\delta
\]
Since $(X^{\bar x,\bar a},I^{\bar a})$ is c\`{a}dl\`{a}g, it is in particular right-continuous at time $0$. Therefore, $\tau>0$, $\P$-almost surely. Then, an application of It\^o's formula to $e^{-\beta t}\varphi(X_t^{\bar x,\bar a},I_t^{\bar a})$ between $0$ and $\tau$, using also \eqref{E:eps}, yields
\begin{align}
\label{E:Inequality_phi}
&e^{-\beta\tau}\varphi(X_\tau^{\bar x,\bar a},I_\tau^{\bar a}) \notag \\
&\geq \ \varphi(\bar x,\bar a) +\eps\frac{1-e^{-\beta\tau}}{\beta} - n \int_0^\tau\int_A e^{-\beta t}\big(\tilde U_t^{\bar x,\bar a,\beta,n}(a')\big)_+ \vartheta(da')dt \notag \\
&\quad \ - \int_0^\tau e^{-\beta t}f(X_t^{\bar x,\bar a},I_t^{\bar a},\beta v^{\beta,n}(X_t^{\bar x,\bar a},I_t^{\bar a}))dt + \int_0^\tau e^{-\beta t}(D_x \varphi(X_t^{\bar x,\bar a},I_t^{\bar a}))\trans\sigma(X_t^{\bar x,\bar a},I_t^{\bar a})dW_t \notag \\
&\quad \ + \int_0^\tau \int_A e^{-\beta t}\tilde U_t^{\bar x,\bar a,\beta,n}(a') \tilde\mu(dt,da'),
\end{align}
where $\tilde U_t^{\bar x,\bar a,\beta,n}(a') = \varphi(X_t^{\bar x,\bar a},a') - \varphi(X_t^{\bar x,\bar a},I_{t^-}^{\bar a})$. On the other hand, applying It\^o's formula to $e^{-\beta t}Y^{\bar x,\bar a,\beta,n}_t$ from $0$ to $\tau$, and using the identification $Y^{\bar x,\bar a,\beta,n}_t$ $=$ $v^{\beta,n}(X^{\bar x,\bar a}_t, I^{\bar a}_t)$, we find
\begin{align}
\label{E:BSDEv^beta,n}
v^{\beta,n}(\bar x,\bar a) \ &= \ e^{-\beta\tau}v^{\beta,n}(X_\tau^{\bar x,\bar a},I_\tau^{\bar a}) + n\int_0^\tau \int_A e^{-\beta t}\big(U_t^{\bar x,\bar a,\beta,n}(a')\big)_+ \vartheta(da') dt \notag \\
&\quad \ + \int_0^\tau e^{-\beta t}f(X_t^{\bar x,\bar a},I_t^{\bar a},\beta v^{\beta,n}(X_t^{\bar x,\bar a},I_t^{\bar a})) dr - \int_0^\tau e^{-\beta t}Z_t^{\bar x,\bar a,\beta,n} dW_t \notag \\
&\quad \ - \int_0^\tau\int_A e^{-\beta t}U_t^{\bar x,\bar a,\beta,n}(a') \tilde\mu(dt,da').
\end{align}
Plugging identity \eqref{E:BSDEv^beta,n} into inequality \eqref{E:Inequality_phi}, we obtain
\begin{align}
\label{E:phi-v^beta,n}
&e^{-\beta\tau}\varphi(X_\tau^{\bar x,\bar a},I_\tau^{\bar a}) - e^{-\beta\tau}v^{\beta,n}(X_\tau^{\bar x,\bar a},I_\tau^{\bar a}) \notag \\
&\geq \ \varphi(\bar x,\bar a) - v^{\beta,n}(\bar x,\bar a) + \eps\frac{1-e^{-\beta\tau}}{\beta} \notag \\
&\quad \ - n\int_0^\tau \int_A e^{-\beta t}\big[\big(\tilde U_t^{\bar x,\bar a,\beta,n}(a')\big)_+ - \big(U_t^{\bar x,\bar a,\beta,n}(a')\big)_+\big] \vartheta(da') dt \notag \\
&\quad \ + \int_0^\tau e^{-\beta t}\big(\sigma\trans(X_t^{\bar x,\bar a},I_t^{\bar a})D_x \varphi(X_t^{\bar x,\bar a},I_t^{\bar a}) - Z_t^{\bar x,\bar a,\beta,n}\big)dW_t \notag \\
&\quad \ + \int_0^\tau\int_A e^{-\beta t}\big(\tilde U_t^{\bar x,\bar a,\beta,n}(a') - U_t^{\bar x,\bar a,\beta,n}(a')\big) \tilde\mu(dt,da').
\end{align}
Define the $[1,n+1]$-valued $\Pc\otimes\Bc(A)$-measurable map $\nu$ as follows
\[
\nu_t(a') \ = \ 1 + n\frac{(\tilde U_t^{\bar x,\bar a,\beta,n}(a'))_+ - (U_t^{\bar x,\bar a,\beta,n}(a'))_+}{\tilde U_t^{\bar x,\bar a,\beta,n}(a')-U_t^{\bar x,\bar a,\beta,n}(a')}1_{\{\tilde U_t^{\bar x,\bar a,\beta,n}(a')-U_t^{\bar x,\bar a,\beta,n}(a')\neq0\}}.
\]
Then, we have $\nu\in\Vc_n$. Let us introduce the probability measure $\P^\nu$ equivalent to $\P$ on $(\Omega,\Fc_T)$, with $T\geq\tau$ (e.g., $T=\delta$), with Radon-Nikodym density given by \eqref{Pnu}. Then, taking the expectation $\E^\nu$ with respect to $\P^\nu$ in \eqref{E:phi-v^beta,n}, (recalling that $\varphi(\bar x,\bar a)=v^{\beta,n}(\bar x,\bar a)$)
\begin{equation}
\label{E:phi-v^beta,n2}
\E^\nu\big[e^{-\beta\tau}\big(\varphi(X_\tau^{\bar x,\bar a},I_\tau^{\bar a}) - v^{\beta,n}(X_\tau^{\bar x,\bar a},I_\tau^{\bar a})\big)\big] \ \geq \ \eps\E^\nu\bigg[\frac{1-e^{-\beta\tau}}{\beta}\bigg].
\end{equation}
Since $\tau>0$, $\P$-a.s., we see that the right-hand side of \eqref{E:phi-v^beta,n2} is strictly positive. On the other hand, from \eqref{E:min_v^beta,n-phi} it follows that the left-hand side of \eqref{E:phi-v^beta,n2} is nonpositive, therefore we get a contradiction.
\ep

\subsubsection{Proof of Proposition \ref{T:Exist}}

Firstly, we prove point (i). To this end, consider $Y^{x,a,\beta,n}$ and $Y^{x,a,\beta,n+1}$. It is useful to fix $T>0$ and to look at the penalized BSDE  \eqref{penalizedbetabsde} on $[0,\infty)$ solved by $Y^{x,a,\beta,n}$ (resp. $Y^{x,a,\beta,n+1}$) as a BSDE on $[0,T]$ with terminal condition $Y_T^{x,a,\beta,n}$ (resp. $Y_T^{x,a,\beta,n+1}$) and generator function $f$. Then, proceeding as in the proof of the comparison Theorem 2.5 of \cite{royer06} for BSDEs with jumps on $[0,T]$, we can find a probability measure $\P^\nu$ equivalent to $\P$ on $(\Omega,\Fc_T)$, such that
\begin{equation}
\label{E:Yn-Yn+1}
Y_t^{x,a,\beta,n}-Y_t^{x,a,\beta,n+1} \ \leq \ \E^\nu\big[e^{-\beta (T-t)}\big(Y_T^{x,a,\beta,n}-Y_T^{x,a,\beta,n+1}\big)\big|\Fc_t\big],
\end{equation}
$\P$-a.s., for all $t\in[0,T]$. Now, from estimate \eqref{estimate_x^2_conditional} and since 
$|Y_T^{x,a,\beta,n}|$, $|Y_T^{x,a,\beta,n+1}|\leq C_{b,\sigma,f}(1+|X_T^{x,a}|)/\beta$, letting $T\rightarrow\infty$ in \eqref{E:Yn-Yn+1}, we obtain $Y_t^{x,a,\beta,n} \leq Y_t^{x,a,\beta,n+1}$, $\P$-a.s., for all $t\geq0$. This shows that the sequence $(Y^{x,a,\beta,n})_n$ is monotone increasing. Since it is bounded by $C_{b,\sigma,f}(1+|X_T^{x,a}|)/\beta$, it converges increasingly to some adapted process $Y^{x,a,\beta}$ satisfying $|Y_t^{x,a,\beta}|\leq C_{b,\sigma,f}(1+|X_t^{x,a}|)/\beta$, for all $t\geq0$.\\
Now, fix again $T>0$ and consider the BSDE with nonpositive jumps \eqref{BSDE}-\eqref{JumpCon}, with $(X,I)=(X^{x,a},I^a)$, on $[0,T]$ with terminal condition $Y_T^{x,a,\beta}$ and generator function $f$. From Theorem 2.1 in \cite{khapha12} we know that there exists a unique minimal solution
\[
(\tilde Y^{x,a,\beta,T},\tilde Z^{x,a,\beta,T},\tilde U^{x,a,\beta,T},\tilde K^{x,a,\beta,T})\in{\bf S^2(0,T)}\times{\bf L^2(W;0,T)}\times{\bf L^2(\tilde\mu;0,T)}\times{\bf K^2(0,T)}
\]
to this BSDE. Moreover, $\tilde Y^{x,a,\beta,T}$ is the increasing limit\footnote{Notice that in Theorem 2.1 in \cite{khapha12}, the terminal condition in the penalized BSDE does not depend on $n$; while, in our case, the penalized BSDE associated with $Y^{x,a,\beta,n}$ has the terminal condition $Y_T^{x,a,\beta,n}$. However, since $Y_T^{x,a,\beta,n}$ converges increasingly $\P$-a.s. to $Y_T^{x,a,\beta}$ as $n\rightarrow\infty$, the results of Theorem 2.1 in \cite{khapha12} are still valid.} of $(Y^{x,a,\beta,n})_n$, so that $\tilde Y_t^{x,a,\beta,T}=Y_t^{x,a,\beta}$, $\P$-a.s., for all $t\in[0,T]$. We also know that $(\tilde Z^{x,a,\beta,T},\tilde U^{x,a,\beta,T})$ is the strong $($resp. weak$)$ limit of $(Z^{x,a,\beta,n},U^{x,a,\beta,n})_n$  in ${\bf L^p(W;0,T)}\times{\bf L^p(\tilde \mu;0,T)}$, with $p\in[1,2)$, $($resp. in ${\bf L^2(W;0,T)}\times{\bf L^2(\tilde \mu;0,T)}$$)$. This implies that $(\tilde Z^{x,a,\beta,T},\tilde U^{x,a,\beta,T})=(\tilde Z^{x,a,\beta,T'}_{|[0,T]},\tilde U^{x,a,\beta,T'}_{|[0,T]})$, for all $T'\geq T$. Then, we define $(Z^{x,a,\beta},U^{x,a,\beta})\in\mathbf{L_{\textup{loc}}^2(W)}\times\mathbf{L_{\textup{loc}}^2(\tilde\mu)}$ as
\begin{equation}
\label{Z,U=tildeZ,tildeU}
(Z^{x,a,\beta}_{|[0,T]},U^{x,a,\beta}_{|[0,T]}) \ = \ (\tilde Z^{x,a,\beta,T},\tilde U^{x,a,\beta,T}), \qquad \forall\,T>0.
\end{equation}
This proves point (ii). Concerning point (iii), from Theorem 2.1 in \cite{khapha12} we have that $\tilde K_t^{x,a,\beta,T}$ is the weak limit of $(K_t^{x,a,\beta,n})_n$ in ${\bf L^2}(\Fc_t)$, for any $0 \leq t \leq T$. It follows that $\tilde K_t^{x,a,\beta,T}=\tilde K_t^{x,a,\beta,T'}$, $\P$-a.s., for all $t\in[0,T]$ and for any $T'\geq T$. Therefore, we define $K^{x,a,\beta}\in{\bf K_{\text{loc}}^2}$ as follows: $K_t^{x,a,\beta}=\tilde K_t^{x,a,\beta,T}$, for all $t\in[0,T]$ and $T>0$. We see that the quadruple $(Y^{x,a,\beta},Z^{x,a,\beta},U^{x,a,\beta},K^{x,a,\beta})$ solves the backward equation \eqref{BSDE} on $[0,\infty)$.\\
Regarding the jump constraint \eqref{JumpCon}, from Theorem 2.1 in \cite{khapha12} we know that
\[
\tilde U_t^{x,a,\beta,T} \ \leq \ 0, \qquad dt\otimes d\P\otimes\vartheta(da)\text{-a.e.}
\]
Then, from the definition \eqref{Z,U=tildeZ,tildeU} of $U^{x,a,\beta}$ we see that \eqref{JumpCon} holds. It remains to prove the minimality condition. Let $(\bar Y^{x,a,\beta},\bar Z^{x,a,\beta},\bar U^{x,a,\beta},\bar K^{x,a,\beta})\in\mathbf{S_{\textup{loc}}^2}\times\mathbf{L_{\textup{loc}}^2(W)}\times\mathbf{L_{\textup{loc}}^2(\tilde\mu)}\times\mathbf{K_{\textup{loc}}^2}$ be another solution to \eqref{BSDE}-\eqref{JumpCon}, with $|\bar Y_t^{x,a,\beta}|\leq C(1+|X_t^{x,a}|)$, for all $t\geq0$ and for some positive constant $C$ (possibly depending on $x$, $a$, and $\beta$). Then, for any $T>0$, $(\bar Y^{x,a,\beta}_{|[0,T]},\bar Z^{x,a,\beta}_{|[0,T]},\bar U^{x,a,\beta}_{|[0,T]},\bar K^{x,a,\beta}_{|[0,T]})$ solves the BSDE \eqref{BSDE}-\eqref{JumpCon} on $[0,T]$. As before, proceeding as in the proof of the comparison Theorem 2.5 of \cite{royer06} for BSDEs with jumps on $[0,T]$, we can find a probability measure $\P^\nu$ equivalent to $\P$ on $(\Omega,\Fc_T)$, such that
\begin{equation}
\label{E:Yn-barY}
Y_t^{x,a,\beta,n}-\bar Y_t^{x,a,\beta} \ \leq \ \E^\nu\big[e^{-\beta (T-t)}\big(Y_T^{x,a,\beta,n}-\bar Y_T^{x,a,\beta}\big)\big|\Fc_t\big],
\end{equation}
$\P$-a.s., for all $t\in[0,T]$. From $|Y_T^{x,a,\beta,n}| \leq C_{b,\sigma,f}(1+|X_T^{x,a}|)/\beta$, $|\bar Y_T^{x,a,\beta}|\leq C(1+|X_T^{x,a}|)$ and estimate \eqref{estimate_x^2_conditional}, letting $T\rightarrow\infty$ in \eqref{E:Yn-barY} we obtain $Y_t^{x,a,\beta,n} \leq \bar Y_t^{x,a,\beta}$, $\P$-a.s., for all $t\geq0$. Then, sending $n\rightarrow\infty$, we find $Y_t^{x,a,\beta} \leq \bar Y_t^{x,a,\beta}$, $\P$-a.s., for all $t\geq0$, which proves the minimality of $(Y^{x,a,\beta},Z^{x,a,\beta},U^{x,a,\beta},K^{x,a,\beta})$ and concludes the proof.
\ep



\newpage

\small

\begin{thebibliography}{}


\bibitem{arendtbatty95} Arendt W. and C. J. Batty (1995): ``A complex Tauberian theorem and mean ergodic semigroups'', in \emph{Semigroup Forum}, \textbf{50}, 351-366, Springer New York.

\bibitem{arilio98} Arisawa M. and P.-L. Lions (1998): ``On ergodic stochastic control",  {\it Comm. in Partial Differential Equations}, {\bf 23}, 2187-2217. 


\bibitem{barsou01} Barles G. and P. Souganidis (2001): ``Space-time periodic solutions and long time behavior of solutions to quasi-linear parabolic equations", {\it SIAM J. Math. Anal.}, 
{\bf 32}, 1311-1323. 




\bibitem{benfre92} Bensoussan A. and J. Frehse (1992): ``On Bellman equations of ergodic control in $\R^n$", {\it J. Reine Angew. Math.}, {\bf 429}, 125-160. 

\bibitem{brihu98} Briand P. and Y. Hu (1998): ``Stability of BSDEs with random terminal time and homogenization of semilinear elliptic PDEs", {\it Journal of Functional Analysis}, {\bf 155}, 455-494. 




\bibitem{crandallishiilions92} Crandall M. G., Ishii H. and P.-L. Lions (1992): ``User's guide to viscosity solutions of second order partial differential equations'', {\it Bulletin of the American Mathematical Society}, \textbf{27}, 1-67.

\bibitem{daprato_zabczyk92} Da Prato G. and J. Zabczyk (1992): \emph{Stochastic Equations in Infinite Dimensions}, Encyclopedia of Mathematics and its Applications \textbf{44}, Cambridge University Press.

\bibitem{daprato_zabczyk96} Da Prato G. and J. Zabczyk (1996): \emph{Ergodicity for infinite dimensional systems}, volume \textbf{229}, Cambridge University Press.

\bibitem{darling95} Darling R. W. R. (1995): ``Constructing gamma-martingales with prescribed limit, using backwards SDE'', {\it The Annals of Probability}, \textbf{23}, 1234-1261.






\bibitem{fleshe99} Fleming W. and S.J. Sheu (1999): ``Optimal long term growth rate of expected utility of wealth", {\it Annals of Applied Probability}, {\bf 9}, 871-903. 


\bibitem{fuhrman_hu_tess09} Fuhrman M., Hu Y., and G. Tessitore (2009): ``Ergodic BSDEs and optimal ergodic control in Banach spaces'', {\it SIAM J. Control Optim.}, \textbf{48}, 1542-1566.


\bibitem{fujishlor06} Fujita Y., Ishii H. and P. Loreti (2006): ``Asymptotic solutions of Hamilton-Jacobi equations in euclidian $n$ space", {\it Indiana Univ. Math. J.}, {\bf 55}, 1671-1700. 


\bibitem{hatnagshe10} Hata H., Nagai H. and S.J. Sheu (2010): ``Asymptotics of the probability minimizing a down side risk",  {\it Ann. Appl. Probab.}, {\bf  20},  52-89.


\bibitem{humadric14} Hu Y., Madec P.-Y. and A. Richou (2014): ``Large time behavior of mild solutions of Hamilton-Jacobi-Bellman equations in infinite dimension by a probabilistic approach", preprint arXiv:1406.5993.


\bibitem{ichihara12} Ichihara N. (2012): ``Large time asymptotic problems for   optimal stochastic  control with superlinear cost'', {\it Stochastic Processes and their Applications}, 
\textbf{122}, 1248-1275.



\bibitem{ichshe13} Ichihara N.  and S.-J. Sheu (2013):  ``Large time behavior of solutions of Hamilton-Jacobi-Bellman equations with quadratic nonlinearity in gradients", 
{\it SIAM J. Math. Anal.}, {\bf  45}(1) 279-306. 


\bibitem{ishii89} Ishii H. (1989): ``On uniqueness and existence of viscosity solutions of fully nonlinear second order elliptic PDE's'', {\it Communications on Pure and Applied Mathematics}, \textbf{42}, 15-45.



\bibitem{jacod_shiryaev03} Jacod J. and A. N. Shiryaev (2003): \emph{Limit Theorems for Stochastic Processes}, second edition, Springer, Berlin.




\bibitem{khapha12} Kharroubi I. and H. Pham (2012): ``Feynman-Kac representation for Hamilton-Jacobi-Bellman IPDEs'', to appear in {\it The Annals of Probability}. 










\bibitem{protter05} Protter P. E. (2005): {\it Stochastic integration and differential equations}, volume 21 of Stochastic Modelling and Applied Probability. Springer-Verlag, Berlin. Second edition. Version 2.1, Corrected third printing.


\bibitem{nag03} Nagai H. (2003): ``Optimal strategies for risk-sensitive portfolio optimization problems with general factors", {\it SIAM J. Cont. and Optim.}, {\bf 41}, 1779-1800. 
 

\bibitem{darpardoux97} Darling R. W. R. and E. Pardoux (1997): ``Backwards SDE with random terminal time and applications to semilinear elliptic PDE'', {\it The Annals of Probability}, \textbf{25}, 1135-1159.


\bibitem{ric09}  Richou A. (2009): ``Ergodic BSDEs and related PDEs with Neumann boundary conditions", to appear in {\it Stochastic Processes and their Applications}. 


\bibitem{robxin13} Robertson S. and H. Xing (2013): ``Large time behavior of solutions to semi-linear equations with quadratic growth in the gradient", preprint. 


\bibitem{roy04} Royer M. (2004): ``BSDEs with a random terminal time driven by a monotone generator and their links with PDEs", {\it Stochastics and stochastic reports}, {\bf 76}, 281-307. 
 
 

\bibitem{royer06} Royer M. (2006): ``Backward stochastic differential equations with jumps and related non-linear expectations'', {\it Stochastic Processes and their Applications}, {\bf 116}, 1358-1376.


\bibitem{safonov88} Safonov M. V. (1988): ``On the classical solution of nonlinear elliptic equations of second order'', {\it Izvestiya Rossiiskoi Akademii Nauk. Seriya Matematicheskaya}, {\bf 52}, 1272-1287.


\bibitem{souzha06} Souplet P. and Q. Zhang (2006): ``Global solutions of inhomogeneous Hamilton-Jacobi equations", {\it  J. Anal. Math.}, {\bf 99}, 355-396. 





\bibitem{tanli94} Tang S. and X. Li (1994): ``Necessary conditions for optimal control of stochastic systems with jumps'', {\it SIAM J. Control and Optimization}, {\bf 32}, 1447-1475.


\end{thebibliography}
\end{document}